\documentclass[11pt]{article}
\usepackage[all]{xy}
\usepackage{graphicx}
\usepackage[dvipsnames]{xcolor}
\usepackage[utf8]{inputenc}
\usepackage[english]{babel}
\usepackage{amsthm}
\colorlet{Jonathan}{ForestGreen}
\colorlet{Alan}{RedViolet}
\colorlet{Christian}{RawSienna}
\usepackage{stackengine}
\usepackage{url}
\usepackage[normalem]{ulem}
\usepackage[colorlinks]{hyperref}
\usepackage[colorinlistoftodos]{todonotes}
\usepackage[margin=3.5cm]{geometry}
\usepackage{amssymb}
\usepackage{amsmath}
\usepackage{mathrsfs}
\usepackage{comment}
\usepackage{enumerate}
\usepackage{latexsym} 
\newtheorem{thm}{Theorem}[section]

\newtheorem{prop}[thm]{Proposition} 
\newtheorem{lemma}[thm]{Lemma}
\newtheorem{cor}[thm]{Corollary} 
\newtheorem{dfn}[thm]{Definition}

\theoremstyle{definition}
\newtheorem{rmk}[thm]{Remark}

\theoremstyle{definition}
\newtheorem{ex}{Example}[section]

\newcommand{\pf}{\noindent{\bf Proof~~}\ }

\newcommand{\reals}{{\mathbb R}}

\newcommand{\id}{{\text{id}}}

\newcommand{\Ker}{{\rm Ker\, }}
\newcommand{\im}{{\rm Im\, }}
\newcommand{\End}{{\rm End}}

\newcommand{\calq}{{\mathcal Q}}

\newcommand{\calv}{{\mathcal V}}

\newcommand{\fieldk}{{\mathbf k}}

\newcommand{\from}{\leftarrow}

\newcommand{\arrows}{\,\lower1pt\hbox{$\longrightarrow$}\hskip-.24in\raise2pt
             \hbox{$\longrightarrow$}\,}

\newcommand{\xc}[1]{}

\usepackage{tikz}
\usepackage{tikz-cd}
\usetikzlibrary{math}
\usetikzlibrary{calc}

\newtheorem{claim}[thm]{Claim}

\newcommand{\rad}{\rm Rad}

\title{(Co)isotropic triples and poset representations}
\author{Christian Herrmann$^*$ \and Jonathan Lorand$^\circ$ \and Alan Weinstein$^\dagger$ 
\\
\ }

\date{%
\footnotesize{
    $^*$Fachbereich Mathematik, Technische Universit\"at Darmstadt, 64289 Darmstadt, Germany, {\tt{herrmann@mathematik.tu-darmstadt.de}}
     \\%
    $^\circ$Institute of Mathematics, University of Zurich, CH-8057 Zurich, Switzerland, {\tt{jonathan.lorand@math.uzh.ch}}
    \\%
    $^\dagger
    $Department of Mathematics, University of California Berkeley, CA 94720-3840 USA. {\tt{alanw@math.berkeley.edu}}
    \\[2ex]%
    }
}

\begin{document}
\maketitle

\begin{abstract}
We study triples of coisotropic or isotropic subspaces in symplectic vector spaces; in particular, we classify indecomposable structures of this kind. The classification depends on the ground field, which we only assume to be perfect and not of characteristic 2. Our work uses the theory of representations of partially ordered sets with (order reversing) involution; for (co)isotropic triples, the relevant poset is ``$2 + 2 + 2$'' consisting of three independent ordered pairs, with the involution exchanging the members of each pair.

A key feature of the classification is that any indecomposable (co)isotropic triple is either ``split'' or ``non-split''. The latter is the case when the poset representation underlying an indecomposable (co)isotropic triple is itself indecomposable. Otherwise, in the ``split'' case, the underlying representation is decomposable and necessarily the direct sum of a dual pair of indecomposable poset representations; the (co)isotropic triple is a ``symplectification''. 

If one ignores the involution, representations of the poset ``$2 + 2 + 2$'' are essentially the same as representations of a certain quiver of extended Dynkin type $\tilde E_6$, and our work relies on the known classification of indecomposable representations of such quivers. In particular from these results it follows that the poset representations underlying (co)isotropic triples come in two types: there are families of indecomposables which depend on a parameter taking a continuum of values (``continuous-type''), and there are those indecomposables which are characterized only by the dimensions of the spaces involved (``discrete-type''). 
This pattern is reflected in the classification of (co)isotropic triples; in particular, there are families of triples which depend on a parameter. Also, indecomposable triples exist in every even dimension. 



In the course of the paper we develop the framework of ``symplectic poset representations'', which can be applied to a range of problems of symplectic linear algebra. The classification of linear hamiltonian vector fields, up to conjugation, is an example; we briefly explain the connection between these and (co)isotropic triples. The framework lends itself equally well to studying poset representations on spaces carrying a non-degenerate \emph{symmetric} bilinear form; we mainly keep our focus, however, on the symplectic side.

\end{abstract}

\tableofcontents

\

\section{Introduction}
In a previous paper \cite{lo-we:coisotropic}, the second two authors have given a complete classification of pairs of coisotropic subspaces in Poisson vector spaces and, equivalently by duality, pairs of isotropic subspaces in presymplectic vector spaces.
Each such pair is uniquely (up to isomorphism and ordering) decomposable as a direct sum of multiples of ten indecomposable pairs, for which there are simple normal forms in ambient spaces of dimension $1$, $2$, or $3$.   
These decomposition problems are special cases of the problem of classifying triples of coisotropic or isotropic subspaces in symplectic vector spaces, with an extra condition relating the third subspace to the first two.

In this paper, we will deal with  \emph{all} (co)isotropic triples in symplectic vector spaces.   The decomposition into indecomposables is still possible 
with summands 
 essentially unique, but there are many more indecomposables.   In dimensions up to 4, there are still only finitely many  isomorphism classes
of indecomposables, while in higher dimensions the moduli space of such classes includes parametrized families as well as single points which may or may not be in the closure of such families.

The classification of pairs in \cite{lo-we:coisotropic} was done by elementary arguments in linear algebra, but the results there (as well as those in 
\cite{lo-we:relations} on the classification of (co)isotropic relations) suggested links with the representation theory of quivers, particularly of those associated with (extended) Dynkin diagrams and
the closely related representation theory of partially ordered sets (posets).
(See e.g. Gabriel and Roiter \cite{garo}.)
We rely on these   to carry out the classification of triples. In fact, we largely reduce our problem to that of studying representations which are maps from a certain 6-element poset with involution to the poset of subspaces of a symplectic vector space with the involution given by symplectic orthogonality.   The classification of these representations, without the involution, is essentially that of certain representations of a quiver associated to the extended Dynkin diagram $\tilde{E_6}$, which is a tree consisting of a central vertex attached to three ``branches" containing two vertices each.  The quiver in question is obtained from $\tilde {E_6}$ by orienting all of its edges toward the central vertex. Here are two depictions of this quiver; the first is the most common, while we will use the latter, to emphasize the poset structure: 

$$
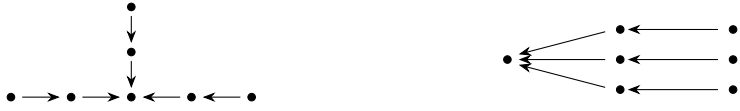
\begin{tikzpicture}[>=Stealth, xscale=1, yscale=1]
\tikzmath{
coordinate \c; 
\c0 = (-4,0);
\c1 = (-4.8,0);
\c2 = (-4,0.6);
\c3 = (-3.2,0);
\c4 = (-5.6,0);
\c5 = (-4,1.2);
\c6 = (-2.4,0);
}

\draw [fill] (\c0) circle [radius=0.05];

\draw [fill] (\c1) circle [radius=0.05];
\draw [fill] (\c2) circle [radius=0.05];
\draw [fill] (\c3) circle [radius=0.05];

\draw [fill] (\c4) circle [radius=0.05];
\draw [fill] (\c5) circle [radius=0.05];
\draw [fill] (\c6) circle [radius=0.05];

\tikzmath{
coordinate \tip, \tail;
\tip1 = (\c0) + (-0.15,0);
\tail1 = (\c1) + (0.15,0);
\tip2 = (\c0) + (0,0.12);
\tail2 = (\c2) + (0,-0.12);
\tip3 = (\c0) + (0.15,0);
\tail3 = (\c3) + (-0.15,0);
\tip4 = (\c1) + (-0.15,0);
\tail4 = (\c4) + (0.15,0);
\tip5 = (\c2) + (0,0.12);
\tail5 = (\c5) + (0,-0.12);
\tip6 = (\c3) + (0.15,0);
\tail6 = (\c6) + (-0.15,0);
}

\draw [<-] (\tip1) -- (\tail1);
\draw [<-] (\tip2) -- (\tail2);
\draw [<-] (\tip3) -- (\tail3);

\draw [<-] (\tip4) -- (\tail4);
\draw [<-] (\tip5) -- (\tail5);
\draw [<-] (\tip6) -- (\tail6);

\tikzmath{
coordinate \c; 
\c0 = (2,0.5) + (-1,0);
\c1 = (3.5,.9) + (-1,0);
\c2 = (3.5,0.5) + (-1,0);
\c3 = (3.5,0.1) + (-1,0);
\c4 = (5,0.9) + (-1,0);
\c5 = (5,0.5) + (-1,0);
\c6 = (5,0.1) + (-1,0);
}

\draw [fill] (\c0) circle [radius=0.05];

\draw [fill] (\c1) circle [radius=0.05];
\draw [fill] (\c2) circle [radius=0.05];
\draw [fill] (\c3) circle [radius=0.05];

\draw [fill] (\c4) circle [radius=0.05];
\draw [fill] (\c5) circle [radius=0.05];
\draw [fill] (\c6) circle [radius=0.05];

\tikzmath{
coordinate \tip, \tail;
\tip1 = (\c0) + (0.15,0.07);
\tail1 = (\c1) + (-0.2,-0.03);
\tip2 = (\c0) + (0.15,0.0);
\tail2 = (\c2) + (-0.2,0);
\tip3 = (\c0) + (0.15,-0.07);
\tail3 = (\c3) + (-0.2,0.03);
\tip4 = (\c1) + (0.1,0.0);
\tail4 = (\c4) + (-0.2,0);
\tip5 = (\c2) + (0.1,0.0);
\tail5 = (\c5) + (-0.2,0);
\tip6 = (\c3) + (0.1,0.0);
\tail6 = (\c6) + (-0.2,0);
}

\draw [<-] (\tip1) -- (\tail1);
\draw [<-] (\tip2) -- (\tail2);
\draw [<-] (\tip3) -- (\tail3);

\draw [<-] (\tip4) -- (\tail4);
\draw [<-] (\tip5) -- (\tail5);
\draw [<-] (\tip6) -- (\tail6);

\end{tikzpicture}
$$
For our results, we rely on the classification of representations of extended Dynkin quivers given in \cite{dl-ri-indecomposable} and \cite{do-fr:representations}, for the case of $\tilde{E_6}$. Representations of this particular quiver have also been studied in quite some detail by Stekolchshik, see \text{e.g.} \cite{raf} and \cite{raf2}. The study of poset representations in spaces equipped with an (anti-)symmetric inner product was first developed, to our  knowledge, by Scharlau and collaborators; see \cite{sc-subspaces} for a concise and enlightening overview. 

For (co)isotropic triples, the connection with the $\tilde {E_6}$-type quiver described above is this: the central vertex corresponds to the ambient symplectic vector space, while each branch corresponds to an isotropic subspace and (adjacent to the central vertex) the coisotropic orthogonal in which it is contained. 

The associated six-element poset consists of the vertices in the branches, with partial ordering given by the arrows connecting them, and (order-reversing)  involution given by exchanging the elements of each pair.  We will use one of the standard notations, $\bf 2+\bf 2+\bf 2$, for this poset.

Since the operation of ``taking the symplectic orthogonal'' induces a one-to-one correspondence between isotropic subspaces and coisotropic ones, in the following we focus on and refer simply to isotropic triples.   The concomitant results for coisotropic triples are implicit. 

A crucial part of our analysis is a result due to Quebbemann, Scharlau, and Schulte \cite{qu-sc-sc:quadratic} and  Sergeichuk \cite{se:classification systems}, which, in our setting, says that every indecomposable isotropic triple
 in dimension $2n$ 
is either already indecomposable as a linear representation of the poset $\bf 2+\bf 2+\bf 2$ or is obtained from an indecomposable linear representation of the same poset in dimension $n$ by a ``doubling construction" known as  hyperbolization
(\cite[p. 267]{qu-sc-sc:quadratic}), and which in our context we will call {\bf symplectification}.
 (It is closely connected to the cotangent bundle construction in symplectic geometry, though the latter always produces lagrangian subspaces.)  This dichotomy reduces our problem to deciding which such indecomposable linear representations actually come from indecomposable isotropic triples, and finding the nonisomorphic isotropic triples which may give rise to the same indecomposable linear representation.\footnote{The simplest triple illustrating this possibility consists of three distinct lines in a symplectic plane.  These are lagrangian, so each one corresponds to both the isotropic and coisotropic subspaces in a nested pair. If we forget the symplectic structure, there is no further invariant, but in the case of real coefficients, there is a symplectic invariant given by the cyclic order of the three lines with respect to the symplectic orientation.    This is sometimes called the Maslov index or Kashiwara-Vergne index of the lagrangian triple.}

\section*{Acknowledgements}

J. L. would like to thank Alberto Cattaneo, Wilhelm Karlsson, Marcel Wild, and Thomas Willwacher; he acknowledges support from the NCCR SwissMAP and SNF Grant No. 200020 172498/1, both funded by the Swiss National Science Foundation. A.W. acknowledges partial support from the UC Berkeley Committee on Research. 

\section{Resum\'e}
 
In this section we give a summary of our results. We hope that placing this summary here, rather than as a final section, will give the reader an initial rough overview and a place to refer back to while reading the paper. Some of the terminology and notation used below is only defined later in the course of the paper; we include here, however, a short list of our most essential terms and notions:

\begin{itemize}
\item A \textbf{linear poset representation} (or just ``poset representation'') of a poset $P$ in a vector space $V$ (always assumed finite-dimensional) is an order-preserving map from $P$ to the poset of subspaces of $V$. Poset representations will usually be denoted by the letter $\psi$. Representations of the poset ${\bf 2} +{\bf 2} +{\bf 2}$ are called sextuples. 
\item If a poset $P$ is equipped with an order-reversing involution, then to each representation $\psi$ of $P$ in $V$, there is a \textbf{dual} representation $\psi^*$ of $P$ in $V^*$. (For the definition, see (\ref{dual representation def}) in Section \ref{section duals comp forms}.) In this paper, we always assume the poset ${\bf 2} +{\bf 2} +{\bf 2}$ to be equipped with the order-reversing involution that exchanges the two elements of each of the three ordered pairs. 
\item An \textbf{isotropic triple} $\varphi$ is an ordered triple of isotropic subspaces of a symplectic vector space. The three isotropics, together with their corresponding symplectic orthogonals, form a sextuple of subspaces which define a linear representation of the poset ${\bf 2} +{\bf 2} +{\bf 2}$, the underlying
linear representation $\hat{\varphi}$ of the triple. The underlying poset representation $\hat \varphi$ of an isotropic triple is necessarily isomorphic to its dual; \text{i.e} $\hat \varphi$ is ``self-dual''. 
\item Given a sextuple $\psi$ in $V$, a \textbf{compatible symplectic form} is a symplectic form on $V$ with respect to which $\psi$ becomes the underlying poset representation of an isotropic triple.  (See Section \ref{section duals comp forms}.)
\item Isotropic triples are examples of \textbf{symplectic poset representations} (see Definition \ref{dfn symp rep}). There are distinct notions of decomposition (and of ``indecomposable'') for  linear poset representations and for symplectic poset representations, respectively (the latter are \emph{orthogonal} decompositions.) In particular, an indecomposable isotropic triple may be decomposable as a linear poset representation.   
\item Any representation $\psi$ of a poset $P$ has an associated \textbf{dimension vector}: this is the function which assigns to each element $x \in P$ the dimension of the associated subspace $\psi(x)$. We usually write the dimension vector not as a function, but as a tuple (\text{i.e.} a ``vector''). 
\item Throughout we assume that the ground field $\fieldk$ is perfect and does not have characteristic $2$. The condition of ``perfectness'' is not a strong one: perfect fields include fields which are algebraically closed, fields of characteristic zero, and finite fields. 
\end{itemize}

\

\noindent
\textbf{Resum\'e}:

\begin{enumerate}
\item 
Symplectically indecomposable isotropic triples $\varphi$  come in two kinds:
\begin{itemize}
\item Split: the underlying sextuple $\hat \varphi$ of the isotropic triple is decomposable as a linear poset representation. By Lemma \ref{sergeichuk}, $\varphi$ is isomorphic to the symplectification of some indecomposable poset representation  $\psi$ of ${\bf 2} +{\bf 2} +{\bf 2}$  which is not self-dual; in particular $\hat{\varphi} \cong \psi \oplus \psi^*$.
\item Non-split: $\varphi$ is such that $\hat \varphi$ indecomposable and self-dual 
\end{itemize}
Our main work is the identification of the duals of indecomposable sextuples and the classification of the indecomposable non-split isotropic triples. Once this is done, by Lemma \ref{sergeichuk}, the classification of the split-type indecomposable isotropic triples is essentially automatic.  Throughout, it is understood that when we speak of classification of indecomposables, we mean the classification of \emph{isomorphism classes} of indecomposables. 
\item In order to classify the non-split indecomposable isotropic triples, we first identify which sextuples $\psi$ are self-dual. We have the following subcases (we use labels (a) through (e), which are also referred to further below):
\begin{itemize}
\item Discrete-type: in these cases, $\psi$ is based on an indecomposable nilpotent endomorphism. Up to isomorphism, $\psi$ is uniquely determined by its dimension vector $\bf d$. There are the following types, with $k \in \mathbb{Z}_{\geq 0}$:
\begin{itemize}
\item[(a)] $A(3k+1,0)$ with ${\bf d}=(3k+1;2k+1,k,2k+1,k,2k+1,k),$
\item[(b)] $A(3k+2,0)$ with  ${\bf d}=(3k+2;2k+1,k+1,2k+1,k+1,2k+1,k+1).$
\end{itemize}
\item Continuous-type: in these cases, $\psi$ has a dimension vector of the form $(3k,2k,k,2k,k,2k,k)$, with $k \in \mathbb{Z}_{> 0}$, and is based on an indecomposable endomorphism $\eta$ (this is a generalization of the $\Delta(k,\lambda)$ of Donovan-Freislich  \cite{do-fr:representations}, and we also use normal forms based on the homogeneous representations of Dlab-Ringel  \cite{dl-ri-indecomposable}). In our encoding, the underlying endomorphism $\eta$ of a self-dual $\psi$ is necessarily such that $\eta$ is similar to $({\sf id} -\eta)^*$. The following types of endomorphism $\eta$ (up to similarity) account for all the indecomposable self-dual continuous-type sextuples:

\noindent
If $\eta$ has an eigenvalue $\lambda$ in the base field:
\begin{itemize}
\item[(c)] $\lambda =\frac{1}{2}.$
\end{itemize}
If $\eta$ has has no eigenvalue in the base field:
\begin{itemize}
\item[(d)] Over the reals: the complex eigenvalues of $\eta$ are $\frac{1}{2}\pm b_\eta \sqrt{-1}$ for unique real  $b_\eta >0$.
\item[(e)] Over perfect fields in general: $\eta=\frac{1}{2} {\sf id} +\zeta$ where $\zeta \neq 0$ is similar to $-\zeta^*$. The characteristic (= minimal) polynomial of $\zeta$ is of the form $r(x)^m$ for an irreducible polynomial $r$ which is of the form $r(x)=p(x^2)$ for some polynomial $p$. 
\end{itemize}
\end{itemize}
\item With the indecomposable self-dual sextuples classified, we then determine which of these admit compatible symplectic forms (and we give such forms explicitly via coordinate matrices). We find that compatible symplectic forms exists as follows: 
\begin{itemize}
\item[(a)] for $A(3k+1,0)$ if and only if $k$ is odd.
\item[(b)] for $A(3k+2,0)$ if and only if $k$ is even. 
\item[(c)] for $\eta$ having an eigenvalue in $\fieldk$: if and only if $k$ is even.
 \item[(d-e)] for $\eta$ having no eigenvalue in $\fieldk$: for all $k$. 
\end{itemize}
\item Following the question of existence of compatible symplectic forms for indecomposable self-dual sextuples, we then address the question of uniqueness. We find that compatible symplectic forms for a sextuple $\psi$ are unique up to isomorphism (i.e. up to an isometry which is an automorphism of $\psi$) and, in
\begin{itemize}
\item[(a-c),] multiplication by a scalar. There are no compatible symmetric forms.
\item[(d),] multiplication by a scalar. There are also compatible symmetric forms.
\item[(e),] multiplication by a `scalar' $Z \in F^+_{H}$.
Here $F =  \fieldk[x]/q(x)$ is considered as a subring of
$\fieldk^{3k\times 3k}$ and $F^+_{H}$ consists of all $Z \in F$ such that $Z^t H=HZ$, where $H$ is the coordinate matrix of a particular compatible form. There are also compatible symmetric forms.
\end{itemize}
\item The following provides a complete list of isomorphism types of isotropic triples in the non-split case. Given a non-split isotropic triple $\varphi$ with  symplectic  form $\omega$ (and associated coordinate matrix $H$) there is an automorphism of $\hat{\varphi}$ which is
\begin{itemize}
\item[(a-d)]
an isometry from $\omega$ to $c\omega$,
for some $0\neq c \in \fieldk$, if and only if $c$ is a square in $\fieldk$. Thus, compatible symplectic forms for a given sextuple $\psi = \hat \varphi$ are parametrized by the square class group $\fieldk^{\times}/(\fieldk^{\times})^2$. 
\item[(e)] an isometry from $H$ to $HZ$, for some $0 \neq Z \in F$, 
if and only if $Z=X^2-Y^2$ where $X\in F^+_{H}$ and
$Y\in F^-_{H}$. Here, $F^{\pm }_{H}=\{Z\mid Z^tH=\pm HZ\}$. 
\end{itemize}
\item For a fixed indecomposable sextuple $\psi$, the set (if non-empty) of all compatible symplectic forms\footnote{Here we mean \emph{all} compatible forms, \text{i.e.} not only up to isometry.} is given via linear expressions involving
\begin{itemize}
\item[(a)] $k$ parameters,
\item[(b)] $k+1$ parameters,
\item[(c-d)] $\frac{k}{2}+1$ parameters.
\end{itemize}
\end{enumerate}

\noindent \textbf{Proof of the Resum\'e}
\begin{enumerate}
\item By Lemma \ref{sergeichuk}.
\item 
For the cases (a-b), in Proposition \ref{DFdiscrete} it is stated which discrete-type indecomposable sextuples are self-dual, and in Theorem \ref{dis} it is shown which ones admit compatible symplectic forms. 

For continuous-type indecomposable sextuples, it is established in Section \ref{duals of continuous sextuples} that self-duality only occurs for framed sextuples, and from Proposition \ref{duals of endo-sextuples} it follows that a framed sextuple $\mathcal{S}_\eta$ is self-dual if and only if $\eta$ is similar to $\id - \eta^*$. Furthermore, those $\mathcal{S}_\eta$ which underly non-split isotropic triples are identified:
\begin{itemize}
\item[(c)] in Theorem \ref{eig} for the case when the underlying endomorphism has an eigenvalue in the ground field. 
\item[(d)] in Theorem \ref{real} for the case when $\fieldk = \mathbb{R}$ and the underlying endomorphism has no eigenvalue in $\mathbb{R}$.
\item[(e)] in Theorem \ref{hom} for the general case when $\fieldk$ is perfect and the underlying endomorphism has no eigenvalue in $\fieldk$. See also Corrollary \ref{paramet underlying}. 
\end{itemize}

\item For (a-b), see Theorem \ref{dis}; for (c), see Theorem \ref{eig}; (d-e), see Theorem \ref{hom}.
\item For (a-b), see Theorem \ref{unidis}; for (c), see Theorem \ref{eig}; for (d), see Theorem \ref{real}; for (e), see Theorem \ref{uniqueness general fields}. 
The uniqueness statements in the theorems above for (a-d) make use of Lemma
\ref{uniqueness up to a scalar}, and in the case (e),  Theorem \ref{uniqueness general fields} makes use of  Lemma \ref{uniqu}. 

\item As in the previous point: for (a-b), see Theorem \ref{unidis}; for (c), see Theorem \ref{eig}; for (d), see Theorem \ref{real}; for (e), see Theorem \ref{uniqueness general fields}. 
\item For (a-b), see Remark \ref{rem}; for (c-d), see Theorem \ref{eig} and Theorem \ref{real}.
\end{enumerate}
\qed \\

\begin{cor}
For an indecomposable representation $\psi$  of the poset $P={\bf 2} +{\bf 2} +{\bf 2}$ in a vector space $V$ the following are equivalent
\begin{enumerate}
\item There is a symplectic representation $\varphi$ such that $\hat{\varphi}=\psi$.
\item $\psi$ is self-dual and $\dim V$ is even.
\end{enumerate}
\end{cor}

\begin{rmk}
The following topics are deferred to possible future work.
\begin{itemize}
\item
For those representations which do {\em not} admit compatible symplectic
structures  (since they are not self-dual or admit only symmetric structure), a detailed description of the isotropic triples resulting from their symplectifications.
\item Analysis of how the classification of isotropic triples changes when the ground field is changed.
\item A description of the isometry groups of indecomposable isotropic triples.
\item A discussion of the classification of isotropic triples
within the category-theoretic framework of
Quebbemann, Scharlau, Schulte
\cite{qu-sc-sc:quadratic}
and Sergeichuk \cite{se:classification systems}. 
\item Analysis of the relation between poset representations which are sextuples and those which involve four subspaces. 
\item The study of how a general (not necessarily indecomposable) isotropic triple decomposes into indecomposable summands; in particular the question of how unique such a decomposition is. 
\item The question of defining invariants for isometry types of isotropic triples
(in particular in relation to the perfect elements established
by Stekolchshik \cite{raf}).
\item An explanation in further detail of the relation between isotropic triples and linear hamiltonian vector fields. 
\end{itemize}
\end{rmk}

\section{General theory}

Here we provide a review of some aspects of symplectic linear algebra and the representation theory of partially ordered sets.  We have tried to do this in such as way as to be accessible to readers who might be familiar with only one of these two areas. In particular,
we give an elementary presentation of the relevant results
from the general theory of Quebbemann, Scharlau, and Schulte
\cite{qu-sc-sc:quadratic} and Sergeichuk \cite{se:classification systems}.  

We work with a fixed ground field, which we assume to be perfect (\text{e.g.} of characteristic zero, finite, or algebraically closed) and not of characteristic 2; otherwise we usually leave the field unspecified. All vector spaces are assumed to be finite-dimensional.

\subsection{Basic notions}

A pair of subspaces $(A,B)$ in a vector space $V$ (without further structure) is completely determined up to isomorphism by four invariants: the dimensions of $V$, $A$, $B$, and $A\cap B$.    For a triple $(A,B,C)\subseteq V$ of subspaces,  one needs to know, in addition to the dimensions of $V$, $A$, $B$, $C$, all pairwise intersections, and the triple intersection, the dimension of one more space, such as $(A+B)\cap C$, giving a total of nine invariants.  (For instance, this ninth invariant is needed to distinguish the two arrangements of three distinct lines in 
3-space, which may or may not be coplanar.)  
 At this point, we have effectively introduced the lattice $\langle A,B,C \rangle$ generated by $A$, $B$, and $C$ as a sublattice of the lattice $\Sigma (V)$ of subspaces in $V$, \text{i.e.} the subspaces generated by $A$, $B$, and $C$ under iterations of the operations of sum and intersection.
The study of such structures dates back to Dedekind and Thrall, immanent in early results of representation theory and leading into abstract lattice theory. The present work is mainly oriented around posets, with $\Sigma (V)$ being a poset with respect to inclusion, and we will not deal with abstract lattice theory. Nevertheless, the lattice structure of $\Sigma (V)$ will play an important role, in particular in calculations. 

A non-degenerate antisymmetric bilinear, i.e. {\bf symplectic}, form $\omega$ on $V$ produces additional structure.  

First of all, there is a naturally associated linear map $\tilde{\omega}$ from $V$ to $V^*$ defined by $\tilde{\omega}(v)(w) = \omega(v,w)$.  Non-degeneracy means that $\tilde{\omega}$ is an isomorphism.

Next, there
is a natural order-reversing involution $A \mapsto A^\perp$ on $\Sigma (V)$, where $A^\perp$, the {\bf
symplectic orthogonal} of $A$, is the subspace $\{v\in V \mid \forall w\in A, \omega (v,w) = 0\}$.  This involution is related to another order-reversing operation, namely the one that maps a subspace $A \subseteq V$ to its annihilator $A^\circ = \{ \xi \in V^* \mid \xi (v) = 0 \ \forall v \in A \}$ in $V^*$. The following result is easy to verify; we state it as a lemma to refer to later.

\begin{lemma}\label{perpannihilator}
For any subspace $A$ in a symplectic space $(V,\omega)$, 
$\tilde{\omega}(A^\perp) \subseteq V^*$ is the annihilator $A^\circ$ of $A$.
\end{lemma}

Using the involution $\perp$, we may define several special types of subspaces.   There are the {\bf isotropic} subspaces, for which $A \subseteq A^\perp$ and the {\bf coisotropic} subspaces, for which $A^\perp \subseteq A$.   Subspaces which are both isotropic and coisotropic, i.e. fixed points of the orthogonality involution, are called {\bf lagrangian}.    Subspaces for which $A\cap A^\perp = \{0\}$ are
called {\bf symplectic}; the restriction of the symplectic form to such a subspace is non-degenerate, and hence again a symplectic form.   For such spaces, the order-reversing property implies that this condition is equivalent to 
$A + A^\perp = V$, so that we can also write the condition as $V = A \oplus A^\perp$.  Involutivity implies that $A^\perp$ is symplectic as well, giving a {\bf symplectic direct sum decomposition} of $V$. Such symplectic direct sums can also be built from the external direct sum of symplectic spaces by equipping the sum with the direct sum form, which is again symplectic. The aim of this paper is to study the decomposition of (co)isotropic triples with respect to symplectic direct sums (not just of two summands).

Given a symplectic direct sum decomposition $V=A_1\oplus A_2$,
we can define in a purely lattice-theoretic way the involution on $\Sigma (A_1)$ associated with the restricted symplectic structure.
Let $C_1$ be any subspace of $A_1$.  Then $C_1^\perp$ contains $A_2$, and the modular law\footnote{The modular law states that, for subspaces $A$, $B$, $C$ of a given vector space, if $B \subseteq A$, then $A \cap (B + C) = B+ A\cap C$.} implies that $C_1^\perp = (C_1^\perp\cap A_1 )\oplus A_2$. So we may define the operation $^{\perp_1}$ by setting $C_1^{\perp_1} := C_1^\perp \cap A_1$.  It is clearly order-reversing, and it is easy to check that it is involutive.   Of course, we can do the same thing in $A_2$.   We see from here that for subspaces $C_j \subseteq A_j$, the direct sum $C_1\oplus C_2$ is (co)isotropic or symplectic if and only if $C_1$ and $C_2$ are.
Similar statements to all of those above hold for orthogonal decompositions with any number of summands.

A {\bf linear representation}\footnote{We use the word ``linear'' here simply to distinguish this notion of representation from the notion of symplectic poset representation, which we define below.} of a partially ordered set, or poset, $P$ in a vector space $V$ is an order preserving map $\psi$ from $P$ to the poset $\Sigma (V)$. The {\bf dimension vector} of $\psi$ is the function which assigns to each $p\in P$ the dimension of $\psi(p)$.  We will also add $\dim V$ as a
first component to the dimension vector.
A {\bf morphism} $f:\psi \to \psi'$ 
between representations $\psi$ in $V$ and $\psi'$ in $V'$
of the same poset $P$ is a linear map $f:V \to V'$ such
that $f(\psi(p)) \subseteq \psi'(p)$
for all $p \in V$; it is an {\bf isomorphism} 
if $f$ is bijective and $f(\psi(p)) =\psi'(p)$
for all $p \in P$. 

The term \textbf{poset with involution} will denote a poset $P$ equipped with an order-reversing involution $ \perp : P \rightarrow P$. The relevance of poset representations and involutions arises when we observe that each isotropic subspace $A$ in a symplectic space $V$ is naturally associated with the subspace $A^\perp$ which contains it and hence with the pair $(A,A^\perp)$. Conversely, if $(A,B)$ is any pair of subspaces for which $A\subseteq B$ and $A^\perp = B$, then $A$ is isotropic and $B$ is its coisotropic orthogonal.    Thus we see that isotropic (or equivalently, coisotropic) subspaces in $V$ may be identified with involution-preserving representations into $\Sigma (V)$ of a  poset consisting of two elements in linear order, equipped with the order-reversing involution which exchanges the two elements.  We will denote this poset by $\bf 2$.  

\begin{dfn}\label{dfn symp rep}
Let $(P,  \perp)$ be a poset with involution, $(V, \omega)$ a symplectic vector space, and $(\Sigma (V), \ \perp)$ the poset of its subspaces equipped with the involution induced by $\omega$\footnote{We use the same symbol to denote two different involutions -- the one on $P$ and the one on $\Sigma (V)$.}. 
A {\bf symplectic representation} of $(P,\ \perp)$ in $V$ is an order-preserving map $\varphi : P \rightarrow \Sigma (V)$ such that 
$$
\varphi(p^\perp) = \varphi(p)^\perp \quad \forall p \in P.
$$
If $\varphi$ is such a representation, 
forgetting the involutions gives us simply a linear representation of $P$ in the vector space $V$, which we denote by $\hat{\varphi}$ and call the \textbf{underlying} linear representation. 

Given another symplectic representation $\varphi'$ in $(V', \omega')$ of the same poset with involution $(P,  \perp)$, 
an {\bf isomorphism} from $\varphi$ to $\varphi'$ is an isomorphism of linear poset representations $\hat \varphi \rightarrow \hat \varphi'$ 
which is also an isometry $(V, \omega) \rightarrow (V', \omega')$. 
\end{dfn}


Now, for any $k$, the $k$-tuples of isotropic-coisotropic pairs in $(V, \omega)$ are simply the symplectic representations in $(V, \omega)$ of the partially ordered set $k\cdot\bf 2  := \bf 2 +\bf 2 + \cdots + \bf 2$ ($k$ times) consisting of $k$ copies of $\bf 2$  which are independent in the sense that $a \leq b$ only when $a$ and $b$ belong to the same copy, with the involution interchanging the elements of each copy of $\bf 2$. 
In the present paper, we derive for $k=0,1,2,3$, up to direct decomposition, the classification of isotropic $k$-tuples from the classification of linear representations of the poset $k\cdot \bf 2$; these essentially correspond to the representations of a quiver associated with
the Dynkin diagrams $A_1$, $A_3$, and $A_5$, and the extended 
Dynkin diagram $\tilde{E}_6$, respectively. The first three cases are easy, with only finitely many indecomposables up to isomorphism, while the latter case is more involved -- in particular there are infinitely many indecomposables, some appearing in 1-parameter families.

Our general strategy for classifying isotropic triples will involve three basic operations related to symplectic representations.  First of all, given a symplectic representation $\varphi$ of $(P,\ \perp)$ in $(V,\omega)$, one may ignore the involution and symplectic structure to obtain a linear representation $\hat \varphi$ of $P$ in $V$.  Second, starting with a linear representation $\psi$ of $(P,\ \perp)$ in $V$, one can ask whether there exists a symplectic form on $V$ such that $\psi$ is in fact a symplectic representation in $(V, \omega)$; in this case $\omega$ is a {\bf compatible} symplectic form. Third,
one can build a symplectic representation out of each linear one by a ``doubling'' construction called {\bf symplectification}, which we  define in subsection \ref{duality} below. 

\begin{rmk}
These operations are close to ones in symplectic geometry.  The first one corresponds to forgetting the symplectic structure on a symplectic manifold, the second is analogous to asking whether a given manifold admits a symplectic structure, while the third is similar to the cotangent bundle construction.
\end{rmk}

\subsection{Decompositions of  linear representations}\label{decomp linear reps}

We fix a poset $P$. Given linear poset representations $\psi$ and $\psi'$, on $V$ and $V'$ respectively, their (external) \textbf{direct sum} is the poset representation on $V \oplus V'$ defined by 
$$(\psi \oplus \psi') (x) = \psi(x) \oplus \psi' (x) \qquad \forall x \in P.$$
A \textbf{subrepresentation} of a linear representation $\psi$ on $V$ is a representation $\psi'$ on a subspace $U \subseteq V$ such that $\psi'(x) \subseteq \psi(x) \ \forall x \in P$. Given subrepresentations $\psi'$ and $\psi''$ of $\psi$, on $U'$ and $U''$ respectively, we say they form an (internal) \textbf{direct sum decomposition} of $\psi$ if $V = U' \oplus U''$ and 
$$
\psi(x) = \psi'(x) \oplus \psi''(x) \qquad \forall x \in P.
$$
In this case, $\psi'(x) = \psi(x) \cap U'$ for all $x \in P$, and similarly for $\psi''$. 

We note that, given a subrepresentation $\psi'$ of $\psi$, there may not exist a subrepresentation $\psi''$ such that $\psi = \psi' \oplus \psi''$. 
For an example, consider $P= \{x_1, x_2, x_3 \}$ endowed with the empty partial order, \text{i.e.} there are no order relations between the elements. Let $\psi$ be a representation of $P$ on a two-dimensional space $V$ such that the subspaces $\psi(x_i)$ are three independent lines, and let $\psi'$ be the subrepresentation on $U' := \psi(x_1)$ such that $\psi'(x_1) = \psi(x_1)$ and $\psi'(x_2) =  \psi'(x_3) = 0$. Suppose there existed a subrepresentation $\psi''$ such that $\psi = \psi' \oplus \psi''$. Then $\psi''$ would be need to be defined on a subspace $U''$ such that $V = U' \oplus U''$; thus $U''$ would be a line in $V$. By the requirement that $\psi(x_i) = \psi'(x_i) \oplus \psi''(x_i)$, it follows that necessarily 
$$\psi (x_2) = 0 \oplus \psi''(x_2) \quad \text{ and } \quad \psi(x_3) = 0 \oplus \psi''(x_3).$$
But since the subspaces $\psi''(x_2)$ and $\psi''(x_3)$ must be contained in the line $U''$, this would imply that $\psi''(x_2) = \psi''(x_3) = U''$, a contradiction to the fact that $\psi''(x_2)$ and $\psi''(x_3)$ are independent. 

For any vector space $V$,  idempotents $\pi$ in the algebra $\End(V)$ of its endomorphisms correspond to direct sum decompositions $V = A \oplus B$ of $V$, where $A =\im \pi$ and $B = \Ker \pi = \im (1-\pi)$.
(The  zero and the identity endomorphisms are denoted by $0$ and $1$, respectively.) The same is true for the endomorphism algebra $\End(V,\psi)$ of a linear representation $\psi$ in $V$ of a fixed partially ordered set $P$. Recall that an endomorphism of a linear representation $\psi$ is a linear map $f:V\to V$ such that $f(\psi(p))\subseteq \psi (p)$ for all $p \in P$.  It is easy to check that $\End(V,\psi)$ is a unital subalgebra of $\End (V)$.  Since the ground field embeds into $\End(V,\psi)$ via its action on the algebra unit, we sometimes refer simply to the
\emph{ring} of endomorphisms of $\psi$.

The basic theory of poset representations parallels that of modules of finite composition length as presented, for example, in \cite[\S 1.4]{lambek}. We nevertheless give some outlines of proofs, for the convenience of the reader. 

We call a subring $E$ of some $\End(V)$ {\bf local} if each of its elements is either invertible or nilpotent. It follows that,  for any $f \in E$, either $f$ or $1-f$ is invertible; namely, if $f$ is not invertible, then  $f^n=0$ for some $n$ and $1-f$ has inverse $\sum_{i=0}^n  f^i$. 
More generally, if $g=\sum_{i=1}^m f_i$ is invertible then so is at least one of the $f_i$. 
The  nilpotent elements form an ideal $\rad E$, the {\bf radical} of $E$. (Namely, given $f$ and $g$, if one of them is nilpotent and if $fg$ were to be invertible, then  $f$ would be surjective and $g$ injective, so actually both would be invertible. And if both $f$ and $g$ are nilpotent then $f+g$ cannot be invertible.)  It follows that $E/\rad E$ is a division ring and $\rad E$ the unique maximal ideal. 

The following result is a version of Fitting's Lemma.
\begin{lemma}\label{local} 
For a linear poset representation $\psi$ in $V$ the following are equivalent
\begin{enumerate}
\item $\psi$ is indecomposable.
\item $\End(V,\psi)$ is local.
\item $\End(V,\psi)$  has only the trivial idempotents $0$ and $1$.
\end{enumerate}
\end{lemma}
\pf
Clearly, \text{2)} $\Rightarrow$ \text{3)} $\Rightarrow$ \text{1).} Now, assume $\psi$ indecomposable.
The rank of $f^j$ is a non-increasing and non-negative function of $j = 1,2,3,\cdots$, so it stabilizes after finitely many 
steps, say $d$ steps.  Then $ f^j(\im f^d) =  \im f^{d+j} = \im f^d$ for all negative $j$; Call this image space $I$.  It follows that the nondecreasing sequence of kernels of these powers also satisfies
$N:=\Ker f^{d+j} = \Ker f^d$; call this null space $N$.   If $v \in I\cap N$, then 
 $v = f^d(w)$ for some $w$, and $0=f^d v$.  Hence, $f^{2d}(w) = 0$, so $v=f^d(w)$ must be zero as well.   So $I\cap N = \{0\}$.
 By dimension counting, we have $V=I \oplus N$, where $I$ and $N$ are both $f$-invariant, $f|_I$ is invertible (since it is surjective), and $f_N$ is nilpotent.

 If $K$ is any $f$-invariant subspace, then $f^k (K)\subseteq K$ is the projection $\pi_I(K)$ of $K$ on $I$.   
 Then the other projection $\pi_N(K) = (1 - \pi_I)(K)$ is contained in $K$ as well.   This shows that any $f$-invariant subspace is the direct sum of its components in $I$ and $N$.   It follows that, if $f$ leaves invariant an indecomposable family of subspaces, either $I$ or $N$ must be zero, and $f$ is either nilpotent or invertible.  
\qed \\

 \begin{rmk}\label{rmk-fitting} 
Even without the indecomposability assumption, the decomposition into invertible and nilpotent parts is unique: for any such decomposition $f_I \oplus f_N$, $I$ and $N$ must be the image and kernel respectively of all sufficiently large powers of $f$.
 \end{rmk}


We now state the Krull-Remak-Schmidt theorem in the form that we will need.

\begin{thm}\label{krs}
Let $\psi=\psi_1 \oplus \psi_2 \oplus \cdots \oplus \psi_n$ and let $\psi'=\psi'_1 \oplus \psi'_2 \oplus \cdots \oplus \psi'_{n'}$ be direct sum decompositions  of isomorphic linear poset representations, e.g., of the same representation, into indecomposable summands. (Such decompositions always exist in finite dimensions.)   Then $n=n'$, and the two decompositions are the same up to isomorphism and permutation of the summands.
\end{thm}
\pf
We may assume $n'\geq n$, and we denote the underlying vector spaces as 
$V_1 \oplus \cdots \oplus V_n$ and $V'_1 \oplus \cdots \oplus V'_{n'}$.
First, consider $n=n'=2$ and an isomorphism $f:\psi \to \psi'$ 
such that $f(v_1,0)= (g(v_1),h(v_1))$ with
morphisms  $g,h:\psi_1 \to \psi'_1$ where $g$ is an isomorphism. 
We claim that $\psi_2$ and $\psi'_2$ are isomorphic, too. 
To prove this, we may assume $h=0$: replace $f$ by $f'$
where $f'(v_1,v_2)= (w_1, w_2- hg^{-1}(w_1))$ if
$f(v_1,v_2)=(w_1,w_2)$.    Then $\psi_2 \cong \psi/(\psi_1 \oplus 0)
\cong (f \psi)/(g\psi_1 \oplus 0) = \psi'/(\psi'_1 \oplus 0) \cong \psi'_2$.

In general, let $f:\psi \to \psi'$ be a
given isomorphism and consider the canonical embeddings
$\varepsilon_i, \varepsilon'_i$ and projections
 $\pi_i,\pi'_i$ given by the decompositions. Put
$g_i= \pi_1' \circ f \circ \varepsilon_i$ and
$h_i=\pi_i \circ f^{-1} \circ \varepsilon'_1$. 
Then $\sum_{i=1}^n g_i \circ h_i =1 \in \End(V'_1,\psi'_1)$. 
Since this ring is local, one of the summands is invertible,
say the first.  Then $ h_1 \circ g_1 \in \End(V_1,\psi_1)$   is not nilpotent,
whence invertible. Thus, $g_1:\psi_1 \to \psi'_1$ is
an isomorphism. Clearly, $f(v_1,0,\ldots ,0)
=(g_1(v_1),w_2, \ldots )$ and it follows, by the special case,
that 
 $\psi_2 \oplus \cdots \oplus \psi_n$
is isomorphic to $ \psi'_2 \oplus \cdots \oplus \psi'_{n'}$.
We repeat the argument until only $\psi_n$ is left on 
one side. Since $\psi_n$ is indecomposable
one has $n=n'$ whence $\psi_n\cong \psi'_n$.    

\qed \\

\subsection{Endomorphisms}

Poset representations of particular interest are ones which can be built from a vector space $U$ and a linear map $\eta : U \rightarrow U$. We will will refer to such a couple $(U, \eta)$ simply as an \textbf{endomorphism} when no conflation with other notions is to be feared. From such an endomorphism one can build a poset representation $(V; U_1, U_2, U_3, U_4)$, \text{i.e.} a representation of the poset $\textbf{1}$+$\textbf{1}$+$\textbf{1}$+$\textbf{1}$: define
\begin{equation}\label{quadruple}
\begin{array}{l}
V = U \times U \\
U_1 = U \times 0 \\
U_2 =  0 \times U \\
U_3 = \{(x,-x)\mid x \in V\} \\
U_4 = \{(x,-\eta(x))\mid x \in V\}.
\end{array}
\end{equation}
We will call this the \textbf{quadruple associated to $(U, \eta)$}. Note that $U_4$ and $U_3$ are the negative graphs, respectively, of $\eta$ and the identity map on $U$. It is straightforward to check that any endomorphism of the poset representation  (\ref{quadruple}) is of the form $f \times f$, where $f : U \rightarrow U$ is such that 
\begin{equation}\label{conjugate}
\eta \circ f = f \circ \eta. 
\end{equation}
The collection of linear maps $f : U \rightarrow U$ satisfying $(\ref{conjugate})$ form an algebra $\text{End}(U, \eta)$, the \textbf{endomorphism algebra of} $(U, \eta)$. More generally, a \textbf{morphism} of endomorphisms $(U, \eta) \rightarrow (U', \eta')$ is a linear map $f : U \rightarrow U'$ satisfying $\eta' \circ f = f \circ \eta$; such a map is an \textbf{isomorphism} if, additionally, $f$ is invertible as a linear map. 

\begin{lemma}
A poset representation $(V; U_1, U_2, U_3, U_4)$ is isomorphic to one of the form $(\ref{quadruple})$ for some $(U, \eta)$ if and only if
\begin{equation}\label{cond for endo}
U_i \oplus U_j = V \text{ for } i<j \leq 3 \quad \text{ and } \quad U_2 \oplus U_4 = V.
\end{equation}
\end{lemma}

\pf
It is easily seen that, given $(U,\eta)$, the associated quadruple satisfies (\ref{cond for endo}). For the converse, suppose $(V; U_1, U_2, U_3, U_4)$ is a quadruple satisfying (\ref{cond for endo}). First, note that the conditions (\ref{cond for endo}) imply that all of the subspaces $U_i$ have the same dimension. Choose any vector space $U$ having the same dimension as the $U_i$ and choose linear isomorphisms $\varphi_1 : U_1 \rightarrow U$ and $\varphi_2: U_2 \rightarrow U$. Next, note that since $U_3$ is independent of both $U_1$ and $U_2$, there exists an isomorphism $\varphi_3 : U_1 \rightarrow U_2$ such that $U_3 = \{ x + \varphi_3 x \mid x \in U_1 \} \subseteq U_1 \oplus U_2$. Finally, using this data we construct the following map 
$$
\varphi: V = U_1 \oplus U_2 \longrightarrow U \times U, \ x + y \longmapsto (\varphi_1 x, - \varphi_1 \varphi_3^{-1}  y).
$$
Clearly $\varphi$ is a linear isomorphism and maps $U_1$ to $U \times 0$ and $U_2$ to $0 \times U$. Furthermore,  
given $x + \varphi_3 x \in U_3$, its image under $\varphi$ is $(\varphi_1 x,-\varphi_1 x)$, as desired. Since $U_4$ is independent of $U_2$, it is the graph of a linear map $\varphi_4 : U_1 \rightarrow U_2$ (though this map may not be an isomorphism). Now given $x + \varphi_4 x \in U_4 \subseteq U_1 \oplus U_2$, we have 
$$
\varphi(x + \varphi_4 x) = (\varphi_1 x, - \varphi_1 \varphi_3^{-1} \varphi_4 x) = (\varphi_1 x, - \varphi_1 \varphi_3^{-1} \varphi_4 \varphi_1^{-1} \varphi_1 x).
$$
Setting $\eta =  \varphi_1 \varphi_3^{-1} \varphi_4 \varphi_1^{-1}$, we have a quadruple of the form $(\ref{quadruple})$ which is isomorphic to $(V; U_1, U_2, U_3, U_4)$. 
\qed \\

An endomorphism $(U, \eta)$ is \textbf{decomposable} if there exist non-zero subspaces $U_1, U_2 \subseteq U$ which are invariant under $\eta$ and form a decomposition $U = U_1 \oplus U_2$. As in the case of poset representations, such decompositions correspond to idempotents in the endomorphism algebra of $(U, \eta)$. For future reference we state:

\begin{lemma}\label{iso of algebras}
The endomorphism algebra of $(U, \eta)$ is isomorphic to the endomorphism algebra of the associated quadruple $(\ref{quadruple})$. In particular, $(U, \eta)$ is indecomposable if and only if its endomorphism algebra is local.
\end{lemma}

\pf
An isomorphism is given by mapping an endomorphism $f$ of $(U, \eta)$ to the endomorphism $f \times f$ of $(\ref{quadruple})$. Its inverse is ``restriction to $U_1$''. 

The characterization of indecomposabilty follows from Lemma \ref{local}.
\qed \\

Another point of view is that an endomorphism $(U, \eta)$ defines a $\fieldk[x]$-module $U_{\fieldk[\eta]}$, where the action of $x$ on $U_{\fieldk [\eta]}$ is defined by the action of $\eta$ on $U$, \text{i.e.} $(\sum_i \lambda_ix^i) \cdot u :=  (\sum_i \lambda_i \eta^i)(u)$ for $u \in U$, $\lambda_i \in \fieldk$. In this case, morphisms $f : (U, \eta) \rightarrow (U', \eta')$ are the same as module homomorphisms $f : U_{\fieldk[\eta]} \rightarrow U'_{\fieldk[\eta']}$ and decompositions of $(U, \eta)$ correspond to direct sum decompositions of $U_{\fieldk[\eta]}$.

\

\begin{rmk}\label{loc}
If $\eta$ is an indecomposable endomorphism of $U$
with an eigenvalue $\lambda$ in the base 
field $\fieldk$, then there is a basis such that $\eta$
is described by a single $\lambda$-Jordan block. With respect to such a basis, the members of 
the endomorphism algebra
$E = \text{End}(U, \eta)$ are given by the matrices $\sum_{i=0}^{d-1} a_i N^i$  
where $d=\dim U$ and $N$ is the nilpotent matrix
with $n_{i,i+1}=1$ for $i=1, \ldots ,d-1$, and $0$ otherwise.
In particular, $E$ is local
$($namely isomorphic to $\fieldk[x]/(x^d))$ and $E = \fieldk \id \oplus \rad E$.
\end{rmk} 
\pf
 Consider the endomorphism $\zeta=\eta -\lambda{\sf id}$ 
and apply the Fitting Lemma.  Since $\zeta$ is non-invertible,
there is $n$ with $\zeta^n=0$, that is $\eta$ admits
Jordan normal form with unique eigenvalue $\lambda$
and there is a single block $J=\lambda I+N$ only. 
Now, $E$ is given by the $A$ such that $AJ=JA$, which is equivalent to the condition $AN=NA$. Since any invariant subspace of $N$ is one of $A$ also, $A$ must be upper triangular and with only a single scalar on each upper diagonal.
In other words, from $AN=B=NA$ one has that $a_{i,j-1}=b_{ij}=a_{i+1,j}$
for $i<j$, and $0$ as entry otherwise.
\qed \\

If an indecomposable endomorphism $(U, \eta)$ does not have an eigenvalue in the ground field $\fieldk$, a generalization of Remark \ref{loc} holds; we recall this in Sections \ref{section genJNF}.

We note that, for $(U, \eta)$ indecomposable, the endomorphism algebra $\text{End}(U, \eta)$ is generated, as a unital $\fieldk$-subalgebra of $\text{End}(U)$,  by the single element $\eta$ (see Proposition \ref{endo alg}). In other words, it consists simply of ``polynomials in $\eta$". In particular, it is a commuatative algebra, and any subspace of $U$ is invariant under $\text{End}(U, \eta)$ if and only if it is invariant under $\eta$. Further details about these endomorphism algebras are discussed in Section \ref{endo alg details}.

\subsection{Decompositions in symplectic spaces}
If $V$ carries a symplectic form $\omega$, we define a transpose operation  $^t$ on its endomorphisms by the usual formula $\omega(f(x), y) = \omega(x,f^t (y))$.  Then the condition $\pi^t\pi = 0$ on an idempotent $\pi$ means that the image of $\pi$ is an isotropic subspace.   In fact, for any $x$ and $y$ in $V$, we have $\omega(\pi(x),\pi(y)) =  \omega(x,\pi^t \pi(y))$.  

Similarly, if $\pi \pi^t = 0$, then the image of $\pi^t$ is isotropic.   It follows that, if $\pi+\pi^t = 1$, then  the images of $\pi$ and $\pi^t$ give a decomposition of $V$ as the direct sum of two isotropic subspaces which must be lagrangian and hence in duality by the bilinear form.  We state this result in the form of a lemma for use below.

\begin{lemma}\label{lagrangian direct sum}
Decompositions of a symplectic vector space $V$ as a direct sum of (two) lagrangian subspaces are in 1-1 correspondence with idempotent endomorphisms $\pi:V\to V$ such that $\pi^t \pi = 0 = \pi\pi^t$ and $\pi+\pi^t = 1$. 
\end{lemma}

Similarly, for symplectic direct sum decompositions, we have the following characterization in terms of idempotents.

\begin{lemma}\label{symplectic direct sum}
Symplectic direct sum decompositions of $V$ into two subspaces are in 1-1 correspondence with idempotent endomorphisms $\pi:V\to V$ which are self-adjoint, i.e. $\pi^t=\pi$.   
\end{lemma}
\pf
It is a standard fact that $\Ker(\pi)=(\im \pi^t)^\perp$, and $\Ker \pi$ is also the image of $1-\pi$.  It follows that $\pi=\pi^t$ if and only if the images of $\pi$ and $1-\pi$ are orthogonal, which means that the corresponding direct sum decomposition is symplectic.
\qed \\

Now suppose that $\varphi$ is a symplectic representation  in $V$ of an involutive poset $P$.

\begin{lemma}\label{transpose endomorphism}
If $f$ is an endomorphism of $\varphi$, then $f^t$ is an endomorphism of $\varphi$ as well.
\end{lemma}
\pf We must show that, for each $p \in P$, $f^t(\varphi(p))\subseteq 
\varphi(p)$. This is equivalent to showing that, for $a \in \varphi(p)$ and $b \in \varphi(p)^\perp$, we have $\omega (f^t(a), b) = 0$, or, equivalently, $\omega(a,f(b))=0$.  Since $\varphi$ is symplectic, we have $b\in \varphi(p^\perp)$; since $f$ is an endomorphism of $\varphi$, $f(b)\in \varphi(p^\perp) = \varphi(p)^\perp$ as well.
\qed \\

\begin{lemma}\label{decomps of symp reps}
Symplectic direct sum decompositions of $\varphi$ into two representations are in $1$-$1$ correspondence with self-adjoint idempotents $\pi = \pi^t$ which are endomorphisms of the underlying linear representation, \text{i.e.} $\pi \in \text{End}(\hat \varphi)$. 
\end{lemma}

\pf 
A direct sum decomposition $V = V_1 \oplus V_2$ of $\varphi$ is nothing else than a decomposition of $\hat \varphi$ where the summands $V_1$ and $V_2$ are also mutual orthogonals. Such a decomposition corresponds to an idempotent $\pi \in \text{End}(\hat \varphi)$ such that $\pi^t = \pi$ (c.f. Lemma \ref{symplectic direct sum}). 
\qed \\

\begin{lemma}\label{symplectic fitting lemma}
If $\varphi$ is an indecomposable symplectic representation and if $f \in \text{End}(\hat \varphi)$ is an endomorphism of $\hat \varphi$ such that $f^t = \pm f$, then $f$ is either nilpotent or an isomorphism. 
\end{lemma}

\pf
From the proof of Lemma \ref{local} we know that there exists a non-negative integer $d$ such that $V =\text{Ker}f^d \oplus \text{Im}f^d$ is a decomposition of $\hat \varphi$. Since $(f^d)^t = (f^t)^d = \pm f^d$, this decomposition of $V$ is one into orthogonal summands; hence it is a symplectic direct sum decomposition of $\varphi$. Since $\varphi$ is assumed to be indecomposable as a symplectic representation, either $\text{Ker}f^d$ or $\text{Im}f^d$ must be zero. 
\qed \\

\subsection{Duals, Compatible forms}\label{section duals comp forms}

Let $\psi : P \rightarrow \Sigma (V)$ be any linear representation in $V$ of a poset $P$ equipped with an order-reversing involution $\perp$.
We can think of $\psi$ as a linear representation of $(P,\perp)$ which doesn't ``see" the involution. Define the \textbf{dual representation} of $\psi$ in $V^*$ by
\begin{equation}\label{dual representation def}
\psi^* : (P,\perp) \longrightarrow \Sigma (V^*), \quad \psi^*(p) = \psi(p^\perp)^\circ  \quad \forall p \in P.
\end{equation}
Note that this definition makes use of the involution on $P$, \text{i.e.} $\psi^*$ is the dual of $\psi$ \emph{with respect to the involution $ \ \perp$} on $P$. In particular, the definition only makes sense viewing $\psi$ as a linear representation of $(P, \perp)$ (rather than only of $P$). The combination of two order inversions - once due to the poset involution and once due to the annihilator operation
 - leads to $\psi^*$ being order preserving. 
 
By Lemma \ref{perpannihilator}, if $\varphi$ is a symplectic representation of $(P,  \perp)$ in $(V,\omega)$,  then $\tilde\omega$ is an isomorphism from $\hat\varphi$ to $\hat\varphi^*$.
This shows that a representation of $P$ can be compatible with a symplectic structure only if it is self-dual, i.e. isomorphic to its dual.  In particular, it has a dimension vector which is {\bf self-dual} in the sense that $\dim \varphi (p) + \dim \varphi (p^\perp) = \dim V$ for any $p\in P$.   

The results above lead to the natural question of determining the relation between isomorphisms of a linear representation to its dual, the self-duality of its dimension vector, and the existence (and uniqueness) of compatible symplectic structures. 

We will see later that many representations of interest to us are characterized up to isomorphism by their dimension vectors. Hence we record the following simple observation. 

\begin{lemma}\label{dualdimension}
If a representation $\psi$ is isomorphic to its dual, then 
it has a self-dual dimension vector. If a representation $\psi$ is characterized up to isomorphism by its dimension vector, and this dimension vector is self-dual, then $\psi$ is isomorphic to its dual.
\end{lemma}
\pf
The first statement has already been noted.   For the second, we observe  that, if the dimension vector of $\psi$ is self-dual, then 
$\psi$ and $\psi^*$ have the same dimension vector, and so by assumption they must be isomorphic.
\qed \\

In studying symplectic structures, it will sometimes be important to consider non-degenerate symmetric bilinear forms, as well. For these, there is an analogous notion of orthogonality and, therefore, an analogous notion of representation of a partially ordered set with involution. To capture both kinds of forms, we'll speak of $\varepsilon$-{\bf symmetric} bilinear forms, where $\varepsilon = 1$ for symmetric forms and $\varepsilon = -1$ for antisymmetric forms; similarly, an $\varepsilon$-{\bf symmetric representation} of a poset with involution is the generalization of the definition of symplectic representation to $\varepsilon$-{\bf symmetric} forms. We will often identify a non-degenerate bilinear form $B$ on $V$ with the linear isomorphism $B: V \rightarrow V^*, v \mapsto B(v, \cdot)$, setting $B^* : V \rightarrow V^*, v \mapsto B(\cdot, v)$. With this notation, $B$ is  $\varepsilon$-symmetric if and only if $B^* = \varepsilon B$\footnote{Note that if $B : V \rightarrow V^*$ is a non-degenerate bilinear form such that $B^* \neq \lambda B$ for some $\lambda \in \textbf{k}$, then in general, for a given subspace $A \subseteq V$, one no longer has ``the'' orthogonal ``$A^\perp$'', but rather one must consider the right- and left-orthogonal of $A$, which in general will not coincide.}. In referring to the parity $\varepsilon$ of $B$ we sometimes write $\varepsilon(B)$. 

For a fixed linear representation $\psi$ on $V$ of a poset $(P,\ \perp)$ with involution, one can ask how many different non-degenerate $\varepsilon$-symmetric forms $B$ exist (if any) which are compatible with $\psi$ in the sense that
$$
\psi(p^\perp) = \psi(p)^\perp \quad \forall p \in P,
$$
where the involution $\ \perp$ on $V$ is the one induced by $B$. 
A non-degenerate $\varepsilon$-symmetric form which is compatible with $\psi$ in this sense will be called a \textbf{compatible form}.

\begin{lemma}\label{compatible forms}
Let $\psi$ be a linear representation of $(P, \perp)$ in $V$, and $B: V \rightarrow V^*$ a non-degenerate $\varepsilon$-symmetric form. Then $B$ is a compatible form (for $\psi$) if and only if $B$ is an isomorphism $\psi \rightarrow \psi^*$. 
\end{lemma}

\pf
That $B$ is compatible means that $\psi(p^\perp) = \psi(p)^\perp$ for all $p \in P$. This is equivalent with $B(\psi(p^\perp)) = \psi(p)^\circ$, and since  $\psi(p)^\circ = \psi^*(p^\perp)$, this is the same as $B(\psi(p^\perp)) = \psi^*(p^\perp)$ for all $p \in P$.
\qed \\

The following is Proposition 2.5 (2) in \cite{qu-sc-sc:quadratic}.

\begin{lemma}\label{selfdualcompatibleform}
Let $\psi$ be an indecomposable linear representation in $V$ of a
poset with involution $(P, \perp)$. Then $\psi$ is isomorphic to its dual if and only if there exists a compatible form.  
\end{lemma}
\pf
If there exists a compatible form, then by Lemma \ref{compatible forms}, such a form defines an isomorphism between $\psi$ and $\psi^*$.  

Conversely, suppose that $B: \psi \rightarrow \psi^*$ is an isomorphism.  Then so is $B^*$, and hence the symmetric and antisymmetric parts $(B+B^*)/2$ and $(B-B^*)/2$ of $B$ are also 
morphisms of representations, as are the endomorphisms
$B^{-1}(B+B^*)/2$ and $B^{-1}(B-B^*)/2$, whose sum is the identity morphism.   By Lemma \ref{local}, the ring of endomorphisms of $\psi$ is local, so the two summands cannot both be degenerate.  It follows that either the symmetric or antisymmetric part of $B$ gives a compatible non-degenerate bilinear form.
\qed \\

\subsection{Symplectification}\label{duality}

\begin{dfn}\label{symplectification}
The \textbf{symplectification}\footnote{This construction is sometimes known in a more general setting
 as
{\bf hyperbolization}, see \cite{qu-sc-sc:quadratic}, because the analogue for symmetric bilinear forms leads to isotropic subspaces in spaces with forms of signature zero, sometimes called ``hyperbolic".   We use ``symplectification'' rather than ``symplectization" because the latter term already refers to the construction of symplectic manifolds from contact manifolds by adding one dimension.}
$\psi^-$ of a linear representation $\psi$ is the representation
$$
\psi^- :(P,\perp) \longrightarrow \Sigma (V^* \oplus V, \Omega), \quad \psi^-(x) = \psi^*(x) \oplus \psi(x), 
 $$
 where
$V^* \oplus V$ is endowed with the canonical symplectic structure
$$\Omega((\xi, v),  (\eta, w)) := \xi(w) - \eta(v) \quad \quad \xi, \eta \in V^* \quad v,w \in V.$$ 
\end{dfn}

\begin{prop}
$\psi^-$ is a symplectic representation.
\end{prop}

\pf
For any $p \in P$, we have $\psi^-(p^\perp)= \psi^*(p^\perp) \oplus \psi(p^\perp) = \psi(p)^\circ \oplus \psi(p^\perp).$  This is the symplectic orthogonal of $\psi(p^\perp)^\circ \oplus \psi(p) = \psi^*(p)\oplus \psi(p) = \psi^-(p).$ \qed \\

Here are some fundamental properties of the symplectification operation.

First of all, given a symplectic vector space $V = (V, \omega)$, we denote by $\overline{V}$ the symplectic vector space $(V, \overline{\omega})$, where we define $\overline{\omega} := - \omega$.  Given a symplectic poset representation $\varphi$ on $(V,\omega)$, we define a symplectic poset representation $\overline{\varphi}$ on $\overline{V}$ by setting $\overline{\varphi}(x) = \varphi(x)$ for all $p \in P$, \text{i.e.} as morphisms of posets, $\varphi$ and $\overline{\varphi}$ are the same, only the form on $V$ has been changed. 


\begin{prop}\label{symp hyperbolization}
For any symplectic representation $\varphi$ on $V$, the symplectic representations $\hat \varphi^- = \hat \varphi^* \oplus \hat \varphi$ on $V^* \oplus V$ and $\varphi \oplus \overline{\varphi}$ on $V \oplus \overline{V}$  are isomorphic.
\end{prop}

\pf An isomorphism of symplectic representations is given by 
 $$\tau: \varphi \oplus \overline{\varphi} \rightarrow \hat \varphi^* \oplus \hat \varphi, \ (v, w) \mapsto (\tfrac{1}{2} \tilde{\omega}(v + w), v - w).$$
Indeed, $\tau$ is a morphism of representations, since when $(v,w) \in  \varphi(x) \oplus \overline{\varphi}(x) = \varphi(x) \oplus \varphi(x)$, then $\tau(v,w)) \in \tilde{\omega}(\varphi(x)) \oplus \varphi(x) = \varphi(x^\perp)^\circ \oplus \varphi(x) = \hat \varphi^-(x)$.
And $\tau$ is a symplectic isomorphism, since  
\begin{align*}
\Omega (\tau(v,w),\tau(v',w')) &= \tfrac{1}{2}\tilde{\omega}(v+w)(v' - w') - \tfrac{1}{2}\tilde{\omega}(v'+w')(v-w) \\
		&= \tfrac{1}{2}[\omega(v,v') + \omega(v,-w') + \omega(w,v') + \omega(w,-w') \\
		&\quad \quad - \omega(v',v) - \omega(v',-w) - \omega(w',v) - \omega(w',-w)] \\
		&=  \tfrac{1}{2}[2\omega(v,v') -2 \omega(w,w')] \\
		&= \omega \oplus \overline{\omega} ((v,v'),(w,w')). 
\end{align*} \qed 

\begin{rmk}
The symplectic isomorphism used in the proposition above is the same as  the one behind the Weyl symbol calculus for pseudodifferential operators (c.f. for instance  \cite{dgp:metaplectic-wigner}, formula (7) and Theorem 8) and the definition of ``Poincare's generating function" in hamiltonian mechanics (c.f. \cite{we:poincare}). 
\end{rmk}

\begin{prop}\label{decomposable symplectification}
The symplectification $\psi^-$ of an indecomposable linear representation $\psi$ is symplectically decomposable if and only if $\psi$ admits a compatible symplectic structure.  
\end{prop}
\pf
$\varphi=\psi^-$ is by definition decomposed linearly into the indecomposables $\psi^*$ and $\psi$.  Suppose that it is also symplectically decomposable into two symplectic representations, $\varphi_1$ and $\varphi_2$.  The latter decomposition is also a linear decomposition of $\hat{\varphi}$, and so by Theorem \ref{krs} the linear representations $\varphi_1$ and $\varphi_2$ must be isomorphic  to $\psi^*$ and $\psi$ in some order.  In particular $\psi$ (as does $\psi^*$) admits a compatible symplectic structure.

Conversely, suppose that $\varphi$ admits a compatible symplectic structure.  Then, by Proposition \ref{symp hyperbolization},
$(\hat{\varphi})^+ =(\hat \varphi)^* \oplus \hat \varphi$ and $\varphi \oplus \overline{\varphi}$ are isomorphic symplectic representations. Since the latter is decomposable, so is the former.
\qed \\

\begin{prop}\label{symplectifications iso classes}\label{isomorphic dual}
If $\psi_1$ and  $\psi_2$ are indecomposable linear representations, then $\psi_1$ is isomorphic to $\psi_2$ or to $\psi_2^*$ if and only if the symplectifications $\psi_1^-$ and $\psi_2^-$ are isomorphic as symplectic representations.  

In particular, two symplectifications of indecomposable linear representations are isomorphic as symplectic representations if and only if they are isomorphic as linear representations.
\end{prop}
\pf
If $\psi_1$ is isomorphic to $\psi_2$ or to $\psi_2^*$, then 
$\psi_1^-$ is isomorphic to  $\psi_2^-$ or to $(\psi_2^*)^-.$
If the former holds, we are done.  For the latter, we must show that 
$\psi_2^-$ is isomorphic to $(\psi_2^*)^-.$  They are clearly isomorphic as linear representations, but under the isomorphism which exchanges the summands, the symplectic structures differ by a factor of $-1$. To correct for this factor, we compose with the antisymplectic isomorphism $(\xi,v)\mapsto (-\xi,v)$ from $\psi_2^-$ to itself.

Conversely, if $\psi_1^*(x) \oplus \psi_1(x)$ and $\psi_2^*(x) \oplus \psi_2(x)$ are isomorphic as symplectic representations, then they are in particular also isomorphic as linear representations. This implies that their indecomposable summands are isomorphic in some order, so either $\psi_1\simeq \psi_2$ or $\psi_1\simeq \psi_2^*$ 
\qed \\

\begin{ex}\label{symplectificationexample}
The following are the symplectifications of the indecomposable representations of the poset $\mathbf 2$, \text{i.e.} nested pairs of subspaces.  Each pair is contained in $\fieldk$, so the symplectification is contained in $\fieldk^*\oplus \fieldk$ and is symplectically indecomposable.
\begin{itemize}
\item
The symplectification of $\fieldk\supseteq 0$ is $\fieldk^*\oplus \fieldk \supseteq 0 \oplus 0$.
\item
The symplectification of $0\supseteq 0$ is $\fieldk^*\oplus 0 \supseteq \fieldk^*\oplus 0$.  
\item
The symplectification of $\fieldk\supseteq \fieldk$ is $0\oplus \fieldk \supseteq 0\oplus \fieldk$.
\end{itemize}
The first example is self-dual, while the latter two examples are dual to one another, and their symplectifications are isomorphic by a ``$90$-degree rotation''.  
\end{ex}

Given a pair of linear poset representations, we say that they are \textbf{mutually dual} or, synonymously, a \textbf{dual pair} if each representation is isomorphic to the dual of the other. On the level of isomorphism classes, symplectification builds symplectic representations by taking the direct sum of dual pairs of linear representations. 

\subsection{Relating symplectic and linear indecomposability}
The following lemma, due to Quebbemann et \text{al.} \cite[Thm.3.3]{qu-sc-sc:quadratic} and  Sergeichuk \cite[Lemma 2]{se:classification systems} will be an essential tool in this paper.  It shows that symplectically indecomposable but linearly decomposable representations arise only through symplectification.
The analogous result holds in the symmetric setting.

\begin{lemma}\label{sergeichuk}
Suppose that  $\varphi : (P,\perp) \rightarrow \Sigma (V)$ is an indecomposable symplectic representation such that $\hat \varphi$ is (linearly) decomposable. Then there exists an indecomposable linear representation $\psi$ such that $\varphi \simeq \psi^-$.   
\end{lemma}

\pf 
Because $\hat \varphi$ is linearly decomposable, there exists a non-trivial idempotent $\pi_1 \in \text{End}(\hat \varphi)$.  After two modifications, $\pi_1$ will be conjugated into an idempotent endomorphism $\pi$ satisfying 
the hypotheses $\pi^t \pi = 0 = \pi\pi^t$ and $\pi+\pi^t = 1$ of Lemma \ref{lagrangian direct sum}, giving the required decomposition.

By
Lemma \ref{transpose endomorphism}, the idempotent $\pi_1^t$ is also an endomorphism.
By Lemma \ref{decomps of symp reps}, $\pi_1^t \neq \pi_1$, since otherwise  $\varphi$ would be decomposable as a symplectic representation. Set $\rho_1 =\pi_1 \pi_1^t$. Note that $\rho_1$ is self-adjoint and lies in $\text{End}(\hat \varphi).$ 
By Lemma \ref{symplectic fitting lemma}, $\rho_1$ must be either nilpotent or an isomorphism. But $\rho_1$ cannot be an isomorphism, since $\pi_1$ and $\pi_1^t$ have nontrivial kernels and cokernels.  So $\rho_1$ must be nilpotent. 

Now set $h_1 := s(\rho_1)$, where $s(X)$ is the binomial series for $(1-X)^{1/2}$; $s(\rho_1)$ is well-defined because $\rho_1$ is nilpotent, which implies that the power series is just a polynomial in $\rho_1$. Note that $h_1 \in \text{End}(\hat \varphi),$ and that $h_1$ is also self-adjoint. Furthermore, $h_1$ is invertible, its inverse being defined by substituting $\rho_1$ in the binomial series for $(1-X)^{-1/2}.$  

Define $\pi_2:= h_1 \pi_1 h_1^{-1}$, and note that $\pi_2$ lies in $\text{End}(\hat \varphi)$ and is again a non-trivial idempotent. Furthermore, 
$$ \pi_2^t \pi_2 = h^{-1} \pi_1 h_1^2 \pi_1^t h_1^{-1} = h_1^{-1} \pi_1 (1 - \pi_1 \pi_1^t) \pi_1^t h_1^{-1} = h_1^{-1} ( \pi_1 \pi_1^t - \pi_1 \pi_1^t) h_1^{-1} = 0.$$  

We are half-way there.   Now $\rho_2 := \pi_2 \pi_2^t$ is a nilpotent, self-adjoint element of $\text{End}(\hat \varphi)$, and $h_2 := s(\rho_2)$ is again  an invertible, self-adjoint endomorphism of $\hat \varphi$. Then $\pi := h_2^{-1} \pi_2 h_2 \in \text{End}(\hat \varphi)$ is a non-trivial idempotent such that 
$$ \pi \pi^t = h_2^{-1} \pi_2 h_2^2 \pi_2^t h_2^{-1} = h_2^{-1} \pi_2 (1 - \pi_2 \pi_2^t) \pi_2^t \tilde  h^{-1} = h_2^{-1} ( \pi_2 \pi_2^t - \pi_2 \pi_2^t) h_2^{-1}  = 0 $$
and 
\begin{align*}
\pi^t \pi = & \ h_2 \pi_2^t (h_2^{-2}) \pi_2 h_2 =  h_2 \pi_2^t (1- \pi_2 \pi_2^t)^{-1} \pi_2 h_2 \\
		& =  h_2 \pi_2^t  (1 + \pi_2 \pi_2^t) \pi_2 h_2 =  h_2 (\pi_2^t  \pi_2 +  \pi_2^t  \pi_2 \pi_2^t \pi_2)   h_2 = 0,
\end{align*}
since $\pi_2^t \pi_2 = 0$. 
Furthermore, $\pi + \pi^t \in  \text{End}(\hat \varphi)$ is idempotent: $(\pi + \pi^t)^2 = \pi^2+ \pi^t \pi + \pi \pi^t + (\pi^t)^2 = \pi + \pi^t$. 
But $\pi + \pi^t$ is also self-adjoint, so, by  Lemma \ref{decomps of symp reps}, $\pi + \pi^t$ must be a trivial idempotent. It cannot be that $\pi + \pi^t = 0$, since this would imply $\pi^t = -\pi$, whence $0 = \pi^t \pi = -\pi^2 = \pi$, a contradiction to $\pi \neq 0$. Thus $\pi + \pi^t  =1$. 
\qed \\ 

\begin{cor}\label{sympose}
If two symplectically indecomposable but linearly decomposable representations are isomorphic as linear representations, then they are isomorphic as symplectic representations.
\end{cor}

\pf
By Lemma \ref{sergeichuk}, each of the two symplectic representations is the symplectification of an irreducible linear representation.  Since the two representations are linearly isomorphic, by Lemma \ref{isomorphic dual}, they are symplectically isomorphic.
\qed \\

Lemma \ref{sergeichuk} tells us that every indecomposable symplectic representation is either linearly indecomposable or the symplectification of a linearly indecomposable representation, but not both. In the former case we say that $\varphi$ is of \textbf{non-split} type; in the latter case we say that $\varphi$ is of \textbf{split} type.

\subsection{Uniqueness of compatible forms}\label{uniqueness of compatible forms}
We briefly discuss the question of uniqueness of compatible bilinear forms and formulate 
a lemma which will be the basis of our analysis of uniqueness of compatible forms.  Later we will verify the hypothesis of the lemma for the representations which concern us, and one particular type will require a generalization. The proof uses ideas from that of Lemma 5 in \cite{se:classification systems}; it is interesting that the proof uses a version of the ``square root'' construction of Lemma 2 of that paper (which is our Lemma \ref{sergeichuk}).  In addition, we use an idea (simplified for our context) from Proposition 2.5
in \cite{qu-sc-sc:quadratic} when showing that any two compatible forms must both be symmetric or antisymmetric.

\begin{lemma}\label{uniqueness up to a scalar}
Let $\psi$ be a linearly indecomposable representation in $V$ of an involutive poset for which the endomorphism algebra $E$ (which is local by Corollary \ref{local}) and has the property that $E = \fieldk \id \oplus \text{Rad}E$.  If $\psi$ admits two compatible bilinear forms, then these forms are equal up to a constant scalar multiple and an automorphism of $\psi$.  In particular, the forms must both be symmetric or antisymmetric.
\end{lemma}
\pf 
If $B_1$ and $B_2$  are the isomorphisms from $\psi$ to $\psi^*$ corresponding to two compatible forms, 
then $C= B_1^{-1} B_2$ is an automorphism of $\psi$.

Let $^\dagger$ denote the antiautomorphism of $E$ given by the operation of adjoint with respect to $B_1$, i.e. 
$B_1(A^\dagger x) (y) = B_1(x)(Ay),$ which is equal to $A^*(B_1(x))(y)$, so we have $A^\dagger = B _1^{-1} A^* B_1.$  

Define the signs $\varepsilon_1$ and $\varepsilon_2$ by $B_i^* = \varepsilon_i B_i$, and let $\epsilon =\varepsilon_1 \varepsilon_2.$   Then we have $C^\dagger = \epsilon C$.  In fact, 
$$C^\dagger =(B_1^{-1} B_2)^\dagger = B_1^{-1}B_2^*(B_1^{-1})^*B_1
    = \epsilon B_1^{-1}B_2B_1^{-1}B_1 = \epsilon B_1^{-1}B_2 = \epsilon C.$$
    
By our hypothesis, we may write $C$ as $c {\sf id} - r_1$, where $r_1\in R$ and $c$ is a scalar
and $c\neq 0$ since $C$ is invertible.  By replacing $B_2$ by   $c^{-1} B_2$ and repeating the argument up to this point, we may assume that $c =1,$ so that $C = 1-r$ for an $r \in R$.

Now if $\epsilon$ were equal to $-1$, we would have $$1-r^\dagger = (1-r)^\dagger = -(1-r),$$ which would imply that $1+1 = r^\dagger - r$, which is impossible since $R$ is closed under addition and contains no invertible elements (and, by assumption, $2\neq 0$ in our ground field).  So $\epsilon = 1$, and both forms are either symmetric or antisymmetric.  Furthermore, $r = r^\dagger$.  

By Lemma \ref{local}, we know that $r$ is nilpotent, and so we can use the formal power series for $\sqrt{1-r}$ to construct an automorphism $h$ such that $h^\dagger = h$ and $h^2 = C.$  

Now we have $$h^* B_1 h = B_1 h^\dagger B_1^{-1} 
             B_1 h = B_1 h^\dagger h = B_1 h^2 =B_1 C = B_2, $$
  which shows that $h$ is an isomorphism between the bilinear forms $B_1$ and $B_2$.
\qed \\ 
\begin{rmk}\label{squares}
A form $\omega$ on $V$ and its scalar multiple $a\omega$ are equivalent by a homothety of $V$ (which is automatically an automorphism of any poset representation) if and only if $a$ is a square in the coefficient field.  This means that the set of equivalence classes of compatible forms under homothety is a principal homogeneous space of the {\bf square class group} of $\fieldk$, defined as the quotient $\fieldk^\times/{\fieldk^\times}^2$ of the multiplicative group of nonzero elements of $\fieldk$ by the perfect squares.

Even if $a$ is not a square, $\omega$ and $a\omega$ might, a priori, still be isomorphic by a linear isomorphism which preserves a particular poset representation. This may also be the case if the representation is linearly decomposable, though in the present paper we will not consider the question of compatible forms for decomposable linear representations.  
\end{rmk}

\begin{rmk}\label{what's ahead}
In Section \ref{uniqueness discrete selfdual} we will see that certain self-dual poset representations fulfill the hypotheses of Lemma \ref{uniqueness up to a scalar}. In these cases it turns out that compatible forms which differ by a scalar $c \in \fieldk$ are in fact equivalent if and only if $c$ is a square. 

For certain other self-dual poset representations, however, only a generalization of Lemma \ref{uniqueness up to a scalar} applies and both symmetric and antisymmetric compatible forms exist for a given such representation; see in particular Sections \ref{uniq:gen} and \ref{unihom}. 
\end{rmk}

\section{Interlude: overview, preview, examples}
\subsection{Next steps}
By  Theorem \ref{krs} and Lemma \ref{sympose}  one obtains the following.

\begin{cor}
For a given symplectic representation,
consider orthogonal decompositions into  symplectically indecomposable summands.
Then any summand is either split or non-split and the
 following are uniquely determined: 
\begin{enumerate}
\item The isomorphism types and multiplicities of split summands.
\item The linear isomorphism types and multiplicities of
non-split summands.
\end{enumerate}
\end{cor}
At this point, we have finished our discussion of the general theory of symplectic poset representations and are ready to move on to specific cases related to isotropic triples.  We will do the following things. 

\begin{itemize}
\item
In section \ref{quivers} we review some essentials of the theory of quiver representations  -- as background for the classification results to be used. We explain how the classification, obtained in \cite{dl-ri-indecomposable,do-fr:representations}, of indecomposable representations of a quiver related to the Dynkin diagram $\tilde{E}_6$ gives the classification of indecomposable linear representations of the poset  ${\bf 2}+{\bf 2}+{\bf 2}$. 

\item As a warmup, we derive, in sections \ref{subspaces} and \ref{pairs}, the classification of isotropic 
$k$-tuples, for $k=0,1,2$, from
the classification of representations of quivers
associated with the Dynkin-diagrams $A_1, A_3,$ and $A_5$. 
Then in section \ref{2+2+2} we give an overview of the quiver representation classification results for the Dynkin diagram $\tilde{E}_6$, which we will use for classifying isotropic triples. 
In section \ref{dim 2}, again as a warmup, we discuss isotropic triples in ambient dimension $2$, and give, in section \ref{higher dim}, a preview of the situation in higher dimensions. 
\item
Using the results of \cite{dl-ri-indecomposable,do-fr:representations}, we identify which indecomposable linear representations of $P= {\bf 2}+{\bf 2}+{\bf 2}$ are dual to one another when we equip ${\bf 2}+{\bf 2}+{\bf 2}$ with the involution $\perp$ which exchanges the respective elements of the three nested pairs in ${\bf 2}+{\bf 2}+{\bf 2}$. From now on we will use the term \textbf{sextuples} to refer to linear representations of this poset with involution.

\item
From the general theory we know that self-dual sextuples admit compatible symmetric or symplectic forms (or both). We determine which self-dual sextuples
admit compatible symplectic forms (thus giving non-split isotropic triples), and we give explicit
constructions of such forms.

\item For the self-dual sextuples, we reduce the classification of compatible forms to a field-theoretic description. When $\fieldk = \mathbb{R}$ or when $\fieldk$ is algebraically closed, compatible forms are parametrized by the square class group $\fieldk^\times / (\fieldk^\times)^2$. For general perfect fields, a similar, though slightly more complicated, description is obtained; see Theorem \ref{uniqueness general fields}. 

\end{itemize}

\subsection{Quiver and poset representations}\label{quivers}

We recall here some basic definitions and results in quiver representation theory, referring to the literature\footnote{For example $\{$\cite{barot}, \cite{benson-book}, \cite{de-we:intro-quivers}, \cite{etingof-et-al},  \cite{kac}, \cite{schiffler}$\}$ is a small sample subset of the available references, to give the reader a starting point.}  for more details.
A {\bf quiver} $\calq$ is simply a directed graph, i.e.~a set $\mathcal V$ of vertices and a set $\mathcal{A}$ of arrows, with source and target maps $s$ and $t$ from $\mathcal{A}$ to $\mathcal{V}$.  We allow multiple arrows with a given source and target, but assume the sets $\mathcal{A}$ and $\mathcal{V}$ to be finite.   With a chosen ground field $\fieldk$,  a {\bf representation} $\rho$ of $\calq$ is simply an assignment to each vertex $v$ a (finite dimensional, for our purposes) $\fieldk$-vector space $\rho(v)$ and to each arrow $a$ a linear map $\rho(a):\rho(s(a))\to \rho (t(a))$.  A morphism $\mu$ from $\rho_1$ to $\rho_2$ is a family of linear maps $\mu_v:\rho_1(v) \to \rho_2(v)$ making the obvious diagrams commute. When the family of linear maps consists of isomorphisms, then $\mu$ is called an isomorphism. The collection of representations of a fixed quiver with their morphisms form a category.

A fundamental problem in the theory of quiver representations (as is the case for representations of just about anything) is to describe the structure of the set of isomorphism classes of representations, and, among these, the indecomposable representations, which are those not decomposable into nontrivial direct sums.

The first basic result is the Krull-Schmidt theorem, which states that each representation of a quiver is isomorphic to a direct sum of indecomposables, and that the summands in this decomposition, with their multiplicities, are unique up to isomorphism and reordering. 
This reduces the classification of representations to the enumeration of those which are indecomposable.  When the set of isomorphism classes of indecomposables is finite, the quiver is of {\bf finite type}. 

As mentioned above, we will be studying isotropic triples by considering them as linear representations of the poset $\textbf{2}+\textbf{2}+\textbf{2}$.  
These poset representations can be identified with particular representations of 
the following quiver; it is obtained by choosing the following orientation on the extended Dynkin diagram $\tilde{E}_6$:
\begin{equation}\label{Q}
\begin{tikzpicture}[>=Stealth, xscale=1, yscale=1]
\tikzmath{
coordinate \c; 
\c0 = (0,0.5);
\c1 = (1.5,.9);
\c2 = (1.5,0.5);
\c3 = (1.5,0.1);
\c4 = (3,0.9);
\c5 = (3,0.5);
\c6 = (3,0.1);
}

%
%

Vertex labels
\node at (\c0) {\small{$v$}}; 
\node at (\c1) {\small{$c_1$}}; 
\node at (\c2) {\small{$c_2$}}; 
\node at (\c3) {\small{$c_3$}}; 
\node at (\c4) {\small{$i_1$}}; 
\node at (\c5) {\small{$i_2$}}; 
\node at (\c6) {\small{$i_3$}};

\tikzmath{
coordinate \tip, \tail;
\tip1 = (\c0) + (0.18,0.08);
\tail1 = (\c1) + (-0.22,-0.04);
\tip2 = (\c0) + (0.18,0.0);
\tail2 = (\c2) + (-0.22,0);
\tip3 = (\c0) + (0.18,-0.08);
\tail3 = (\c3) + (-0.22,0.03);
\tip4 = (\c1) + (0.15,0.0);
\tail4 = (\c4) + (-0.22,0);
\tip5 = (\c2) + (0.15,0.0);
\tail5 = (\c5) + (-0.22,0);
\tip6 = (\c3) + (0.15,0.0);
\tail6 = (\c6) + (-0.22,0);
}

\draw [<-] (\tip1) -- (\tail1);
\draw [<-] (\tip2) -- (\tail2);
\draw [<-] (\tip3) -- (\tail3);

\draw [<-] (\tip4) -- (\tail4);
\draw [<-] (\tip5) -- (\tail5);
\draw [<-] (\tip6) -- (\tail6);

\end{tikzpicture}
\end{equation}
We'll refer to this quiver also as $\tilde{E}_6$. The labels on the vertices are in principle arbitrary (and sometimes unnecessary); we fix this choice since it is suggestive for our application to (co)isotropic triples. We identify the vertices other than $v$ with the elements of the poset $\textbf{2}+\textbf{2}+\textbf{2}$, with the partial order indicated by the arrows, \text{i.e.} $i_1 \leq c_1$, $i_2 \leq c_2$ and $i_3 \leq c_3$. 
%

When considering representations of this quiver, we will denote the space associated to the vertex $v$ by $V$, and the spaces associated to the vertices $c_1,c_2, c_3$ and $i_1, i_2, i_3$ will be denoted by $C_1, C_2, C_3$ and $I_1, I_2, I_3$, respectively. When it is clear what the maps are, a representation of (\ref{Q}) will be denoted by the 7-tuple of spaces $(V; C_1, I_1; C_2, I_2; C_3, I_3)$ or just $(V;C_i,I_i)$. We will also call such representations $\textbf{sextuples}$, just as we do poset representations of $\textbf{2}$+$\textbf{2}$+$\textbf{2}$. Abusing notation slightly, dimension vectors of sextuples will be denoted $(v, c_1, c_2, c_3, i_1, i_2, i_3)$, where the entries denote the dimensions of the (sub)spaces of a representation $(V; C_1, I_1; C_2, I_2; C_3, I_3)$. 

It is straightforward to see that, on the level of isomorphism classes, linear poset representations of $\textbf{2}+\textbf{2}+\textbf{2}$ are in one-to-one correspondence with quiver representations of the quiver $\tilde E_6$ where all arrows are represented by injective linear maps (we will call these injective representations). This correspondence is compatible with the notions of direct sum for poset and quiver representations, respectively. Furthermore, one can prove that an indecomposable quiver representation of $\tilde E_6$ is an injective representation if and only if the space assigned to the central vertex is non-zero (see \cite{raf2}, Proposition A.7.1). 

Using this, one can read off the indecomposable representations of
$\textbf{2}$+$\textbf{2}$+$\textbf{2}$ from those of $\tilde{E}_6$
as given in \cite{dl-ri-indecomposable,do-fr:representations}.  We will mainly use the explicit normal forms given in \cite{dl-ri-indecomposable}. Since that reference (in contrast to \cite{do-fr:representations}) only treats the case when the ground field is algebraically closed, we will use a straightforward generalization of their normal forms for more general fields. By inspection of the proofs \cite{dl-ri-indecomposable}, it is only the normal forms for continuous-type representations that must be generalized. These are discussed in Section \ref{continuous}.

%
%
%

An essential tool for the study of quiver representations is the {\bf Tits form} of a quiver $\calq$.  This is the quadratic form $q$ on the 
$\mathbb{Z}$-module generated by the vertices defined by $q(\mathbf n)=
 \sum _\calv n_v^2 - \sum_{\calv \times \calv} a_{v,w}n_v n_w$, where $n_v$ is the coefficient of $v$ in $\mathbf{n}$, and $a_{v,w}$ is the number of arrows from $v$ to $w$.  Note that this form does not depend on the direction of the arrows.
 
The idea behind the Tits form is that, if the coefficients of $\mathbf n$ are the dimensions of vector spaces assigned to the vertices, then the second term is (the negative of) the dimension of the linear space of all representations of $\calq$ in this family of vector spaces, while the first term is the dimension of the group, acting on the space of representations, whose orbits are the isomorphism classes. In fact, the scalar multiples of the identity act trivially, so we may say that the ``virtual dimension'' of the moduli space of isomorphism classes of representations with dimension vector $\mathbf n$ is $1-q(\mathbf n)$.  The actual dimension is at least this large (and larger if more of the group acts trivially), so if our ground field is, e.g., the real or complex numbers, the only way in which there can be finitely many isomorphism classes with dimension vector $\mathbf n$ is if $q(\mathbf n)$ is at least $1$.  
 
This suggests (but does not prove, for various reasons), that a quiver is of finite type if and only if its Tits form is positive definite.  In fact, this is true (for any ground field!) and is part of what is known as {\bf Gabriel's Theorem}.   The other part of the theorem states that the connected quivers of finite type are exactly those for which the associated undirected graph is a Dynkin diagram of type $A$, $D$, or $E$ \cite{gab}.
 For these quivers, it turns out that the nonnegative solutions of $q(\mathbf{n}) = 1$, known as the {\bf positive roots}, are precisely the dimension vectors of indecomposable representations, and there is exactly one isomorphism class corresponding to each positive root.  
 
The quiver (\ref{Q}) relevant to the classification of isotropic triples is not of finite type.   It does belong, however, to the ``next best'' class, that of the so-called {\bf tame} quivers.  For these, the Tits form is positive {\em semi}definite, with one-dimensional null space which we will denote by $N$.  $N$ has a smallest positive element, which is the dimension vector of a family of representations whose isomorphism classes are parametrized, in the case of $\fieldk$ algebraically closed, by a 1-dimensional variety.  The positive roots thus fall into lines parallel to $N$.  They still correspond to indecomposable representations, which now belong to families indexed by the positive integers.

An extension of Gabriel's theorem (\text{c.f.} \cite{dl-ri-indecomposable}, \cite{do-fr:representations}) tells us that a quiver is of (infinite) tame type if and only if the  corresponding undirected graph is an {\bf extended Dynkin diagram}; these are obtained from certain Dynkin diagrams by the addition of an edge attached to an extremal vertex.  Among these, for instance, is $\widetilde{D_4}$, consisting of four edges attached to a central vertex.   If the edges are all oriented to point toward the vertex, the representations of the quiver are closely related to those of the partially ordered set $\textbf{1}$+$\textbf{1}$+$\textbf{1}$+$\textbf{1}$ consisting of four incomparable elements.
Representations of this poset are just quadruples of subspaces in vector spaces.  These were classified by Gel'fand and Ponomarev \cite{ge-po}, who showed that many classification problems in linear algebra reduce to the classification of certain subspace quadruples.
For instance, endomorphisms of a vector space $V$ in dimension $n$ correspond to certain kinds of quadruples of subspaces of dimension $n$ in $V\oplus V$ (the ``axes,'' the diagonal, and the graph of the endomorphism).  Indecomposable representations of this kind correspond to indecomposable endomorphisms which, in the case of an algebraically closed field, are those given by a single Jordan block.   Since the diagonal element and the size of such a block is arbitrary, one sees immediately the presence of 1-parameter families with arbitrarily large dimension vectors.  This example also shows the possible complications for fields which are not algebraically closed, where indecomposable endomorphisms are parametrized by irreducible polynomials which are no longer necessarily linear.
We will see later that there are close connections between the representations of $\widetilde{D_4}$ and those of $\widetilde{E_6}$, the latter also being connected to the classification of endomorphisms.

\medskip

In the next two sections, we reinterpret, via the theory of quiver and poset representations, the easy classifications of symplectic vector spaces and isotropic subspaces, followed by the classification of isotropic pairs, c.f. \cite{lo-we:coisotropic}.  The relevant Dynkin diagrams are $A_1$, $A_3$, and $A_5$.


\subsection{Symplectic spaces, isotropic subspaces}\label{subspaces}

We start with the classification of symplectic vector spaces with {\em no} distinguished isotropic subspaces.   We may think of these as the symplectic representations of the empty poset, or of the Dynkin diagram $A_1$, whose quiver consists of a single vertex with no arrows, and whose Tits form in terms of the one-entry dimension vector $(v)$ is the positive definite form $v^2$.  It is well known that any symplectic vector space admits a basis (in fact many such bases) of the form $(e_1,f_1,\ldots e_n,f_n)$, with the symplectic form determined by the conditions that  
$\omega(e_j,f_j)= 1$ for all $j$ and $\omega(a,b)=0$ for all other pairs $(a,b)$ of basis elements.   We call this a {\bf symplectic basis}.  As a consequence,  any symplectic vector space can be decomposed as a direct sum of copies of the space $\fieldk^2$ with symplectic basis $(e,f)$.  The only invariant of a symplectic space is its dimension, which must be an even nonnegative integer.  This is consistent with the fact that the Tits form has a single positive root, which is $(1)$.  

Another viewpoint here is that there is one indecomposable representation of the empty poset, the 1-dimensional space $\fieldk$.
The symplectification of this representation is $\fieldk^*\oplus \fieldk$.   In fact, we will use this description of the 2-dimensional symplectic space, rather than $\fieldk^2$.

We now move on to the example of individual isotropic subspaces, which correspond to certain symplectic representations of the poset $\bf 2$, i.e. nested pairs of subspaces, as noted above, where isotropic $I$ in $V$ corresponds to the representation $V\supseteq I^\perp\supseteq I$.

The quiver associated to this poset is 

$$\begin{tikzpicture}[>=Stealth, xscale=1, yscale=1]
\tikzmath{
coordinate \c; 
\c0 = (0,0);
\c1 = (1.5,0);
\c2 = (3,0);
}

\draw [fill] (\c0) circle [radius=0.05];
\draw [fill] (\c1) circle [radius=0.05];
\draw [fill] (\c2) circle [radius=0.05];

\tikzmath{
coordinate \tip, \tail;
\tip1 = (\c0) + (0.15,0);
\tail1 = (\c1) + (-0.2,0);
\tip2 = (\c1) + (0.15,0);
\tail2 = (\c2) + (-0.2,0);
}

\draw [<-] (\tip1) -- (\tail1);
\draw [<-] (\tip2) -- (\tail2);

\end{tikzpicture}
$$
with underlying Dynkin diagram $A_3$.
For a dimension vector of the form $(v,c,i)$, corresponding to a representation $V\from C\from I$, the Tits form is $q(v,c,i) = v^2+c^2+i^2 - vc - ci$, which can be rewritten as $$\frac{1}{2}(v^2 +(v-c)^2+(c-i)^2+i^2),$$ which is clearly positive definite.  The dimension vectors of indecomposable representations of the quiver are the positive roots, i.e. nonnegative integer solutions of $q(v,c,i)=1$.  These must be vectors such that exactly two of the four squared summands are equal to 1.
Of the six such solutions, those for which the arrows are represented by injective maps are $(1,0,0)$, $(1,1,0)$, and $(1,1,1).$   
Over a ground field $\fieldk$, the corresponding representations are  $(1;0,0):= \fieldk\supseteq 0\supseteq 0$, $(1;1,0):= \fieldk\supseteq \fieldk \supseteq 0$, and $(1;1,1):= \fieldk\supseteq \fieldk \supseteq \fieldk$.  $(1;1,0)$ is self-dual, while $(1;0,0)$ and $(1;1,1)$ are (isomorphic to) the duals of one another. 

On $(1;1,0)$, there are no symplectic forms.\footnote{The compatible nondegenerate bilinear forms on $(1;1,0)$ are symmetric; the group $\fieldk^*/(\fieldk^*)^2$, where $\fieldk^*$ is the multiplicative group of $\fieldk$, acts simply and transitively on the isomorphism classes of such forms.  This quotient group is $\mathbb Z_2$ in the case of a finite field or the real numbers, and the trivial group for an algebraically closed field.  It can be much larger, for instance in the case of the rational numbers. In any case, the isotropic subspace is the zero subspace.}

As we have seen in Example \ref{symplectificationexample}, the symplectifications of these representations, contained in $\fieldk^*\oplus\fieldk$, are $(1;1,0)^-=\fieldk^*\oplus \fieldk \supset 0 \oplus 0$,  $(1;0,0)^- =\fieldk^*\oplus 0 \supseteq \fieldk^*\oplus 0$,  
and $(1;1,1)^-= 0\oplus \fieldk \supseteq 0\oplus \fieldk$.
 
The isotropic subspaces in $\fieldk^*\oplus\fieldk$ are the zero subspace in the first case, and (lagrangian) lines in the latter two cases. The latter two symplectic representations are isomorphic, but we will use both of them in the classification of multiple isotropic subspaces.

\subsection{Isotropic pairs}\label{pairs}
The relevant poset here is $\bf 2 + \bf 2$, to which we associate the quiver 
$$
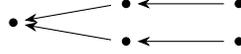
\begin{tikzpicture}[>=Stealth, xscale=1, yscale=1]
\tikzmath{
coordinate \c; 
\c0 = (0,0.25);
\c1 = (1.5,0.5);
\c2 = (1.5,0);
\c3 = (3,0.5);
\c4 = (3,0);
}

\draw [fill] (\c0) circle [radius=0.05];
\draw [fill] (\c1) circle [radius=0.05];
\draw [fill] (\c2) circle [radius=0.05];
\draw [fill] (\c3) circle [radius=0.05];
\draw [fill] (\c4) circle [radius=0.05];

\tikzmath{
coordinate \tip, \tail;
\tip1 = (\c0) + (0.15,0.03);
\tail1 = (\c1) + (-0.2,-0.02);
\tip2 = (\c0) + (0.15,-0.03);
\tail2 = (\c2) + (-0.2,0.02);
\tip3 = (\c1) + (0.15,0);
\tail3 = (\c3) + (-0.2,0);
\tip4 = (\c2) + (0.15,0.0);
\tail4 = (\c4) + (-0.2,0);
}

\draw [<-] (\tip1) -- (\tail1);
\draw [<-] (\tip2) -- (\tail2);
\draw [<-] (\tip3) -- (\tail3);
\draw [<-] (\tip4) -- (\tail4);

\end{tikzpicture}
$$
with underlying Dynkin diagram $A_5$.

For a representation of this quiver given by maps 
$ I_1\to C_1\to V\from C_2\from I_2$, we will write the dimension vector in the form $(v;c_1, i_1  ;c_2 ,i_2)$.

The Tits form may be written as 
 $$q(v;c_1, i_1  ;c_2 ,i_2)=\frac{1}{2}[(v-c_1)^2 + (c_1-i_1)^2 +i_1^2 +(v-c_2)^2+(c_2-i_2)^2+i_2^2]$$
 Again, this is positive definite, and the positive roots are those vectors making exactly two of the squared terms equal to 1.   Those giving representations by injective maps (which can be found by consulting a table of the positive roots of $A_5$) are as follows.  The only self-dual root is  $(1;1,0;1,0)$.  The others, arranged together with their duals, are
 $(1;0,0;0,0)$ and $(1;1,1;1,1)$, $(1;1,0;0,0)$ and $(1;1,0;1,1)$, $(1;0,0;1,0)$ and $(1;1,1; 1,0)$, and $(1;1,1;0,0)$ and $(1;0,0;1,1)$.  Since the self-dual root has an odd-dimensional ambient space, the corresponding representation does not admit a symplectic structure.   We then have five indecomposable isotropic pairs by symplectification, all lying in $\fieldk^*\oplus \fieldk$; they are
$(0\oplus \fieldk,0\oplus \fieldk)$ (two equal lines), $(0\oplus 0,0\oplus \fieldk)$ (zero and a line), $(0\oplus \fieldk,0\oplus 0)$ (a line and zero), $(\fieldk^*\oplus 0,0\oplus \fieldk)$ (two distinct lines), and $(0\oplus 0,0\oplus 0)$ (two zero subspaces).  They correspond exactly (in a different order) to the five symplectic indecomposables numbered 6 through 10 in Theorem 2 of \cite{lo-we:coisotropic}.

\subsection{Overview of indecomposable  representations of the poset $\bf 2+\bf 2+\bf 2$}\label{2+2+2}
We come now to the central object of this paper. As noted earlier, the quiver which we associate with the poset $\bf 2+\bf 2+\bf 2$ governing isotropic triples is 

$$
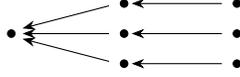
\begin{tikzpicture}[>=Stealth, xscale=1, yscale=1]
\tikzmath{
coordinate \c; 
\c0 = (0,0.5);
\c1 = (1.5,.9);
\c2 = (1.5,0.5);
\c3 = (1.5,0.1);
\c4 = (3,0.9);
\c5 = (3,0.5);
\c6 = (3,0.1);
}

\draw [fill] (\c0) circle [radius=0.05];

\draw [fill] (\c1) circle [radius=0.05];
\draw [fill] (\c2) circle [radius=0.05];
\draw [fill] (\c3) circle [radius=0.05];

\draw [fill] (\c4) circle [radius=0.05];
\draw [fill] (\c5) circle [radius=0.05];
\draw [fill] (\c6) circle [radius=0.05];

\tikzmath{
coordinate \tip, \tail;
\tip1 = (\c0) + (0.15,0.07);
\tail1 = (\c1) + (-0.2,-0.03);
\tip2 = (\c0) + (0.15,0.0);
\tail2 = (\c2) + (-0.2,0);
\tip3 = (\c0) + (0.15,-0.07);
\tail3 = (\c3) + (-0.2,0.03);
\tip4 = (\c1) + (0.1,0.0);
\tail4 = (\c4) + (-0.2,0);
\tip5 = (\c2) + (0.1,0.0);
\tail5 = (\c5) + (-0.2,0);
\tip6 = (\c3) + (0.1,0.0);
\tail6 = (\c6) + (-0.2,0);
}

\draw [<-] (\tip1) -- (\tail1);
\draw [<-] (\tip2) -- (\tail2);
\draw [<-] (\tip3) -- (\tail3);

\draw [<-] (\tip4) -- (\tail4);
\draw [<-] (\tip5) -- (\tail5);
\draw [<-] (\tip6) -- (\tail6);

\end{tikzpicture}
$$
and the  corresponding extended
 Dynkin diagram is $\widetilde{E_6}$.
An explicit description of the indecomposable $\widetilde{E_6}$ representations has been given by Donovan and Freislich in 
\cite{do-fr:representations}, organized into families described
as follows: The dimension vectors of indecomposables are arranged in lines parallel to $N=\mathbb{N}(3;2,1;2,1;2,1)$.  Each of these lines contains a least positive element, followed
 by vectors obtained by adding the elements of $N$.   These give sequences in increasing dimensions. The dimension vectors in $N$ are referred to
as {\bf continuous}, the others as {\bf discrete}; we also use this terminology for the associated indecomposable representations. 

\begin{rmk}\label{ground field}
The classification of discrete-type indecomposables is in fact independent of the ground field $\fieldk$, while the classification of continuous-type indecomposables does (partially) depend on $\fieldk$. This follows from DR \cite{dl-ri-indecomposable}, or by inspection of the proofs in DF \cite{do-fr:representations}. 
\end{rmk}

\begin{rmk}\label{discrete}
In view of Corollary \ref{local} and the fact that indecomposable discrete-type sextuples are uniquely determined up to isomorphism by their dimension vectors, to show that a sextuple is of a given isomorphism type of discrete type it is sufficient to show that it has the required
dimension vector and that its endomorphism ring is local.
\end{rmk}

   

An important step in classifying symplectic representations of the poset $P = \bf 2+\bf 2+\bf 2$ with our chosen involution is to identify which linear representations are self-dual, since these may admit compatible symplectic forms. For the discrete-type dimension vectors, self-duality of the dimension-vector implies self-duality of the (uniquely) corresponding indecomposable linear representation, \text{c.f.} Lemma \ref{dualdimension}. In view of Remark \ref{ground field}, the self-dual discrete-type sextuples may be read off from the classification in \cite{do-fr:representations}:

%
%
%
%

\begin{prop}\label{DFdiscrete}
The self-dual discrete
dimension vectors of indecomposable representations 
of ${\bf 2}+{\bf 2}+{\bf 2}$ are of the form
$(3k+1;2k+1,k,2k+1,k,2k+1,k)$ and $(3k+2;2k+1,k+1,2k+1,k+1,2k+1,k+1)$.
For each of these there is, up to isomorphism, a unique indecomposable, named $A(3k+1,0)$ resp. $A(3k+2,0)$.
In particular, these are the only self-dual discrete sextuples.
\end{prop}


Explicit descriptions of the isotropic triples associated to the self-dual sextuples $A(3k+1,0)$ and $A(3k+2,0)$ are given in Sections \ref{impl A(3k+1,0)} and \ref{impl A(3k+2,0)}. In Section \ref{existence adm forms} compatible symplectic forms are constructed explicitly, and in Section \ref{unidis} we discuss their uniqueness.
The (more difficult) question of duality for continuous-type indecomposable sextuples is discussed in Sections \ref{continuous}, \ref{continuous real and complex}, and \ref{continuous perfect}.


In the present section we give some further explanations of the DF-classification \cite{do-fr:representations}.
 The discrete indecomposable sextuples are labeled in the form $L_s(3k +i, d)$, where 
$k$ can be any non-negative integer and $i \in \{0,1,2,3\}$ (but does not always run over that whole set). $3k+i$ gives the dimension of the ambient space $V$, and $d$ is the defect $\sum c_j + \sum i_j - 3 v$.   
$L$ is a letter in $\{A,B,C,D\}$ which encodes the degree of symmetry of the dimension vector, with $A =$ fully symmetric, i.e. all ``arms" equal, $B$ or $C =$ exactly two arms equal, $D =$ no arms equal. 
The subscript $s$ is either empty (when $L =A$), an integer in $\{1,2,3\}$ when $L = B$ or $C$, or a pair of unequal integers in $\{1,2,3\}$ (when $L = D$). This subscript encodes ``where the asymmetry is".  In cases $B$ or $C$, it tells which arm has a different dimension vector.  In case $D$, it tells how two asymmetric dimension vectors can be related via permutation of the arms. For example, $D_{12}(3k+1,0)$ is related to $D_{31}(3k+1,0)$ via the permutation $1 \rightarrow 3$, $2 \rightarrow 1$, $3 \rightarrow 2$ (but only two indices are needed to specify a permutation of 3 elements). 

There is no essential difference between the cases  $L=B$ and $L=C$.  The reason for using two different letters seems to simply be that, in the defect $-1$ and $+1$ cases, if the parameters  $k$, $i$, $d$, and $s$ are fixed, then there are still two distinct dimension vectors of indecomposables.

The lowest dimensional members of each family of discrete-type indecomposable sextuples are listed in the following table (with $i=1,2,3$):

\begin{equation*}
\begin{array}{c||l|l|l}
\text{defect} & \dim 1 & \dim 2 & \dim 3  \\ 
\hline
-3 & A(1,-3) & A(2,-3)&  \\
-2 & C_i(1, -2) & C_i(2, -2)& C_i(3, -2)   \\ 
-1 & B_i(1, -1) & B_i(2, -1) & B_i(3, -1)  \\ 
- 1& C_i(1, -1) & C_i(2, -1) & C_i(3, -1)   \\ 
0 & D_{i+2,i+1}(1,0) & D_{i+2,i+1}(2,0) &   \\ 
\hline
0 & A(1,0) & A(2,0) &   \\ 
\hline
0 & D_{i,i+1}(1,0) & D_{i,i+1}(2,0) &   \\ 
1 & C_i(1, 1) & C_i(2, 1) & C_i(3, 1) \\ 
1 & B_i(1, 1)  & B_i(2, 1) & B_i(3, 1)    \\ 
2 & C_i(1, 2) & C_i(2, 2) & C_i(3, 2)  \\ 
3 & A(1,3) & A(2,3) &  
\end{array}
\end{equation*}

The row in the middle contains the lowest dimensional members of the two families of discrete indecomposables which are self-dual, namely the families $A(3k+1,0)$ and $A(3k+2,0)$. For any other entry in the table, its dual indecomposable is found by reflecting along the horizontal middle axis, \text{e.g.} $A(1,-3)$ and $A(1,3)$ are mutually dual, $C_i(1, -2)$ and $C_i(1, 2)$ are mutually dual, etc..

In contrast to the discrete-type indecomposable sextuples, the classification of continuous-type indecomposable sextuples is dependent on the ground field (again, this follows from the classification in DF \cite{do-fr:representations} and DR \cite{dl-ri-indecomposable}). Although we will ultimately work in the setting where the ground field is only assumed to be perfect, for illustrative purposes, we assume for the moment that the ground field is the complex numbers. 

The indecomposable continuous-type sextuples can be arranged into families whose lowest dimensional members are listed in the following table. Here, $i = 1,2,3,$ and the parameter $\lambda$ is understood as ranging in the disjoint union of the sets $\mathbb{C} \backslash \{0,1\}$ and  $\{0_i, 1_1, 1_2, \infty_i \}$, where the latter $8$ (formal) elements are labels for certain ``exceptional" indecomposables.

\begin{equation*}
\begin{array}{c||l}
\text{defect} & \dim 3  \\ 
\hline
0 & \Delta(1,0_i)  \\ 
0 & \Delta(1,\lambda), \ 0 < \vert \lambda \vert < 1, \text{ or } \vert \lambda \vert = 1 \text{ with } \text{Im} \lambda > 0 \\ 
0 & \Delta(1,1_1)  \\ 
\hline
0 & \Delta(1,-1)   \\ 
\hline
0 & \Delta(1,1_2)  \\ 
0 & \Delta(1,\lambda^{-1}),  \ 0 < \vert \lambda \vert < 1, \text{ or } \vert \lambda \vert = 1 \text{ with } \text{Im} \lambda > 0 \\ 
0 & \Delta(1,\infty_{i+1})  
\end{array}
\end{equation*}

As with the previous table, this one is also arranged so that dual sextuples are placed symmetrically to each other with respect to reflection around the central horizontal row (and leaving this row fixed); in particular, the only self-dual element of the table is $\Delta(1,-1)$\footnote{We note that, for indecomposable sextuples of the type $\Delta(1,\lambda)$ with $\lambda \in \mathbb{C} \backslash \{0,1\}$, the separation into two families (according to the absolute value of $\lambda$ and the sign of its imaginary part) is something we have introduced ``artificially" for this table in order to emphasize a separation into dual pairs. In fact, the values of $\lambda$ are not intrinsic; the values used here are simply one of many possible ways to parametrize the ``moduli space'' of continuous-type indecomposable sextuples.}.

Finally, we review  these results in the context of Section \ref{quivers}.
The quiver we work with, $\tilde{E}_6$, is tame but not of finite type, and its indecomposable representations, and hence those of the underlying poset, are  infinite in number and of arbitrarily high
dimension, with some of them of discrete type and others of continuous type. 
Following the pattern in the preceding subsections, we will denote a dimension vector by $(v;c_1,i_1;c_2,i_2;c_3,i_3)$.
The Tits form 
\begin{align*}
q(v;c_1, i_1  ;c_2 ,i_2;c_3,i_3)=&\frac{1}{2}[-v^2+(v-c_1)^2 + (c_1-i_1)^2 +i_1^2 \\
	& +(v-c_2)^2+(c_2-i_2)^2+i_2^2+(v-c_3)^2 + (c_3-i_3)^2 +i_3^2]
\end{align*}
is now positive {\em semi}definite; its null space $N$ is 1-dimensional, generated over the integers by
$\nu := (3;2,1;2,1;2,1)$. Dlab-Ringel \cite{dl-ri-indecomposable} give a nice presentation of (a constant times) this form as a sum of six squares.  Since there are 7 variables, this shows positive semidefiniteness.

The discrete indecomposable sextuples are those with dimension vectors for which the Tits form takes the value $1$.   


\subsection{Triples in dimension 2}\label{dim 2}
We will enumerate here the isotropic triples in dimension 2, which are all symplectically indecomposable.  Since many of them are linearly decomposable, we will be using the description in Example \ref{symplectificationexample} of the symplectifications of the three nested pairs in $\fieldk$.

The following are the possibilities for a triple, in the symplectic plane $\fieldk^*\oplus \fieldk$, of pairs each consisting of a coisotropic subspace containing its isotropic orthogonal.  These are necessarily symplectically indecomposable, but only the final example is linearly indecomposable.  The rest are  symplectifications of 1-dimensional linear representations.
\begin{enumerate}
\item 
Three copies of $\fieldk^*\oplus\fieldk\supseteq 0\oplus 0$.  (The isotropics are all zero.)  This is the symplectification of the sextuple in $\fieldk$ consisting of three copies of $\fieldk\supseteq 0$, whose dimension vector is  $(1;1,0;1,0;1,0)$.  In the DF classification, this is $A(1,0)$.  It is self-dual, but all compatible bilinear forms are symmetric rather than symplectic.  
\item
Two copies of  $\fieldk^*\oplus\fieldk\supseteq 0\oplus 0$ and one copy of $0\oplus \fieldk \supseteq 0\oplus \fieldk$.  (Two isotropics are zero, and one is a line.)  This is the symplectification of two copies of $\fieldk\supseteq 0$ and one copy of $\fieldk\supseteq\fieldk$ (or its dual $0\supseteq 0$). In terms of the DF classification, the linear representations being symplectified are of the type $B_r(1,1)$ (or its dual $B_r(1,-1)$), for $r=1$, $2$, or $3$.  Thus there are three possibilities here, depending upon the value of $r$ (i.e., on which of the three isotropics is a line).  The $B_r(1,\pm 1)$ are the first members of the families $B_r(3k + 1,\pm 1)$.  $B_3(1, 1)$, for example, has the dimension vector $(1;1,0;1,0;1,1)$, and its dual $B_3(3k + 1, -1)$ has dimension vector $(1;1,0;1,0;0,0)$.  

\item
One copy each of $\fieldk^*\oplus\fieldk\supseteq 0\oplus 0$, $\fieldk^*\oplus 0 \supseteq \fieldk^*\oplus 0$, and $0\oplus\fieldk\supseteq 0\oplus\fieldk$.
(One isotropic is zero, and the other two are different lines.)
This is the symplectification of one copy each of $\fieldk\supseteq 0$, $0\supseteq 0$, and $\fieldk\supseteq\fieldk$. 
In terms of the DF classification, the linear representation being symplectified is $D_{12}(1,0)$, $D_{31}(1,0)$, or $D_{23}(1,0)$ (or the dual 
$D_{32}(1,0)$, $D_{21}(1,0)$, or $D_{13}(1,0)$ respectively).  Again, there are three possibilities, depending upon which of the three isotropics is zero.  The $D_{ij}(1,0)$ are the first members of the families $D_{ij}(3k+1,0)$.
\item
One copy of $\fieldk^*\oplus\fieldk\supseteq 0\oplus 0$ and two copies of $0\oplus\fieldk\supseteq 0\oplus\fieldk$.  (One isotropic is zero, and the other two are identical lines.)
This is the symplectification of one copy of  $\fieldk\supseteq 0$ and two copies of $\fieldk\supseteq\fieldk$. In the DF classification, this is  $C_r(1,2)$ (or its dual $C_r(1,-2)$) for $r=1$, $2$, or $3$, leading again to three possibilities, depending upon which isotropic is zero.  The $C_r(1,\pm 2)$ are the first members of the families $C_r(3k+1,\pm 2)$.
\item 
Three copies of $0\oplus \fieldk \supseteq 0\oplus\fieldk$.  (All three isotropics are the same line.)  This is the symplectification of three copies of $\fieldk\supseteq \fieldk$ (or its dual $0\supseteq 0$).  
In the DF classification, this is $A(1,3)$ (or $A(1,-3))$. The corresponding dimension vector is  $(1;1,1;1,1;1,1)$ (or $(1;0,0;0,0;0,0)$).
The $A(1,\pm 3)$ are the first members of the families
$A(3k+1,\pm 3)$.  
\item
Two copies of $0\oplus \fieldk$  $\supseteq 0\oplus \fieldk$  and one copy of $\fieldk^*\oplus 0$  $\supseteq \fieldk^*\oplus 0$.
(Two of the isotropics are the same line, and the third one is a different line.)  This is the symplectification of two copies of $\fieldk\supseteq \fieldk$ and one copy of $0\supseteq 0$ (or vice versa). In the DF classification, this is $C_r(1,1)$ (or its dual $C_r(1,-1)$), for $r=1$, $2$, or $3$.
Again we have three possibilities, depending upon which line is distinct from the other two.  The $C_r(1,\pm 1)$ are the first members of the families $C_r(3k+1,\pm 1)$.
\item
The final case, that of three distinct lines in a plane, is linearly indecomposable.  In the DF classification, as a linear representation, it is $A(2,0)$, which is self-dual.  Its dimension vector is $(2;1,1;1,1;1,1)$.  The number of isomorphism classes of symplectic representations with this underlying linear representation depends on the ground field $\fieldk$.  We may find a symplectic basis 
$(e_1,f_1)$ whose elements span $I_1$ and $I_2$, respectively. $I_3$ is then spanned by $e_1 + af_1$ for some nonzero $a\in \fieldk$.  If we change the basis to $be_1,b^{-1}f_1$, then $I_3$ is spanned by $be_1 + baf_1 = be_1 +b^2 a(b^{-1}f_1)$.  This implies that the set of isomorphism classes of triples of lines may be parametrized (taking the case $a=1$ as ``basepoint'') by the square class group $\fieldk^\times/{\fieldk^\times}^2$ introduced in Remark \ref{squares}.  Thus, when $\fieldk$ is algebraically closed, there is just one isomorphism class of this type, while in the case $\fieldk = \reals$, there are two.\footnote{For information about the square class group of other fields, we refer to \cite{becher}, \cite{lam-book}, and \cite{rajwade}.  For instance, the square class group of a finite field has order 1 or 2 according to whether the characteristic is even (\text{i.e.} 2) or odd. For the p-adic numbers, the order of the square class group is 8 for p = 2, and 4 otherwise.}  In this case, the isomorphism class is invariant only under cyclic permutations of the three lines, and the {\bf Maslov index} of the triple distinguishes the two possibilities\footnote{Here we are referring to the Maslov index for Lagrangian triples, also known as the Kashiwara index, see  \cite{duist:morse-index}, \cite{li-ve:weil-maslov}.}. 
$A(2,0)$ is the first member of the family $A(3k+2,0)$ whose members for $k$ even admit compatible symplectic structures.  Those for odd $k$ require symplectification.
\end{enumerate}

We conclude that the number of isomorphism classes of isotropic triples in dimension $2$ is $1+3+3+3+1+3+\#(\fieldk^\times/{\fieldk^\times}^2)$, or $14+\#(\fieldk^\times/{\fieldk^\times}^2)$, where the last term is the order of the  square class group.

\subsection{Higher dimensions: A preview}\label{higher dim}
All of the isotropic triples in dimension $2$ were listed in the previous subsection.  Again, they are the symplectifications, for the case $k=0$, of $A(3k+1,0)$,
$B_r(3k+1,\pm 1)$, $D_{ij}(3k+1,0)$, $C_r(3k+1,\pm 2)$, $A(3k+1,\pm 3)$, $C_r(3k+1,\pm 1)$, along with the isotropic triples which arise from compatible forms for the self-dual sextuple $A(2,0)$. 

For higher $k$, the members of the family $A(3k+2,0)$ are always self-dual (see Section \ref{section on discrete non-split} below), and they admit symplectic forms if only if $k$ is even.  Similarly, the members of the family $A(3k+1,0)$ are all self-dual and admit symplectic forms if and only if $k$ is odd. Thus, for $k=1$ or $k=2$ one finds non-split indecomposable isotropic triples arising from compatible forms for $A(3+1,0)$ and $A(6+2,0)$, respectively. These are in ambient dimension $4$ and $8$, respectively. The symplectifications of $A(3+2,0)$ and $A(6+1,0)$ give indecomposable isotropic triples in ambient dimension $10$ and $16$, respectively.

In dimension $4$, beside the non-split isotropic triples associated with $A(3+1,0)$, we have, for $k=0$, the symplectifications of $A(3k+2,\pm 3)$, $C_r(3k+2,\pm 2)$, $B_r(3k+2,\pm 1)$, $C_r(3k+2,\pm 1)$, and $D_{ij}(3k+2,0)$. 

In dimension $6$, we have, for $k=0$, the symplectifications of the discrete sextuples $C_r(3k+3,\pm 2)$, $B_r(3k+1, \pm 1)$, and $C_r(3k+3,\pm 1)$, and for $k=1$ the symplectifications of the non self-dual continuous-type indecomposable sextuples. In addition,  for $k=2$ there are the non-split isotropic triples arising from compatible symplectic forms for the self-dual continuous-type sextuples.  As we will see later in this paper, non-split continuous-type isotropic triples pose the most intricate case to study, and the number of such isotropic triples depends in particular also on the ground field. 

We end our preview of the ``higher landscape'' of isotropic triples with the remark that indecomposable isotropic triples exist in every even dimension. This follows from the fact that there exist indecomposable non-self-dual sextuples in every given dimension; their symplectifications therefore give indecomposable isotropic triples in every even dimension.

\subsection{Hamiltonian vector fields}\label{Hamiltonian vector fields}

In this section we outline briefly how another problem in linear symplectic geometry can be treated using symplectic poset representations of $\textbf{1}$+$\textbf{1}$+$\textbf{1}$+$\textbf{1}$, and how this problem can in turn be encoded in isotropic triples. As a result, we will get our first example of isotropic triple of continuous type. 

Set $P = \textbf{1}$+$\textbf{1}$+$\textbf{1}$+$\textbf{1}$ and consider the involution $\perp$ on $P$ which exchanges the first two elements and leaves the last two elements fixed. This involution is trivially order-reversing, since all elements of $P$ are incomparable. Symplectic poset representations of $(P, \perp)$ are then subspace systems $(V; U_1, U_2, U_3, U_4)$ where $V$ is symplectic, $U_1$ and $U_2$ are mutually orthogonal, and $U_3$ and $U_4$ are lagrangian.

The problem we will discuss is that of classifying linear hamiltonian vector fields up to conjugation by linear symplectomorphisms; in other words, the classification of the orbits of the Lie algebra $\text{sp} (V, \omega)$ under the adjoint action of the symplectic group $\text{Sp}(V, \omega)$. This is a problem whose solution is well-known and has a long history (going back to Williamson \cite{Williamson} in the 1930s). It has since been treated by various authors, in particular with a view toward finding special normal forms adapted to applications; see, for example, \cite{kocak}, \cite{laub-meyer}. 

Let $(V, \omega)$ be a symplectic vector space. A linear hamiltonian vector field $X$ on $V$ is an element of the Lie algebra $\text{sp}(V, \omega)$; \text{i.e.} it is a linear map $X: V \rightarrow V$ such that $\tilde \omega \circ X = - X^* \circ \tilde \omega$. One wishes to understand equivalence classes, where one linear hamiltonian vector field $(V_1, \omega_1, X_1)$ is equivalent to another, $(V_2, \omega_2, X_2)$, if there exists a linear symplectomorphism $\phi : V_ 1 \rightarrow V_2$ such that $X_2 \circ \phi = \phi \circ X_1$. There is a natural notion of direct sum, and any linear hamiltonian vector field is the direct some of indecomposable pieces. We wish, here, to point out how one may view linear hamiltonian vector fields as symplectic poset representations of $(P, \perp)$. For this we proceed in two steps: first, in the following lemma, we reformulate hamiltonian vector fields in terms of certain kinds of subspaces.

\begin{lemma}\label{reformulation ham}
There is a bijective correspondence between linear hamiltonian vector fields  $(V, \omega, X)$ and linear maps $f :V \rightarrow V^*$ such that $\text{graph}(f) \subseteq V \times V^*$ is a symplectic subspace with respect to the canonical symplectic form on $V \times V^*$.
\end{lemma}

\pf
We give only a sketch. Given $(V, \omega, X)$, it is readily checked that the graph of $f_X := (\tilde \omega \circ X) + \tilde \omega$ is a symplectic subspace. 

Conversely, if $f: V \rightarrow V^*$ is a linear map whose graph is a symplectic subspace, then the asymmetric part $f_a$ will be invertible, and hence defines a symplectic structure $\omega_f$ on $V$.  Setting $X_f := f_a^{-1} f_s$, one finds that $X_f$ is hamiltonian with respect to $\omega_f$.

It is straightforward to check that the two operations are mutually inverse to one another. \qed \\

\begin{cor} A linear hamiltonian vector field $(V, \omega, X)$ can be encoded in the symplectic representation of $(P, \perp)$
$$
(V \times V^*;  \text{graph}(f_X), \text{graph}(f_X)^\perp, V \times 0, 0 \times V^*).
$$
\end{cor}

Although we do not show it here, the passage from a linear hamiltonian vector field to the associated symplectic representation of $(P, \perp)$ is functorial and compatible with the respective notions of direct sum. 

Next we show how the objects above can be encoded in isotropic triples. Observe that the symplectic poset representation of $(P, \perp)$ which we associated to a linear hamiltonian vector field is such that the first two subspaces, which are mutually orthogonal, are symplectic subspaces; in other words, they are independent to each other. The last two subspaces, which are lagrangian, are also independent, and all four subspaces have the same dimension. In the following we will consider  only those symplectic poset representations of $(P, \perp)$ which are of this kind. 

Given such a symplectic representation $\varphi = (V; S, S^\perp, L_1, L_2)$, let $\omega_V$ denote the symplectic form on $V$, $\omega_S$ the restriction of $\omega_V$ to the symplectic subspace $S$, and let $\overline{S}$ denote a copy of $S$ equipped with the symplectic form $- \omega_S$. From $\varphi$ we construct the following isotropic triple in the ambient symplectic space $V \times \overline{S}$, with form $\omega := \omega_V  \times - \omega_S$:
\begin{equation}\label{triple from quadruple}
\begin{array}{ll}
I_1 = L_1 \times 0 &  \quad C_1 = L_1 \times \overline{S}	\\
I_2 = L_2 \times 0  &  \quad C_2 = L_2 \times \overline{S}	\\
I_3 = \{(x,x) \mid x \in S \}  &  \quad C_3 = I_3 + (S^\perp \times 0)
\end{array}
\end{equation}
The passage from $\varphi$ to (\ref{triple from quadruple}) is also functorial and compatible with direct sums. Thus, combining the above, we obtain a way of turning any linear hamiltonian vector field into an isotropic triple. 
We will not analyze in full detail here exactly which types of isotropic triples can be built from linear hamiltonian vector fields. The following, though, describes a large class of isotropic triples which can.

\begin{prop}\label{non-exceptional hamiltonian iso triples}
Let $\varphi = (V; C_i, I_i)$ be an isotropic triple such that 
\begin{enumerate}
\item $V = I_1 \oplus I_2 \oplus I_3$, with $\dim I_i = 1/3 \dim V$ for $i=1,2,3$,
\item $I_i + I_j$ is a symplectic subspace for all $i\neq j$. 
\end{enumerate}
Then we can construct from $\varphi$ a symplectic form $f_a$ and a linear hamiltonian vector field $X$ on $I_2$ such that $\varphi$ is isomorphic to the isotropic triple (\ref{triple from quadruple}) obtained from $(I_2, f_a, X)$.
\end{prop}

\begin{rmk}
Since the geometric description of the isotropic triples in this proposition is invariant under permutation of the indices, the choice of $I_2$ for carrying the hamiltonian vector field is arbitrary. We make this particular choice in order to have coherence with certain types of normal forms which we will use later. 
\end{rmk}

\pf
We can write $V$ as the sum $S \oplus S'$ of $I_1 \oplus I_2$ and its symplectic orthogonal.   
$I_1$ and $I_2$ being a lagrangian decomposition of S, we can identify $I_1$ with $I_2^*$ and, hence, $S$ with the ``cotangent bundle" $T^*I_2$.  

Now $I_3$ cannot intersect $S'$, or its sums with $I_1$ and $I_2$ would not be symplectic subspaces.  So $I_3$ is the graph of a map $g: S \leftarrow S'$ which is antipresymplectic, since $I_3$ is isotropic.  This means that this map $g$ pulls back the symplectic form on $S$ to the negative of that on $S'$.  This implies that $g$ is injective and that its image $g(S')$  is a symplectic subspace of $S$ which is independent of $I_1$ and $I_2$.  Thus, $g(S')$ is the graph of a map 
$f: I_2 \rightarrow I_2^*$, which can be considered as a bilinear form on  $I_2$.   Since $g(S')$ is symplectic rather than isotropic, the form is not symmetric; in fact, it has a nondegenerate antisymmetric part $f_a$. Writing $f = f_a + f_s$ as the sum of its antisymmetric and symmetric parts, since $f_a $ is invertible, we can form
the product $f_a^{-1}f_s$, which is a linear hamiltonian vector field on the symplectic space $(I_2, f_a)$. 

To see that the isotropic triple $\varphi$ is isomorphic to the one of the form (\ref{triple from quadruple}) associated to $(I_2, f_a, X)$, observe that $g$ defines a symplectomorphism $\overline S \leftarrow S'$. It then easy to check that the direct sum of $g$ with the `identity map' on $S = I_1 \oplus I_2$ defines a symplectomorphism from $\varphi$ to (\ref{triple from quadruple}). 
\qed \\

Let us look at an example, to see that isotropic triples of the kind in Proposition \ref{non-exceptional hamiltonian iso triples} do exist. In fact the indecomposable ones come in families dependent on a parameter taking a continuum of values. Indeed, the isomorphism class of the isotropic triple (\ref{triple from quadruple}) built from a linear hamiltonian vector field $X$ depends on $X$ up to conjugation of $X$ by symplectomorphisms, so the spectrum of $X$ is an invariant of the isotropic triple. 

\begin{ex}
Let $V = \mathbb{R}^6$ with symplectic basis $(f_1,f_2,f_3,e_1,e_2,e_3)$. Set
$I_1 = <f_3,e_1>, I_2 = <f_1,e_3>$, and $I_3 = <-\lambda f_1 + (\lambda -1)f_2+f_3, e_1+_2+e_3>$, and
let $\lambda$ vary. Then some computation shows that $X = \tilde \omega^{-1}\tilde \sigma$ as above is given by the matrix whose two diagonal entries are  $1-\lambda$ and $1+\lambda$, so the associated triples for different values of $\lambda$ are non-isomorphic and give a nontrivial 1-parameter family.
\end{ex}

\begin{rmk}
In ambient dimension $6$, it is easy to see directly that isotropic triples of the kind in Proposition \ref{non-exceptional hamiltonian iso triples} are symplectically indecomposable. If there were a decomposition, it would be an orthogonal splitting of the form $\fieldk^2\oplus \fieldk^4$.  Looking
at the possible ways in which each of the isotropics decomposes, it is not hard to check that the induced 2 form on one of the sums must have rank 2 rather than 4, a contradiction.
\end{rmk}

\section{Discrete non-split isotropic triples}\label{section on discrete non-split}

Having given an overview of some background material in the representation theory of posets and quivers, and its connection to isotropic triples, we begin now with the details of our classification. In this section, we study those sextuples which are self-dual and of discrete type, and how they may give rise to non-split isotropic triples. Discrete-type sextuples are somewhat simpler to study than the continuous-type sextuples; the latter are studied in the subsequent, remaining sections of this paper (except for the last section). 

Recall that the indecomposable discrete-type sextuples are uniquely characterized by their dimension vector. In particular, self-dual indecomposable discrete-type sextuples are precisely those whose dimension vector is self-dual, which means here that  $c_j+i_j=v$ for $j=1,2,3$\footnote{See Lemma \ref{dualdimension}.}. Thus, such self-dual sextuples may be read off from the classification in \cite{do-fr:representations}. As stated already in Proposition \ref{DFdiscrete}, the discrete indecomposable sextuples with self-dual dimension vector are denoted 
$A(3k+i,0)$, for $k \in Z_+$ and $i$ equal to $1$ or $2$.  The dimension vectors are   
\begin{align*}
& (3k+1;2k+1,k;2k+1,k;2k+1,k) \quad \quad & \text{for} \ A(3k+1,0), \text{ and} \\
&(3k+2;2k+1,k+1;2k+1,k+1;2k+1,k+1)  \quad \quad &\text{for} \ A(3k+2,0).
\end{align*}
By Lemmas \ref{dualdimension} and \ref{selfdualcompatibleform}, each of these representations admits a compatible form. The degree of uniqueness of such forms is specified in Theorem \ref{unidis} below. Of course, for the symplectic case, non-split istropic triples can only arise in cases of $k$ odd for $A(3k+1,0)$ and $k$ even for $A(3k+2,0)$, since only then is the ambient vector space $V$ even-dimensional.  In fact, whenever $V$ \emph{is} even-dimensional, the compatible forms granted by Lemma \ref{selfdualcompatibleform}
 {\em are} symplectic; this follows from Theorem \ref{dis} and Theorem \ref{unidis}. The self-dual discrete-type sextuples having odd ambient dimension, on the other hand, lead to isotropic triples via their symplectification.

\subsection{Small dimensions}

We give here brief geometric descriptions of the lowest-dimensional non-split discrete-type isotropic triples. These follow, for example, from the normal forms given in the subsections below.

In the previous section, we already saw the first example of a non-split isotropic triple: the underlying sextuple is of type $A(2,0)$ (it belongs to the $A(3k+2,0)$-family), consisting of three distinct lines in a plane.

Next, the non-split isotropic triple arising from the sextuple $A(4,0)$, which belongs to the $A(3k+1,0)$-family. Here, the dimension vector is $(4;3,1;3,1;3,1)$. The three smaller subspaces (in this case they are lines) are independent, and for each line, the 3-dimensional subspace containing it is independent of the other two lines. With a compatible symplectic form, the isotropic triple we obtain has the following form: the isotropic subspaces are three lines $I_i$ in a 4-dimensional space $V$, and each of their orthogonal subspaces $C_i$ is independent from the other two isotropics. Furthermore, $(C_i \cap C_j ) \cap (I_i + I_j) =0$ for any $i\neq j$, so any two of the isotropics span a symplectic subspace. Since the three isotropics are independent, their sum $C$ has codimension 1 and is hence coisotropic. Its orthogonal $C^\perp = C_1 \cap C_2 \cap C_3 \subseteq I_1 + I_2 + I_3$ is a line which must be pairwise independent with each of the $I_i$: if not, then $C^\perp=I_i$ would be the case for some $i$, and hence $C = C_i$ and so $C_i$ would contain all the istropics, a contradiction. Thus the $I_i$ and $C^\perp$ are four lines in general position in $V$. Symplectic reduction via the coisotropic $C$ gives a non-split isotropic triple of the type $A(2,0)$.

  

Moving on, the next case is the sextuple $A(8,0)$, which is in the $A(3k+2,0)$-family. We find the dimension vector to be $(8;5,3;5,3;5,3)$, so the corresponding isotropics are a triple of 3-spaces $I_i$ in an 8 dimensional symplectic space $V$. Though the isotropics $I_i$ are pairwise independent, they are not fully independent: for each distinct triple of indices, $Q_i := I_i \cap (I_j + I_k)$ is a line. The three lines $Q_1, Q_2, Q_3$ are themselves pairwise-independent, and span the 2-dimensional space
$$
I = (I_1 + I_2) \cap (I_2 + I_3) \cap (I_3 + I_1). 
$$ 
This space $I$ is contained in its orthogonal, 
$$
C = (C_1 \cap C_2) + (C_2 \cap C_3) + (C_3 \cap C_1)
$$
and is hence isotropic. Symplectic reduction by the coisotropic $C$ gives a non-split isotropic triple of the type $A(4,0)$. 

We consider one more case. The underlying sextuple is of type $A(10,0)$, which is in the $A(3k+1,0)$-family. The dimension vector is $(10; 7, 3; 7, 3; 7, 3)$. Thus we are again dealing with 3-dimensional isotropics, but this time in a 10-dimensional ambient space. As in the 4-dimensional example above, the isotropics here are completely independent, and their sum is a codimension 1 subspace $C$ which is therefore coisotropic. Its orthogonal, 
$$
C^\perp = (I_1 + I_2 + I_3)^\perp = C_1 \cap C_2 \cap C_3
$$
is a line which is pairwise independent with each of the $I_i$. Symplectic reduction via $C$ gives a non-split isotropic triple of the type $A(8,0)$.

\subsection{Implementing the $A(3k+ 1,0)$}\label{impl A(3k+1,0)}

Given $k\geq 0$ and  a
 basis $\beta = (e_1, \ldots e_{k+1},f_1, \ldots ,f_k,g_1,\ldots ,g_k)$
of a  vector space $V^\beta$,
a sextuple 
$(V^\beta;C_i^\beta,I_i^\beta )$
of type  $A(3k+1,0)$ is given by
\begin{equation}\label{A(3k+1,0)}
\begin{array}{ll}
I^\beta_1= \langle  e_1-f_1, ...,  e_k - f_k \rangle,  \\
I^\beta_2= \langle  g_1, ..., g_k \rangle,  \\
I^\beta_3=\langle  f_1-g_1, ...,  f_k - g_k \rangle,   \\
 \\
C^\beta_1 = \langle e_1,..., e_{k+1}, f_1,..., f_k \rangle  \\
C^\beta_2 = \langle e_1,..., e_{k+1}, g_1,..., g_k \rangle  \\
C^\beta_3 = \langle e_1, e_2 - f_1, ..., e_{k+1} - f_k, f_1 - g_1, ..., f_k - g_k  \rangle.
\end{array}
\end{equation}
Note that 
\[
\begin{array}{ll}
C^\beta_1  = I_1^\beta+ \langle e_1,..., e_{k+1} \rangle & \\
C^\beta_2 = I^\beta_2 + \langle e_1,..., e_{k+1} \rangle & \\
C^\beta_3= I^\beta_3 + \langle e_1, e_2 - f_1, ..., e_{k+1} - f_k \rangle = I^\beta_3 + \langle e_1, e_2 - g_1, ..., e_{k+1} - g_k \rangle. 
\end{array}
\]

In view of Remark \ref{discrete}, to show that this really does define an isotropic triple of type $A(3k+ 1,0)$ it suffices
to observe that the above sextuple has the required dimension vector
 \[
 \quad \dim V^\beta=3k+1,\quad \dim I_i^\beta=k,\quad \dim C_i^\beta=2k+1
 \]
and that its endomorphism algebra is local.
The latter follows from Remark \ref{loc} and the following. 

\begin{lemma}\label{endo alg A(3k+1,0)}
Let $\psi$ be an indecomposable sextuple of type $A(3k+1,0)$. Then 
$$\text{End}(\psi) \simeq  \text{End}((U, \eta))$$
where $\eta$ is an indecomposable nilpotent endomorphism and $\dim U = k+1$. In particular, $\text{End}(\psi)$ is local and $\text{End}(\psi) = \fieldk \id \oplus \text{Rad}$. 
\end{lemma}

\pf
Let $\psi$ be given in the normal form (\ref{A(3k+1,0)}). 
Consider the $\text{End}(\psi)$-invariant subspace
$$
U := C_1^\beta \cap C_2^\beta = \langle e_1, ..., e_{k+1} \rangle
$$
and the indecomposable nilpotent endomorphism $\eta$ of $U$ defined by $\eta(e_1) = 0$ and $\eta(e_{i+1}) = e_i$, for $i=2,..., k$. The endomorphism algebra of $(U, \eta)$, \text{i.e.} the algebra of endomorphisms of $U$ which commute with $\eta$, is local because $\eta$ is indecomposable\footnote{See Remark \ref{loc}. There it is also noted that this endomorphisms algebra $E$ is such that $E = \fieldk \id \oplus \text{Rad} E$.}. We will see now that this algebra is isomorphic to the endomorphism algebra of the sextuple (\ref{A(3k+1,0)}), hence the latter is local as well. To do this, we use the following map. Given an endomorphism $a$ of $U$ which commutes with $\eta$, we can extend it to an endomorphism $\overline a$ of (\ref{A(3k+1,0)}) by defining linear isomorphisms
\[
\begin{array}{ll}
f: \langle e_1,..., e_k \rangle \to \langle f_1,...,f_k \rangle, \qquad f(e_i) := f_i \quad \forall i = 1, ... , k, \\
g: \langle f_1,..., f_k \rangle \to \langle  g_1, ..., g_k \rangle, \qquad g(f_i) := g_i \quad \forall i = 1, ..., k,
\end{array}
\]
and setting
\begin{equation}\label{extension conditions}
\begin{array}{ll}
{\bar a} := f a f^{-1} \qquad  \text{ on }  \langle f_1,..., f_k \rangle, \\
{\bar a} :=g {\bar a} g^{-1}     \qquad  \text{ on }  \langle g_1,..., g_k \rangle.
\end{array}
\end{equation}
(Note that domain of $f$ is invariant under $a$.)

Since $\overline a$ is defined via its action on the basis $\beta$, it is easily checked directly that $\overline a$ is an endomorphism of (\ref{A(3k+1,0)}). To see this, note in particular that $I_1^\beta=\{x- f(x)\mid x  \in \langle f_1,..., f_k \rangle \}$ and $I_3^\beta =\{x - g(x)\mid x \in \langle g_1,..., g_k \rangle \}$. 

The map $a \mapsto \overline a$ has as its inverse the operation of taking an endomorphism $b$ of (\ref{A(3k+1,0)}) and restricting it to $U = C_1^\beta \cap C_2^\beta$ (which will necessarily be an invariant subspace of $b$, since by assumption $C_1^\beta$ and  $C_2^\beta$ are $b$-invariant). To see this, notice that such a $b$ necessarily decomposes as the direct sum of its restrictions to the subspaces
$$
U = \langle e_1, ..., e_{k+1} \rangle, \quad F := \langle f_1,...,f_k \rangle, \quad G:= \langle g_1,...,g_k \rangle
$$
since these subspaces sum to $V$ and must be invariant under $b$:
$$
U = C^\beta_1 \cap C^\beta_2, \quad F = I^\beta_2, \quad G = C^\beta_1 \cap (I^\beta_1 + I^\beta_2).
$$
The invariance of $I^\beta_1$ and $I^\beta_3$ under $b$ then enforces that $b_{\vert_{U}}$ is related to $b_{\vert_{F}}$ and $b_{\vert_{G}}$ via (\ref{extension conditions}), and together with the invariance of $C^\beta_3 \cap C^\beta_1$ it is ensured that $b_{\vert_{U}}$ commutes with $\eta$. 

Thus if we restrict $b$ to $U$ and then extend to $\overline b$, we recover $b$. Conversely, if we start with an endomorphism $a$ of $U$, the restriction of $\overline a$ to $U$ is of course, by definition, again $a$. 

Finally, it is clear from (\ref{extension conditions}) that the operation $a \mapsto \overline a$ is a morphism of algebras. \qed

\subsection{Implementing the $A(3k+2,0)$}\label{impl A(3k+2,0)}
We will in fact work with
$A(3(k-1)+2,0)=A(3k-1,0)$, which
is a subquotient of $A(3k+1,0)$ for $k>0$, so that  the construction
of compatible forms for sextuples of both types $A(3k+1,0)$ and $A(3k+2,0)$ can be treated uniformly. 

\medskip

Given $k\geq 1$ and a basis $\gamma = (e_2, \ldots e_{k},f_1, \ldots ,f_k,g_1,\ldots ,g_k)$
of a  vector space $V^\gamma$, a sextuple 
$(V^\gamma;C_i^\gamma,I_i^\gamma)$
of type  $A(3k-1,0)$ is given by
\begin{equation}\label{A(3k-1,0)}
\begin{array}{ll}
I^\gamma_1 = \langle  f_1, f_2 - e_2, ..., f_k - e_k \rangle \\
I^\gamma_2 = \langle g_1, ..., g_k \rangle, \\
I^\gamma_3 = \langle f_1 - g_1, ..., f_k - g_k \rangle, \\
\\
C^\gamma_1=  \langle e_2, ..., e_k, f_1,..., f_k  \rangle \\
C^\gamma_2 = \langle e_2, ..., e_k, g_1, ..., g_k \rangle \\
C^\gamma_3= \langle f_1 - e_2, ..., f_{k-1} - e_k, f_1 - g_1, ..., f_k - g_k \rangle.
\end{array}
\end{equation}
Note that 
\[
\begin{array}{l}
C_1^\gamma  = I_1^\gamma + \langle e_2, ..., e_k \rangle= I_1^\gamma+ \langle f_2, ..., f_k \rangle  \\
C^\gamma_2 = I^\gamma_2 + \langle e_2, ..., e_k \rangle  \\
C^\gamma_3= I^\gamma_3 + \langle f_1 - e_2, ..., f_{k-1} - e_k \rangle.  
\end{array}
\]

Again, in view of Remark \ref{discrete}, it suffices
to observe that such a sextuple has the required dimension vector
  \[
\dim V^\gamma=3k-1,\quad \dim I^\gamma_i=k, \quad \dim C^\gamma_i=2k-1\]
and local endomorphism algebra.  
For the latter, we proceed similarly as for sextuples of type $A(3k+1,0)$, combining Remark \ref{loc} and the following.  

\begin{lemma}
Let $\psi$ be an indecomposable sextuple of type $A(3k-1,0)$. Then 
$$\text{End}(\psi) \simeq  \text{End}((U, \eta))$$
where $\eta$ is an indecomposable nilpotent endomorphism and $\dim U = k$. In particular, $\text{End}(\psi)$ is local and $\text{End}(\psi) = \fieldk \id \oplus \text{Rad}$. 
\end{lemma}

\pf
Let $\psi$ be given in the normal form (\ref{A(3k-1,0)}). Consider the $\text{End}(\psi)$-invariant subspace
$$
U := \langle f_1,..., f_k  \rangle = C_1^\gamma\cap (I_2^\gamma +I_3^\gamma)
$$ 
and the indecomposable nilpotent endomorphism $\eta$ of $U$ defined by $\eta(f_1) = 0$ and $\eta(f_{i+1}) = f_i$, for $i=2,..., k$.
Similarly as in the previous section, the algebra of endomorphisms of $U$ which commute with $\eta$ is isomorphic to the endomorphism algebra of (\ref{A(3k-1,0)}). 

To show this, we begin with an endomorphism $a$ of $U$ which commutes with $\eta$ and extend it to an endomorphism $\overline a$ of $V^\gamma$ by defining maps $f$, $g$, and $h$ by
\[
\begin{array}{ll}
f : \langle e_2, ..., e_k \rangle \to \langle f_1,..., f_k \rangle, \quad f(e_i) = f_i \ \text{ for } i =2, ..,k, \\
g: \langle f_1,..., f_k \rangle \to \langle  g_1, ..., g_k \rangle, \quad g(f_i) = g_i \ \text{ for } i =1, ..,k, \\
h: \langle f_1, ..., f_k \rangle \to \langle e_2,..., e_k \rangle, \quad h(f_1) = 0, \ \  h(f_{i}) = e_i  \ \text{ for } i =2, ..,k.
\end{array}
\]
and setting
\begin{equation}\label{relations for extension}
\begin{array}{ll}
\overline a  := h a f \qquad \text{ on } \langle e_2, ..., e_k \rangle. \\
\overline a := g a g^{-1} \ \quad  \text{ on } \langle  g_1, ..., g_k \rangle.
\end{array}
\end{equation}

Note that $I_1^\gamma=\{x - h(x)\mid x  \in \langle f_1, ..., f_k \rangle  \}$ and $I_3^\gamma =\{x - g(x)\mid x \in  \langle f_1, ..., f_k \rangle \}$, and that 
$$
E = \langle e_2, ..., e_{k} \rangle, \quad U = \langle f_1,...,f_k \rangle, \quad G = \langle g_1,...,g_k \rangle
$$
are subspaces invariant under the endomorphism algebra of (\ref{A(3k-1,0)}) since 
$$
E = C_1^\gamma \cap C_2^\gamma, \quad U = C_1^\gamma \cap (I_2^\gamma + I_3^\gamma), \quad G = I_2^\gamma. 
$$

For an endomorphism $b$ of (\ref{A(3k-1,0)}), the operation $b \mapsto b_{\vert_{F^\gamma}}$ is inverse to the operation $a \mapsto \overline a$. Indeed, such a $b$ decomposes as the direct sum $b =  b_{\vert_{E}} \oplus  b_{\vert_{U}} \oplus b_{\vert_{G}}$, and the relations (\ref{relations for extension}) are enforced by the invariance of $I_1^\gamma$ and $I_3^\gamma$ under $b$.

Finally, that the map $a  \mapsto \overline a$ is a morphism of algebras is evident from (\ref{relations for extension}).
\qed

\

The bases $\beta$ and $\gamma$ are referred to as {\bf standard bases}
for the respective types of sextuple. For each fixed $k > 0$ one can view (\ref{A(3k-1,0)}) as a subquotient of (\ref{A(3k+1,0)}) by considering the subspaces $I^\beta \subseteq C^\beta \subseteq V^\beta$, 
$$
I^\beta = \langle e_1 \rangle, \quad C^\beta = \langle e_1,...,e_k, f_1,..., f_k, g_1,..., g_k \rangle, 
$$
and letting $V^\gamma = C^\beta/I^\beta$, where the basis $\gamma$ is induced (modulo re-indexing) by the elements of the basis $\beta$ which span $C$ (i.e. all elements except $e_{k+1}$). Then we have the identifications 
$$
I_i^\gamma = (I_i^\beta \cap C^\beta)/I^\beta \quad \text{and} \quad C_i^\gamma = (C_i^\beta \cap C^\beta)/I^\beta \quad \text{for} \ i = 1,2,3. 
$$

In fact, we can also view a sextuple of type $A(3(k-2)+1,0) = A(3k - 5, 0)$ as a subquotient of a sextuple of type $A(3k-1,0)$. If a sextuple of the latter type is given in the form $(\ref{A(3k-1,0)})$, then choosing 
$$
I^\gamma = \langle f_1, g_1 \rangle, \quad C^\gamma = \langle e_2,...,e_k, f_1,..., f_{k-1}, f_1,..., g_{k-1} \rangle, 
$$
and passing to the subquotient $C^\gamma /I^\gamma$ gives a sextuple of type $A(3(k-2)+1,0)$, with standard basis induced from the standard basis of $(\ref{A(3k-1,0)})$.

In total, the (isomorphism classes of) discrete non-split sextuples form two ``chains of subquotients'' (arrows indicate the passage to a subquotient):
\[
\begin{tikzcd}
			& k=0	& k=1			& k=2			& k=3			& \dots\\
A(3k+1,0) 		& A(1,0)	& A(4,0) \arrow[dl]	& A(7,0) \arrow[dl]	& A(10,0) \arrow[dl]	& \dots \arrow[dl] \\
A(3k+2,0)		& A(2,0)	& A(5,0) \arrow[ul]	& A(8,0) \arrow[ul]	& A(11,0) \arrow[ul]	& \dots  \arrow[ul]
\end{tikzcd}
\]

\subsection{Existence of compatible forms}\label{existence adm forms}
\newcommand{\vep}{\varepsilon}

\begin{thm}\label{dis} Any sextuple of type  $A(3k\pm1,0)$ admits $\vep$-symmetric forms  where  $\vep=(-1)^k$.
With respect to standard bases there is a recursion providing  compatible $\vep$-symmetric forms  with coefficients in the prime subfield. 
\end{thm}

The existence will be shown by constructing matrices for these compatible forms
with respect to standard bases. These matrices will be of the shape $(\ref{comp block form})$, \text{i.e.} having a block structure and built using a (smaller) matrix which we call ``$A$''. 

Fix $k\in \mathbb{N}$, let $\vep=(-1)^k$, and let $A=A_k\in \fieldk^{(k+1)\times (k+1)}$. 
We consider the following relations on the entries of $A$:
\begin{enumerate}[(1)]
\item \label{1} $a_{ij}=0$ \ for $i+j<k+2$, \quad $a_{ij}\neq 0$ \ for $i+j=k+2$
\item \label{2} $A=\vep A^t$
\item \label{3}  $a_{i,j+1}=a_{ij}-a_{i+1,j}$ \ for $i,j=1,\ldots,k$ 
\item \label{4}  $a_{kk}=a_{k,k+1}+a_{k+1,k}$.
\end{enumerate}
Further we define $\mathring A =\mathring A_k$ as the minor of $A$ given by restricting to row and column indices in $\{ 2,\ldots,k \}$ and we set 
\[
c_i=a_{i,k+1},\quad c'_i=a_{k+1,i}. 
\]
If (\ref{1}) holds, then the matrix $A$ has the following form
\[A=
\left[\begin{array}{c|cccccc|c}
0&0&0&\cdots&\cdots&0&0&a_{1,k+1}\\ \hline
0&0&0&\cdots&\cdots&0&a_{2,k}&a_{2,k+1}\\
0&0&0&\cdots&\cdots&a_{3,k-1}&a_{3,k}&a_{3,k+1}\\
\vdots&\vdots&\vdots&&&\vdots&\vdots&\vdots\\
\vdots&\vdots&0&&&\vdots&\vdots&\vdots\\
0&0&a_{k-1,3}&\cdots&\cdots&a_{k-1,k-1}&a_{k-1,k}&a_{k-1,k+1}\\
0&a_{k2}&a_{k3}&\cdots&\cdots&a_{k,k-1}&a_{kk}&a_{k,k+1}\\\hline
a_{k+1,1}&a_{k+1,2}&a_{k+1,3}&\cdots&\cdots&a_{k+1,,k-1}&a_{k+1,k}&a_{k+1,k+1}
 \end{array}\right].
\]
The vertical and horizontal lines inside of this matrix are intended solely as visual aids. 

Now we define
\begin{equation}\label{comp block form}
H=H_k:=
\left[\begin{array}{c|ccc|c||c|ccc||ccc|c}
0&0&\cdots&0&c_1&0&0&\cdots&0&0&\cdots&0&0\\ \hline
0&&&& c_2&0&&&&&&&0\\
\vdots&&\mathring A&&\vdots &\vdots&&\mathring A&&&0&&\vdots\\
0&&&& c_k&0&&&&&&&0\\\hline
c'_1&c'_2&\ldots&c'_k& c_{k+1}&c'_1&c'_2&\cdots&c'_k&0&\cdots&0&0 \\ \hline\hline
0&0&\cdots&0&c_1&0&0&\cdots&0&0&\cdots&0&c_1\\ \hline
0&&&& c_2&0&&&&&&&c_2\\
\vdots& &\mathring A& &\vdots &\vdots& &\mathring A& & & \mathring A&&\vdots\\
0&&&& c_k&0&&&&&&&c_k\\ \hline\hline
0&&&&0&0&&&&&&&0\\ 
\vdots&&0&&\vdots&\vdots&&\mathring A&&&0&&\vdots\\
0&&&&0&0&&&&&&&0\\ \hline
0&0&\cdots&0& 0&c'_1&c'_2&\cdots&c'_k&0&\cdots&0&0\\
\end{array}\right]
\end{equation}
and we interpret this matrix to be the coordinate matrix of a bilinear form $B$ on $V^\beta$ (\text{c.f.} Section \ref{impl A(3k+1,0)}) with respect to the standard basis $\beta = \{e_1,..., e_{k+1}, f_1,..., f_k, g_1, ..., g_k\}$. The double lines in the matrix are visual aids for seeing the block structure 
\[
H=
\left[ \begin{array}{ccc}
H_{11} & H_{12} & H_{13} \\
H_{21} & H_{22} & H_{23} \\ 
H_{31} & H_{32} & H_{33}
\end{array}
\right]
\]
related to the subspaces $\langle e_1,..., e_{k+1} \rangle$, $\langle f_1,..., f_k \rangle$, and   $\langle g_1, ..., g_k \rangle$. 

\begin{claim} 
If (\ref{1}) and (\ref{2}) hold, then 
$B$ is non-degenerate and $\vep$-symmetric
\end{claim}
\pf
That $B$ is $\vep$-symmetric follows directly from (\ref{2}). To see that $B$ is non-degenerate, note that
the blocks of the matrix $H$ are such that $H_{33}$, $H_{31}$ and $H_{13}$ are zero and where  
$H_{11}$, $H_{23}$, and $H_{32}$ are square matrices
having non-zero entries on the anti-diagonal, and zeros above the anti-diagonal;
in particular, the latter blocks are invertible. Because $H_{13}$ is zero and because a `copy' of $H_{21}$ is contained in $H_11$, we can use row operations to turn transform $H_{21}$ to zero in such a way that only $H_{22}$ is additionally changed under these operations. In a similarly manner we can also turn $H_{12}$ to zero using column operations. At this point our block matrix has the following form (tilde indicates that there are changes)
\[
\tilde H=
\left[ \begin{array}{ccc}
H_{11} & 0 & 0\\
0 & \tilde H_{22} & H_{23} \\ 
0 & H_{32} & 0
\end{array}
\right]
\]
and where $\tilde H_{22}$ is
\[
\left[\begin{array}{c|ccc}
0&0&\cdots&0\\ \hline
0&&& \\
\vdots&&\mathring A &\\ 
0&&&
\end{array}\right]
-
\left[\begin{array}{c|ccc}
0&0&\cdots&0\\ \hline
0&&& \\
\vdots&&\mathring A &\\ 
0&&&
\end{array}\right]
-
\left[\begin{array}{c|ccc}
0&0&\cdots&0\\ \hline
0&&& \\
\vdots&&\mathring A &\\ 
0&&&
\end{array}\right]
=
- 
\left[\begin{array}{c|ccc}
0&0&\cdots&0\\ \hline
0&&& \\
\vdots&&\mathring A &\\ 
0&&&
\end{array}\right].
\]
Now, clearly we can add columns from the third block column to the second block column to transform $\tilde H_{22}$ into the zero matrix, and this leaves all other blocks unchanged, since $H_{13}$ and $H_{33}$ are zero. (Equivalently, we could have used row transformations using rows from the bottom block row.) This puts our block matrix in the form
\[
\left[ \begin{array}{ccc}
H_{11} & 0 & 0\\
0 & 0 & H_{23} \\ 
0 & H_{32} & 0
\end{array}
\right]
\]
which shows non-degeneracy, since the non-zero blocks are non-degenerate. 
\qed \\

\begin{claim}\label{str} If (\ref{2}) holds, then 
\[
B(I_1^\beta,C_1^\beta)=B(I_2^\beta,C_2^\beta)=B(I_3^\beta, \langle e_1,
e_{2} - g_1,..., e_{k+1} - g_k  \rangle )=0.
\]
\end{claim}

\pf
\begin{itemize}
\item $B(e_i,\;e_j-f_j) =a_{ij}-a_{ij}=0$
for $i=1,\ldots,k+1$  and $j=1,\ldots,k$. 
\item $B(f_i,\;e_j-f_j)=a_{ij}-a_{ij}=0$ for $i,j=1,\ldots,k$. 
\item $B(e_i,g_j)=0$ for $i=1,\ldots,k+1$, $j=1,\ldots,k$.
$B(g_i,g_j)=0$ for  $i,j=1,\ldots,k$.
$B(e_1,\;f_j-g_j)=B(e_1,f_j)-B(e_i,g_j)=0-0=0$ for  $i,j=1,\ldots,k$.
\item $B(f_i-g_i,\;e_{j+1}-g_j)=   
B(f_i,e_{j+1})-B(f_i,g_j) -0+0= a_{i,j+1}-a_{i,j+1}=0$ for $i,j=1,\ldots,k$.
\end{itemize}
\qed \\

\begin{claim}\label{ad} 
If (\ref{1}) through (\ref{4}) hold, 
then
 $B$ is $\vep$-symmetric and a compatible form for the sextuple
$(V^\beta;C_i^\beta,I_i^\beta)$. 
\end{claim}

\pf
It remains to show  $B(I_3^\beta,C_3^\beta)=0$. Consider 
$B(f_i-g_i,\;f_j-g_j)
=B(f_i,f_j)-B(f_i,g_j)- B(g_i,f_j)+B(g_i,g_j)=:x_{ij}$. Direct checking gives:
\begin{itemize}
\item $x_{11}=0+0+0+0=0$. 
\item $x_{1j}=0+0+a_{2j}+0=0$ for $j=2,\ldots,k-1$.
\item $x_{ij}
=  a_{ij}- a_{i,j+1}-a_{i+1,j}+0=0$
for $i,j=2,\ldots,k-1$ by (3).
\item $x_{ik}= a_{ik}- c_i-a_{i+1,k}+0=0$  
for $i=2,\ldots,k-1$ by (3). 
\item $x_{1k}=0-c_1-a_{2k}+0=0$ by (3).
\item $x_{kk}= a_{kk}-c_k-c'_k=0$ by (4).
\end{itemize}
In view of Claim \ref{str} and that 
\begin{align*}
I^\beta_3&=\langle  f_1-g_1, ...,  f_k - g_k \rangle \\
C^\beta_3&= \langle e_1, e_2 - f_1, ..., e_{k+1} - f_k, f_1 - g_1, ..., f_k - g_k  \rangle
\end{align*}
we have compatibility of $B$.
\qed \\

\begin{claim}\label{ex} 
 For each $k = 0,1,2,...$
there exist matrices $A_k$ satisfying (\ref{1})--(\ref{4}).
\end{claim}

\pf We proceed by induction on $k$. 
Let \[A_0=(c_1), \quad A_1= \left(\begin{array}{cc}
0&c_1 \\
-c_1&0
\end{array}\right),\quad c_1\neq 0\]
Assume that $A_{k-2}$ is given satisfying (\ref{1})--(\ref{4}).
We define $A=A_k$ to have minor $\mathring A=A_{k-2}$
with respect to row and column indices $2,\ldots,k$.
Thus, we have 
\begin{enumerate}[(1')]
\item for $2\leq i,j\leq k$: $a_{ij}=0$ if $i+j<k+2$, \ $a_{ij}\neq 0$ if $i+j=k+2$ 
\item $\mathring A=\vep \mathring A^t$
\item   $a_{i,j+1}=a_{ij}-a_{i+1,j}$ for $i,j=2,\ldots,k-1$, 
\item  $a_{k-1,k-1}=a_{k-1,k}+a_{k,k-1}$.
\end{enumerate}
It remains to grant
\begin{enumerate}[(1'')]
\item $a_{1j}=a_{i1}=0$ for $i,j=1,\ldots,k$ and $a_{1,k+1}\neq 0$,
$a_{k+1,1}\neq 0$ 
\item $a_{k+1,i}=\vep a_{i,k+1}$ for $i=1,\ldots ,k+1$ 
\item   $a_{i,j+1}=a_{ij}-a_{i+1,j}$ for $j=k,\,i=1,\ldots$  resp. $i=k,\,j=1,\ldots,k-1$
\item  $a_{kk}=a_{k,k+1}+a_{k+1,k}$.
\end{enumerate}
Define
\[ \begin{array}{rcll}
a_{1,k+1}&=&-a_{2,k}\\ a_{i,k+1}&=&a_{i,k}-a_{i+1,k} &\mbox{ for } i=2,\ldots,k-1\\
a_{k,k+1}&=&\frac{1}{2} a_{k,k}&\mbox{ if $k$ is even}\\
a_{k,k+1}&=& \mbox{arbitrary } &\mbox{ if $k$ is odd} \\
a_{k+1,k+1}&=&
 \mbox{arbitrary } &\mbox{ if $k$ is even}\\
a_{k+1,k+1}&=&0&\mbox{ if $k$ is odd}\\ 
a_{k+1,i}&=& \vep a_{i,k+1} &\mbox{ for } i=1,\ldots ,k. 
\end{array} \]
Then (1'') and (2'') are obvious, as is (3'') for the cases when $j=k$, 
while for those cases with $i=k$ we have $a_{k,j+1}=\vep a_{j+1,k}=\vep(a_{j,k}-a_{j,k+1})=a_{kj}-a_{k+1,j}$. Finally, if $k$ is odd then $a_{kk}=\vep a_{kk}=0
= a_{k,k+1}+\vep a_{k,k+1} =a_{k,k+1}+a_{k+1,k}$; and if $k$ is even
then $a_{kk}= \frac{1}{2}a_{kk}+\frac{1}{2}a_{kk}=a_{k,k+1}+
a_{k+1,k}$. \qed \\

\

\noindent 
\textbf{Proof of Theorem \ref{dis}}.
For the $A(3k+1,0)$, Claims \ref{ad} and \ref{ex} and the
proof of the latter give the recursive construction 
of such forms with respect to standard bases. For $A(3k-1,0)$, we view $V^\gamma$ as a subspace of $V^\beta$ and define on it a form given by a matrix $H' = H_k'$ obtained from $H = H_k$ by omitting all rows and columns indexed by
$e_1$ or $e_{k+1}$. Then $H'$ defines an $\vep$-symmetric
form $B'$ on $V^\gamma$. To see this, let $S^\beta:= \langle e_1, e_{k+1}\rangle$ and note that $H_{\vert_{S^\beta}}$ is non-degenerate and $\vep$-symmetric, and that we can make the identification $V^\gamma = (S^\beta)^\perp \subseteq V^\beta$. Hence $H' = H_{\vert_{(S^\beta)^\perp}}$ is also non-degenerate and $\vep$-symmetric.  It remains to show that $H'$ is compatible with the $A(3k-1,0)$ sextuple in $V^\gamma$ obtained from the $A(3k+1,0)$ sextuple in $V^\beta$. 

Recall that if $X \subseteq V^\beta$ is an element of the sextuple in $V^\beta$, then the corresponding subspace of the sextuple in $V^\gamma$ is given by $\underline{X} := \pi(X^\beta \cap C^\beta)$, where $\pi$ denotes the projection onto the second factor of the (orthogonal) decomposition $V^\beta = S^\beta \oplus V^\gamma$. Since $C^\beta = (C^\beta)^\perp \oplus V^\gamma$ is coisotropic, this is an instance of coisotropic reduction; in particular $$
B'(\underline{X},\underline{Y}) = B(X \cap C^\beta ,Y \cap C^\beta)
$$
for any $X,Y \subseteq V^\beta$. Hence $B(X,Y) = \{0\}$ implies that $B'(\underline{X},\underline{Y})= \{0\}$.  Since we know the dimensions of all the subspaces involved in our sextuples, this shows that for any element $X$ of the sextuple in $V^\beta$, the orthogonal of $\underline{X}$ in $V^\gamma$ is the same subspace as $\underline{X^\perp}$.
\qed \\

\begin{rmk}\label{rem} Suppose we are given sextuples of the types $A(3k+1,0)$ and $A(3k-1,0)$ in terms of standard bases. 
\begin{enumerate}[(i)]
\item With respect to the standard basis, any compatible form on $A(3k+1,0)$
is given by a matrix $H_k$ which is built from a matrix $A$ as above, and such that $A$
satisfies the relations (\ref{1})--(\ref{4}). The analogous statement is true for compatible forms on $A(3k-1,0)$ in terms of matrices of the form $H_k'$, where $H_k'$ is obtained from $H_k$ by deleting the 1st and $(k+1)th$ rows and columns. 
\item A parametrization of 
such matrices $H_k$ resp. $H_k'$ is given 
by the parameters 
$$
a_{k+1,j}, \quad j=1,\ldots,k+1 \quad \quad \text{resp.} \quad \quad a_{k,j}, \quad  j=2,\ldots,k 
$$
Moreover, the matrix entries are obtained via  linear expressions
from the parameters.
\end{enumerate}
\end{rmk}

\pf
(i) According to the proofs of Claims \ref{str} and \ref{ad},
the structure  of $B$ and the relations (\ref{1})--(\ref{4}) are forced by the requirement of admissibility.

(ii) This follows from the recursive construction of the matrices $A_k$; in each step one free parameter can be chosen. 

\qed 

\begin{ex}
Let $k=1$. A sextuple of type $A(3k+1, 0)$ has ambient dimension $4$. To construct a compatible symplectic form (with respect to a standard basis) we begin with the matrix
$$
A_1= \left(\begin{array}{cc}
0&c_1 \\
-c_1&0
\end{array}\right)
$$
where $c_1$ is any non-zero scalar. From this we obtain the matrix (\ref{comp block form}) of a compatible symplectic form; in this example it is
$$
H_1 = \left(\begin{array}{cc|c|c}
0&c_1 & 0 & 0 \\
-c_1&0 & -c_1 & 0 \\
\hline 
0 & c_1 & 0 & c_1 \\
\hline 
0 & 0 & -c_1 & 0
\end{array}\right)
$$
(note that $\mathring  A_1$ is the empty matrix). To obtain a compatible symplectic form for the sextuple $A(3k-1,0)$, with respect to a standard basis induced from $A(3k+1,0)$, we only need to drop the $1st$ and $(k+1)th$ rows and columns of the compatible form $H_1$ given above. This gives
$$
H_1' = \left(\begin{array}{cc}
0& c_1  \\
-c_1&0  
\end{array}\right).
$$
\end{ex}

\begin{ex}
Let $k=3$. A sextuple of type $A(3k+1, 0)$ has ambient dimension $10$. To construct a compatible symplectic form (with respect to a standard basis) we proceed similarly as in the previous example. Following the recursion recipe given in Claim \ref{ex} above, we build from $A_1$ the matrix
$$
A_3= \left(\begin{array}{cccc}
0 & 0 & 0 & c_1 \\
0 & 0 & -c_1 & -c_1 \\
0&c_1 & 0 & c_2 \\
-c_1& c_1 & -c_2 & 0
\end{array}\right)
$$
where $c_2$ is a scalar that we may freely choose. 
Now from $A_3$ we obtain the matrix (\ref{comp block form}) of a compatible symplectic form; in this example it is
$$
H_3 = \left(\begin{array}{cccc|ccc|ccc}
0 & 0 & 0 & c_1 & 0 & 0 & 0 & 0 & 0 & 0 \\
0 & 0 & -c_1 & -c_1  & 0 & 0 & -c_1 & 0 & 0 & 0 \\
0&c_1 & 0 & c_2  & 0&c_1 & 0 & 0 & 0 & 0\\
-c_1& c_1 & -c_2 & 0 & -c_1&c_1 & -c_2 & 0 & 0 & 0 \\
\hline
0 & 0 & 0 & c_1 & 0 & 0 & 0 & 0 & 0 & c_1 \\
0 & 0 & -c_1 & -c_1  & 0 & 0 & c_1 & 0 & -c_1 & -c_1\\
0&c_1 & 0 & c_2  & 0&-c_1 & 0 &c_1 & 0 & c_2\\
\hline
0 & 0 & 0 & 0  & 0 & 0 & -c_1  & 0 & 0 & 0 \\
0 & 0 & 0 & 0  & 0 & c_1 & 0 & 0 & 0 & 0 \\
0 & 0 & 0 & 0  & -c_1& c_1  & -c_2 & 0 & 0 & 0 
\end{array}\right).
$$
To obtain a compatible symplectic form for the sextuple $A(3k-1,0)$, we drop the $1st$ and $(k+1)th$ rows and columns of $H_3$. This gives
$$
H_3' =  \left(\begin{array}{cc|ccc|ccc}
 0 & -c_1   & 0 & 0 & -c_1 & 0 & 0 & 0 \\
c_1 & 0   & 0&c_1 & 0 & 0 & 0 & 0\\
\hline
 0 & 0  & 0 & 0 & 0 & 0 & 0 & c_1 \\
 0 & -c_1   & 0 & 0 & c_1 & 0 & -c_1 & -c_1\\
 c_1 & 0   & 0&-c_1 & 0 &c_1 & 0 & c_2\\
\hline
 0 & 0  & 0 & 0 & -c_1  & 0 & 0 & 0 \\
 0 & 0   & 0 & c_1 & 0 & 0 & 0 & 0 \\
 0 & 0   & -c_1& c_1  & -c_2 & 0 & 0 & 0 
\end{array}\right).
$$
\end{ex}

\subsection{Uniqueness of compatible forms}\label{uniqueness discrete selfdual}

\begin{thm}\label{unidis}
The $(-1)^k$-symmetric forms on an $A(3k\pm1,0)$ are 
 unique up to 
isometric automorphism and multiplication by scalars;
there are  no $(-1)^{k+1}$-symmetric forms
on an $A(3k\pm1,0)$.

If $B$ is a compatible $\vep$-symmetric form on $A(3k\pm1,0)$
and $c\in \fieldk$, then there is an automorphism $\eta$
of the sextuple which is an isometry, in the sense $\eta^*B = cB$, if and only if $c$ is a square in $\fieldk$.
\end{thm}
\pf
Dealing with $A(3k+1,0)$, we continue the discussion
from Subsection \ref{uniqueness of compatible forms}. By Lemma \ref{endo alg A(3k+1,0)}, we may apply Lemma \ref{uniqueness up to a scalar}; this yields the first claim.

Now, suppose there exists an automorphism $\eta$ of $A(3k+1,0)$ such that 
$f^*B = cB$, and let $\eta_U$ to be its restriction to $U= C_1^\beta \cap C_2^\beta$, which is invariant under $\eta$. We may assume that $B$ is given by 
a matrix $B_k$ in terms of a standard basis as above. Thus, $U = \langle e_1, ..., e_{k+1} \rangle$
is non-degenerate and the restriction of $B$ to $U$
has matrix $A$ which is zero above the
anti-diagonal.
By Observation \ref{loc}, $\eta_U$ 
has  upper triangular matrix $(\eta_{ij})$ with
 diagonal entries  all the same. 
It follows that
\begin{align*}
cB(e_{k+1},e_1)&=B(\eta e_{k+1},\eta e_1)  
=B( \sum_{i=1}^{k+1}\eta_{i,k+1}e_i,\eta_{11} e_1) \\
	&= B( \eta_{k+1,k+1}e_{k+1},\eta_{11} e_1)
= B( \eta_{11}e_{k+1},\eta_{11} e_1) = \eta_{11}^2 B(e_{k+1},e_1)  
\end{align*}
whence $c = \eta_{11}^2$. (Conversely, given a square $c = b^2 \in \fieldk$, one can of course always find an isometry between $B$ and $cB$: simply take $b \cdot \text{id}$.)

For $k>1$, a similar reasoning works for $A(3k-1,0)$, with
$U=C_1^\gamma\cap C_2^\gamma =  \langle e_2, ..., e_{k} \rangle$.
\qed \\

\section{Continuous types: sextuple classification and duality}\label{continuous}

By definition, continuous-type indecomposable sextuples are those indecomposable sextuples which have dimension vectors of the form $(3k;2k,k;2k,k;2k,k)$, for integers $k \geq 1$.  Note that this is a self-dual dimension vector.

As mentioned above, and proven in \cite[Section 4.6]{do-fr:representations}, when $\fieldk$ is algebraically closed  the (isomorphism types of) continuous-type indecomposable sextuples in a given dimension (indicated by $k$) consist of:
\begin{itemize}
 \item a collection
$\Delta(k;\lambda)$, labeled by scalars $\lambda \in \fieldk \backslash \{ 0,1\}$, 
\item so-called \textbf{exceptional} continuous-type sextuples, which are labeled
$$
\begin{array}{l}
\Delta_j(k,1) \text{ for } j\in \{1,2\}, \\
\Delta_j(k,0) \text{ for } j\in \{1,2,3\},  \\
\Delta_j(k,\infty), \text{ for } j\in \{1,2,3\}.
\end{array}
$$
\end{itemize}  

In the case of more general fields, the continuous-type indecomposable sextuples $\Delta(k;\lambda)$ must be replaced by indecomposable sextuples $\Delta(k;\gamma)$, where $\gamma$ ranges over all indecomposable endomorphisms of $\fieldk^k$ which do not have $0$ or $1$ as an eigenvalue. The exceptional continuous-type indecomposable sextuples remain the same. Thus the isomorphism classes of continuous-type indecomposable sextuples in a given dimension $3k$ are parametrized by the disjoint union of the sets
$$
\{ \gamma \mid \gamma \text{ indecomposable endomorphism of } \fieldk^k, \text{ with } 0, 1 \notin \text{spec}(\gamma) \} 
$$
and 
$$
 \{0_1, 0_2, 0_3, 1_1, 1_2, \infty_1, \infty_2, \infty_3\},
$$
where the latter set consists of formal labels for the exceptional types. 

We review the classification and structure of continuous-type sextuples in Subsection \ref{class cont sextuples} below. To further analyze continuous-type sextuples, we recall the notion of a frame in Subsection \ref{frm} and use these in Subsection \ref{framed sextuples subsection}  to build continuous-type sextuples from underlying linear endomorphisms. In Subsection \ref{identifying framed sextuples} we give a detailed description of morphisms of such sextuples, and in Subsection \ref{duals of continuous sextuples} we identify which indecomposable continuous-type sextuples are dual to which. This sets the stage for the analysis of self-dual continuous indecomposable sextuples and their compatible forms in Sections \ref{continuous real and complex} and \ref{continuous perfect}.

Naturally, when $\fieldk$ is not algebraically closed, the classification of indecomposable sextuples and associated isotropic triples becomes more complicated. First of all, the indecomposable endomorphisms $\gamma$ underlying the sextuples have a richer structure. Furthermore, self-dual continuous-type indecomposable sextuples may admit more compatible forms when $\fieldk$ is not algebraically closed. 
The general question of uniqueness of compatible forms is treated in Subsection \ref{uniq:gen}.

As always, if we are looking for compatible symplectic structures, we can restrict our attention to the cases where the integer $k$ is even.

\subsection{Classification of continuous sextuples}\label{class cont sextuples}

%
With the sole exception being the type $\Delta_2(k; 1)$, all indecomposable continuous-type sextuples (up to isomorphism) may be obtained via a functor $S$ from the ``continuous-type'' indecomposable representations of the extended Dynkin quiver $\tilde{A}_5$ (with a certain orientation, as reflected in the diagram (\ref{diag A5}))\footnote{See \cite{do-fr:representations}, pages 44 and 46}.
Indecomposable representations of this $\tilde{A}_5$ quiver have the following normal form
\begin{equation}\label{diag A5}
\begin{tikzcd}
 	& X_2  & \\
Y_1 \arrow[dd, "\alpha_1"']  \arrow[ru, "\beta_1"]   &  	    & Y_2 \arrow[lu, "\alpha_2"']  \arrow[dd, "\beta_2"]  \\
	&   		&     \\
X_1  &  	   & X_3 \\	
	  &  Y_3 \arrow[ru, "\alpha_3"']  \arrow[lu, "\beta_3"]   & 
\end{tikzcd}
\end{equation}
where $X_i = Y_i = \fieldk^k$ for some $k \in \mathbb{Z}_{>0}$, and $i \in \{1,2,3\}$, and the linear maps $\alpha_i$ and $\beta_i$ are all the identity map on $\fieldk^k$, with one exception. The map which is an exception -- we call it $\gamma$ -- is a linear endomorphism of $\fieldk^k$, and the representation of $\tilde A_5$ is indecomposable if and only $\gamma$ is an indecomposable linear endomorphism. 
The following are the continuous-type indecomposable representations of $\tilde A_5$, up to isomorphism (in each case, $k$ runs over $\mathbb{Z}_{>0}$): 
\begin{itemize}
\item $\Xi(k; \gamma)$:  $\gamma$ runs over isomorphism classes of automorphisms of $\fieldk^k$ and we assume $\gamma = \beta_1$. Any other choice of ``position'' of $\gamma$, e.g. $ \gamma = \alpha_1$, leads to a representation which is isomorphic to $\Xi(k; \gamma')$ for some  $ \gamma' = \beta_1$. 
\item $\Xi_{i}(k; 0)$ and $\Xi_{i}(k; \infty)$:  $\gamma$ is the unique indecomposable nilpotent endomorphism on $\fieldk^k$ in Jordan normal form (for nilpotent $ \gamma$ this normal form always exists), and  $i \in \{1,2,3\}$. The convention is that $\Xi_{i}(k; 0)$ denotes the case when $ \gamma = \beta_i$ and $\Xi_{i}(k; \infty)$ denotes the case when $\gamma = \alpha_i$. 
\end{itemize} 
Although DF \cite{do-fr:representations} work over algebraically closed fields, their classification of the indecomposable representations of $\tilde A_5$ does not depend on this, and admits the straightforward generalization above, where single Jordan blocks are replaced with the condition of indecomposability. Our notation is a slight modification of their notation.   

The image under $S$ of a continuous-type $\tilde A_5$ representation is the sextuple
\begin{equation}\label{sextuple normal form}
V=X_1 \oplus  X_2 \oplus X_3, \ \ C_i=X_i \oplus X_{i+1}, \ \ I_i=\im (\alpha_i \times \beta_i),  
\end{equation}
where indices are understood modulo 3. 
It will be convenient for us to cast these normal forms in slightly different notation.  We set $X = \fieldk^k$, 
$$
X_1 = X \times 0 \times 0, \ X_2 = 0 \times X \times 0, \ X_3 = 0 \times 0 \times X,
$$ 
and $V = X \times X \times X$. We call an endomorphism $\textbf{exceptional}$ if it has $0$ or $1$ as eigenvalue. Otherwise, it is \textbf{non-exceptional}. Note that any direct summand of a non-exceptional endomorphism is again non-exceptional. 

For non-exceptional $\gamma$, the sextuples which are isomorphic to $S\Xi(k; - \gamma)$\footnote{The change of sign in front of $\gamma$ here follows the conventions of DF \cite{do-fr:representations}.} will be called of type $\Delta(k; \gamma)$. Normal forms for these sextuples are 
\begin{equation}\label{normal form sextuple}
\begin{array}{ll}
I_1 =  \{ (x, -\gamma x, 0) \mid x \in X \} & \quad C_1 = X \times X \times 0 \\
I_2 = \{ (0,  x, x) \mid x \in X \} & \quad C_2 = 0 \times X \times X \\
I_3 =  \{ (x, 0, x) \mid x \in X \}  & \quad C_3 =  X \times 0 \times X.
\end{array}
\end{equation}
For indecomposable $\gamma$ having eigenvalue $1$, the sextuple $S\Xi(k; - \gamma)$ is called of type $\Delta_1(k; 1)$; its normal form is the same as above.  The only continuous-type indecomposable sextuple not isomorphic to a sextuple in the image of the functor $S$ is the type $\Delta_2(k; 1)$. It is obtained from $\Delta_1(k; 1)$ via certain functor $\theta^+$ (see \cite{do-fr:representations}, p. 38 and 46); a normal form for $\Delta_2(k; 1)$ is:
\begin{equation}\label{normal form second exceptional type}
\begin{array}{ll}
I_1  = \{ (0,  0, x) \mid x \in X \} & \quad C_1 = \{ (y,  \gamma y, x) \mid x, y \in X \} \\
I_2  = \{ (x,  0, 0) \mid x \in X \} & \quad C_2 = \{ (x, y, -y) \mid x, y  \in X \} \\
I_3  = \{ (0, x, 0) \mid x \in X \}  & \quad C_3 = \{ (-y, x, y) \mid x, y \in X \},
\end{array}
\end{equation}
where $\gamma$ is indecomposable and with eigenvalue $1$. 
Finally, for the cases when $\gamma$ is nilpotent we set $\Delta_{i}(k; 0):= S \Xi_{i}(k; 0)$ and $\Delta_{i}(k; \infty):= S \Xi_{i}(k; \infty)$, for $i =1,2,3$. For normal forms we take the same spaces $C_1, C_2, C_3$ as in (\ref{normal form sextuple}), and
\begin{itemize}
\item for $\Delta_{1}(k; 0)$ and $\Delta_{1}(k; \infty)$, respectively:
\begin{equation*}
\begin{array}{lcl}
I_1 = \{ (x, \gamma x, 0) \mid x \in X \} &  & \quad I_1 = \{ (\gamma x,  x, 0) \mid x \in X \} \\
I_2 = \{ (0,  x, x) \mid x \in X \} & \quad  \text{ and } \quad & \quad I_2 = \{ (0,  x, x) \mid x \in X \} \\
I_3 = \{ (x, 0, x) \mid x \in X \}  & &  \quad I_3 = \{ (x, 0, x) \mid x \in X \}
\end{array}
\end{equation*}
\item for $\Delta_{2}(k; 0)$ and $\Delta_{2}(k; \infty)$: 
\begin{equation*}
\begin{array}{lcl}
I_1 = \{ (x,  x, 0) \mid x \in X \} &  & \quad I_1 = \{ (x,  x, 0) \mid x \in X \} \\
I_2 = \{ (0,  x, \gamma x) \mid x \in X \} & \quad  \text{ and } \quad & \quad I_2 = \{ (0,  \gamma x, x) \mid x \in X \} \\
I_3 = \{ (x, 0, x) \mid x \in X \}  & &  \quad I_3 = \{ (x, 0, x) \mid x \in X \}
\end{array}
\end{equation*}
\item for $\Delta_{3}(k; 0)$ and $\Delta_{3}(k; \infty)$:
\begin{equation*}
\begin{array}{lcl}
I_1 = \{ (x,  x, 0) \mid x \in X \} &  & \quad I_1 = \{ (x,  x, 0) \mid x \in X \} \\
I_2 = \{ (0,  x, x) \mid x \in X \} & \quad  \text{ and } \quad & \quad I_2 = \{ (0,  x, x) \mid x \in X \} \\
I_3 = \{ (\gamma x, 0, x) \mid x \in X \}  & &  \quad I_3 = \{ (x, 0, \gamma x) \mid x \in X \}.
\end{array}
\end{equation*}
\end{itemize}

The following will be useful for identifying isomorphism types.

\begin{lemma}\label{continuous dim vectors}
Let $(V; I_i, C_i)$ be an indecomposable continuous-type sextuple with $\dim V = 3k$. Consider the following eight subspaces: $I_1 + I_2 + I_2$ and $C_1 \cap C_2 \cap C_3$, $I_1 \cap C_3$ and $I_3 \cap C_1$, $I_2 \cap C_1$ and $I_1 \cap C_2$, $I_3 \cap C_2$ and $I_2 \cap C_3$. Let $\epsilon = (\epsilon_1, ..., \epsilon_8)$ be the corresponding $8$-vector of the dimensions of these spaces. 
The different possible types of indecomposable continuous-type sextuple have the following associated $8$-vectors $\epsilon$. 
\begin{enumerate}
\item $\Delta(k; \gamma)$, then $\epsilon = (3k,0,0,0,0,0,0,0)$
\item $\Delta_1(k; 1)$, then $\epsilon = (3k-1,0,0,0,0,0,0,0)$
\item $\Delta_2(k; 1)$, then  $\epsilon = (3k,1,0,0,0,0,0,0)$
\item $\Delta_{1}(k; 0)$, then $\epsilon = (3k,0,1,0,0,0,0,0)$
\item $\Delta_{3}(k; \infty)$, then $\epsilon = (3k,0,0,1,0,0,0,0)$
\item $\Delta_{2}(k; 0)$, then $\epsilon = (3k,0,0,0,1,0,0,0)$
\item $\Delta_{1}(k; \infty)$, then $\epsilon = (3k,0,0,0,0,1,0,0)$
\item $\Delta_{3}(k; 0)$, then $\epsilon = (3k,0,0,0,0,0,1,0)$
\item $\Delta_{2}(k; \infty)$, then $\epsilon = (3k,0,0,0,0,0,0,1)$
\end{enumerate}
\end{lemma}

\pf
Consider first the sextuples for which $\gamma$ is an isomorphism. It is straightforward to check, \text{e.g.} using the normal forms above, that for such sextuples $\dim I_j \cap C_l =0$ for all $j \neq l$. Thus $\epsilon_2$ through $\epsilon_8$ are zero for the types $\Delta(k; \gamma)$, $\Delta_1(k; 1)$ and $\Delta_2(k; 1)$. 

Furthermore, if a sextuple is of type $\Delta(k; \gamma)$ or $\Delta_1(k; 1)$, then from the normal form (\ref{normal form sextuple}) we see that $C_1 \cap C_2 \cap C_3 = 0$ and that
$$
(x, y, z) \in I_1 \cap (I_2 + I_3) \  \Leftrightarrow \ (x,y,z) = (x, -\gamma x,0) \text{ with } x - \gamma x =0,
$$
so $I_1 \cap (I_2 + I_3) \neq 0$ if and only if $\gamma$ has $1$ as eigenvalue. In the case $\Delta_1(k; 1)$ when $\gamma$ does have $1$ as eigenvalue, the corresponding eigenspace has dimension $1$ (because $\gamma$ is indecomposable) and so $\dim I_1 \cap (I_2 + I_3) = 1$. Thus in this case
$$
\dim (I_1 + I_2 + I_3) = \dim I_1 + \dim (I_2 + I_3) - 1 = 3k - 1 = \dim V - 1. 
$$
So, we have found that $(\epsilon_1, \epsilon_2) = (3k, 0)$ for $\Delta(k; \gamma)$ and $(\epsilon_1, \epsilon_2) = (3k-1, 0)$ for $\Delta_1(k; 1)$.

For sextuples of type $\Delta_2(k; 1)$, it follows from the normal form (\ref{normal form second exceptional type}) that $I_1 \cap (I_2 + I_3) = 0$ and that
$$
C_1 \cap C_2 \cap C_3 = \{ (x, \gamma x, - \gamma x) \mid x = \gamma x \}. 
$$
Since for the type $\Delta_2(k; 1)$ the map $\gamma$ has a $1$-dimensional eigenspace for the eigenvalue $1$, we find that $\dim C_1 \cap C_2 \cap C_3 = 1$. So, in this case $(\epsilon_1, \epsilon_2) = (3k, 1)$.

Now consider the type $\Delta_{1}(k; 0)$. The same arguments as for the case $\Delta(k; \gamma)$ show here that $(\epsilon_1, \epsilon_2) = (3k, 0)$. Note that 
$$
I_1 \cap C_3 = \{ (x, \gamma x, 0) \in I_1 \mid \gamma x = 0 \} = \ker \gamma,
$$ 
which is $1$-dimensional since $\gamma$ is an indecomposable nilpotent map. From the normal form for $\Delta_{1}(k; 0)$ is easily check that the other intersections $I_j \cap C_l$, $j \neq l$, are zero. Thus $\epsilon_3 =1$, and $\epsilon_4$ through $\epsilon_8$ are zero.  

The remaining cases are very similar to the case $\Delta_{1}(k; 0)$ and may be treated analogously.
\qed \\

\begin{cor}\label{char via cont dim vectors}
Suppose we are given an indecomposable continuous-type sextuple with ambient dimension $3k$. 
The sextuple is of type
\begin{enumerate}
\item $\Delta(k; \gamma)$ if and only if $\dim (I_1 + I_2 + I_2) = \dim V$ and $\dim  (C_1 \cap C_2 \cap C_3) = 0$.
\item $\Delta_1(k; 1)$ if and only if $ \dim ( I_1 + I_2 + I_2 ) = \dim V - 1$ 
\item $\Delta_2(k; 1)$ if and only if $ \dim (C_1 \cap C_2 \cap C_3) = 1$
\item $\Delta_{i}(k; 0)$ if and only if $\dim ( I_i \cap C_{i-1}) = 1$.
\item $\Delta_{i}(k; \infty)$ if and only if $\dim ( I_i \cap C_{i+1}) = 1$.
\end{enumerate}
\end{cor}

\subsection{Frames}\label{frm}
Following von Neumann \cite{vNeumann},
we introduce an abstract kind of coordinate system. Given a vector space $V$, a \textbf{frame} for $V$ is a collection of five subspaces $A_1, A_2, A_3, A_{12}, A_{23}$ satisfying the following relations: 
\begin{align}
&V=A_1 \oplus A_2 \oplus A_3 \label{three pieces} \\
&A_1 + A_2 =A_1 \oplus A_{12}=A_2 \oplus A_{12} \label{A_12} \\
&A_2+A_3 =A_2 \oplus A_{23}=A_3 \oplus A_{23} \label{A_23} 
\end{align}
As a shorthand notation, we refer to a frame as $\bar A$. The notions of morphism and isomorphism of frames are the obvious ones (\text{i.e.} view a frame as a special kind of poset representation.)

The following shows that the definition could also be phrased in a way that is more symmetrical.

\begin{lemma}\label{frame symmetry}
Suppose we are given a frame $A_1, A_2, A_3, A_{12}, A_{23} \subseteq V$. 
We define $A_{31}$ by 
\begin{equation}\label{A_31 def} 
A_{31} =(A_3+A_1)\cap (A_{12}+A_{23}).
\end{equation}
Then 
\begin{equation}\label{A_{31}} 
A_3+A_1 =A_3 \oplus A_{31}=A_1 \oplus A_{31}.
\end{equation}
\end{lemma}

\pf
To see that $ A_3 \cap A_{31} = 0$, note that $\dim V = 3n$ for some $n \in \mathbb{N}$, and 
$$
A_3 \cap A_{31}  = A_3 \cap [(A_1+A_3)\cap (A_{12}+A_{23})] = A_3 \cap  (A_{12}+A_{23}). 
$$
Thus
\begin{align*}
\dim (A_{31} \cap A_3) &= \dim A_3 + \dim (A_{12}+A_{23}) - \dim (A_3 + A_{12} + A_{23}) = 3n - 3n = 0
\end{align*}
since, via (\ref{three pieces}), (\ref{A_12}), (\ref{A_23}), and (\ref{A_31 def}), 
$$
A_{12} \cap A_{23}= 0 \quad \text{and} \quad A_3 + A_{12} + A_{23} = A_3 + A_{12} + A_2 = A_3 + A_{1} + A_2 = V.
$$
Similar dimension arguments can be used to show that  $A_1 \cap A_{31} = 0$ and that $\dim A_{31} = n$. 
\qed \\

\begin{cor}\label{more frame symmetry}
If $A_1, A_2, A_3, A_{12}, A_{23} \subseteq V$ is a frame, and $A_{31}$ defined as above, then 
\begin{align}
\label{second sym} A_{12} &=(A_1+A_2)\cap (A_{23}+A_{31}), \\
\label{third sym} A_{23} &=(A_2+A_3)\cap (A_{31}+A_{12}).
\end{align}
\end{cor}

\pf
If $\{ A_1, A_2, A_3, A_{12}, A_{23} \}$ is a frame, then by Lemma \ref{frame symmetry} also 
$\{ A_1, A_2, A_3, A_{23}, A_{31} \} $
is a frame. Application of Lemma \ref{frame symmetry} to this latter frame gives (\ref{second sym}); an analogous argument gives (\ref{third sym}). 
\qed \\

The relations (\ref{three pieces}), (\ref{A_12}), (\ref{A_23}), and (\ref{A_{31}}) imply that $\dim A_1 = \dim A_2 = \dim A_3 = 1/3 \dim V$ and that each $A_{ij}$ can be interpreted as the negative\footnote{Following von Neumann, we use negative graphs for symmetry reasons when dealing with frames.} graph of a linear isomorphism $h_{ij}: A_i \rightarrow A_j$, i.e.
$$
A_{ij}=\{x -h_{ij}(x)\mid x \in A_i\},
$$
where $ij \in \{12, 23, 31 \}$. We set $h_{ii} := id$ and $h_{ji} := h_{ij}^{-1}$.  

\begin{lemma}\label{frame facts}
\
\begin{enumerate}
\item Given a frame in $V$, the associated maps satisfy $h_{jk} \circ h_{ij}= h_{ik}$ for any indices $i, j, k \in \{1,2,3\}$. 
\item Any frame is isomorphic to one built from a vector space $U$ in the following way: 
\[
\begin{array}{ll}
V= U \times U \times U  & \\
A_1= U \times 0 \times 0 & \quad A_{12}= \{(x,-x,0)\mid x \in U\} \\
A_2 = 0 \times U \times 0 & \quad A_{23}=\{(0,x,-x)\mid x \in U\} \\
A_3 = 0 \times 0 \times U & \quad  A_{31}=\{(x,0,-x)\mid x \in U\}
\end{array}
\]
\end{enumerate}
\end{lemma}

\pf
\begin{enumerate}
\item Since we are dealing only with invertible maps, equations of the form $h_{jk}\circ h_{ij}= h_{ik}$ are equivalent to ones obtained by applying, to both sides of an equation, the operations of inversion, or pre- or post-composition with one of the ``$h$'' maps. This allows one to reduce to the case of showing a single identity, say 
\begin{equation}\label{cocycle condition}
h_{23} \circ h_{12} = h_{13}. 
\end{equation}
The negative graph of $h_{13}$ is $A_{31}$ (this subspace is the negative graph of $h_{31}$, and hence also of the inverse $h_{13}$). So it is sufficient to show that the negative graph of $h_{23} \circ h_{12}$ is $A_{31}$. But
\begin{align*}
\text{graph}(- h_{23} \circ h_{12}) &= \{ x - z  \mid x \in A_1, z \in A_3, z = (h_{23} \circ h_{12})(x) \} \\
			&= \{ x - h_{12}(x) + h_{12}(x) - h_{23} (h_{12}(x))  \mid x \in A_1, z \in A_3 \} \\ 
			&\subseteq (A_1+A_3)\cap (A_{12}+A_{23}) = A_{31}, 
\end{align*}
and the last inclusion is, for dimension reasons, actually an equality.
\item Suppose we are given a frame $\bar A$ in some vector space $W$. Set $U := A_2$ and let $v_1,...,v_n$ be a basis of $A_2$. A basis of $A_1$ is defined via $u_i := h_{21}(v_i)$, $i = 1,...,n$, and a basis of $A_3$ is defined by $w_i := h_{23}(v_i)$. Since (\ref{cocycle condition}) is equivalent to 
\begin{equation}\label{trivial monodromy}
h_{31} \circ h_{23} \circ h_{12}  = \text{id}_{A_1},
\end{equation}
it follows that $u_i = h_{31}(w_i)$ for each $i$. Now the linear isomorphism which sends the basis $u_1,...,u_n, v_1,..,v_n,w_1,..,w_n$ to the basis of $V := U \times U \times U$ built canonically from $v_1,...,v_n$ has, as its image, a frame of the desired form. 
\end{enumerate}
\qed \\

\noindent
\begin{rmk}
One might think of (\ref{trivial monodromy}) as a kind of cocycle condition which says that that the endomorphism obtained from the ``loop'' $A_1 \overset{h_{12}}{\longrightarrow} A_2 \overset{h_{23}}{\longrightarrow} A_3 \overset{h_{31}}{\longrightarrow} A_1$ is trivial. 
\end{rmk}

%
%
%

Given a frame $\bar A$ on $V$, with $\dim V = 3k$, we define a \textbf{frame basis} for $\bar A$ to be an ordered basis 
$$\{u_1,...,u_k, v_1,...,v_k, w_1,...,w_k \}$$ of $V$ such that 
\begin{itemize}
\item $\{u_1,...,u_k \}$ is a basis of $A_1$, $\{v_1,...,v_k \}$ is a basis of $A_2$, $\{w_1,...,w_k \}$ is a basis of $A_3$, 
\item  $h_{12}(u_i) = v_i$, \ $h_{23}(v_i) = w_i$, \  $h_{31}(w_i) = u_i$ for all $i =1, ..., k$.
\end{itemize}
A frame basis always exists (c.f. the proof of Lemma \ref{frame facts}). 

We define an \textbf{augmented frame} in $V$ to be a frame $\bar A$ in $V$ together with a subspace $C \subseteq V$ such that $A_1 + A_2 = A_1 \oplus C$. This latter condition says that $C$ is the negative graph of a linear map $h :A_2 \rightarrow A_1$ (which is uniquely determined by $C$). In particular, an augmented frame determines uniquely an endomorphism 
$$\eta := h_{12} \circ h : A_2 \longrightarrow A_2$$ 
which we call the \textbf{underlying endomorphism} of the augmented frame. If $\eta$ is the underlying endomorphism of an augmented frame we write this as $(\bar A, \eta)$. As was the case for frames, augmented frames can be viewed as special kinds of poset representations, with the inherited notions of morphism and isomorphism. 

\begin{lemma}\label{endos frames} Let $(\bar A, \eta)$ in $V$ and $(\bar A', \eta')$ in $V'$ be augmented frames.
There is a bijective correspondence between morphisms $(A_2, \eta) \rightarrow (A_2', \eta')$ and morphisms $(\bar A, \eta) \rightarrow (\bar A', \eta')$. 

When $(A_2, \eta) = (A_2', \eta')$, this gives an isomorphism between the endomorphism algebras of $(A_2, \eta)$ and  $(\bar A, \eta)$. 
\end{lemma}

\pf
Suppose $f: V \rightarrow V'$ is a morphism of augmented frames. In particular then $f(A_1) \subseteq A_1'$, $f(A_2) \subseteq A_2'$, $f(A_{12}) \subseteq A_{12}'$, and $f(C) \subseteq C'$, which implies that the diagrams
\[
\begin{tikzcd}
 A_2 \arrow[r,"f_{\vert{A_2}}"] \arrow[d,"h"']& A_2' \arrow[d,"h'"]& &A_1 \arrow[r,"f_{\vert{A_1}}"]\arrow[d,"h_{12}"']& A_1' \arrow[d,"h_{12}'"] \\
 A_1\arrow[r,"f_{\vert{A_1}}"]& A_1'& &A_2\arrow[r,"f_{\vert{A_2}}"]& A_2'
\end{tikzcd}
\]
commute. Stacking the second diagram under the first, we obtain that 
\[
\begin{tikzcd}
 A_2 \arrow[r,"f_{\vert{A_2}}"] \arrow[d,"\eta"']& A_2' \arrow[d,"\eta'"] \\
 A_2\arrow[r,"f_{\vert{A_2}}"]& A_2'
\end{tikzcd}
\]
commutes, as desired. 

Conversely, if we have a linear map $f: A_2 \rightarrow A_2'$ such that $f \circ \eta = \eta' \circ f$, then it extends to a morphism $\hat f$ of augmented frames by setting $\hat f = h'_{21} \circ f \circ h_{12}$ on $A_1$ and $\hat f = h'_{23} \circ f \circ h_{32}$ on $A_3$.

It is easy to see that the thus defined operations ``restriction from $V$ to $A_2$'' and ``extension from $A_2$ to $V$'' are mutually inverse.

Now assume $(A_2, \eta) = (A_2', \eta')$. To show that the mutually inverse ``extension'' and ``restriction'' maps define algebra isomorphisms between the respective endomorphism algebras of $(A_2, \eta)$ and $(\bar A, \eta)$, it is sufficient to check that one of these maps is a morphism of algebras. This is easiest to check for the restriction map: the operation of restriction is clearly compatible with composition, addition, scalar multiplication, and the units in the respective endomorphism algebras.  
\qed \\

\subsection{Framed sextuples}\label{framed sextuples subsection}

Given an augmented frame $(\bar A,C)=(\bar A,\eta)$ in $V$ we define an associated sextuple $\mathcal{S}_\eta$ in $V$ by  
\begin{equation}\label{sextuple from aug frame}
\begin{array}{ll}
I_1 = A_1 &  \quad C_1 = A_1+A_2	\\
I_2 = A_3  &  \quad C_2 = A_2+A_3	\\
I_3 = (C + A_3) \cap (A_1+ A_{23})  &  \quad C_3 = A_{12} + I_3	\\
\end{array}
\end{equation}
A sextuple which is isomorphic to one of this type will be called a \textbf{framed sextuple}.

\begin{lemma}\label{augm frames} 
Let $(\bar A, \eta)$ and $(\bar A', \eta')$ be augmented frames in $V$ and $V'$, respectively .
\begin{enumerate}
\item\label{augm frames p1} Given a sextuple $\mathcal{S}_\eta$, the underlying augmented frame can be recovered via
\begin{equation}\label{underlying frame}
\begin{array}{ll}
A_1 = I_1 & \quad A_{12} = C_1 \cap C_3 \\
A_2 = C_1 \cap C_2 & \quad A_{23} = (I_3+I_1)\cap C_2  \\
A_3 = I_2 & \quad A_{31} = (A_1+A_3)\cap (A_{12}+A_{23}) \\
C = (I_2 + I_3) \cap C_1 & 
\end{array}
\end{equation}
\item\label{augm frames p2}  If $\mathcal{S}$ is a sextuple such that the expressions (\ref{underlying frame}) define an augmented frame, then $\mathcal{S}$ is a framed sextuple, \text{i.e.} of the form (\ref{sextuple from aug frame}). 
\item\label{augm frames p3} A linear map $f: V \rightarrow V'$ is a morphism $(\bar A, \eta) \longrightarrow (\bar A', \eta')$ if and only if it is a morphism $\mathcal{S}_\eta \longrightarrow \mathcal{S}_{\eta'}$. 
\item \label{augm frames part 3} Any sextuple $\mathcal{S}_\eta$ built from an augmented frame $(\bar A, \eta)$  is isomorphic to one with the following form, where $U = A_2$:
\begin{equation}\label{S_eta normal form}
\begin{array}{ll}
V = U \times U \times U &   \\
I_1 = U \times 0 \times 0 & \quad C_1 = U \times U \times 0 \\
I_2 = 0 \times 0 \times U & \quad C_2 = 0 \times U \times U \\
I_3 = \{(-\eta x,x,-x)\mid x \in U\} & \quad C_3 = \{(x,-x,0)\mid x\in U\} + I_3. 
\end{array}
\end{equation}
\end{enumerate} 
\end{lemma}

\pf
\begin{enumerate}
\item This can be checked using elementary linear algebra resp. modular lattice calculations. For example, 
$C_1 \cap C_2 = (A_1 + A_2) \cap (A_2 + A_3) = A_2$, since by assumption $V = A_1 \oplus A_2 \oplus A_3$. 
To see that $C = (I_2 + I_3) \cap C_1$, we plug in the definitions of $I_2$, $I_3$ and $C_1$ and calculate
\begin{align*}
(I_2 + I_3) \cap C_1 &=  (A_3 + [(C + A_3) \cap (A_1+ A_{23})]) \cap (A_1 + A_2) \\
			&\subseteq [(C+ A_3) \cap (A_1 + A_{23} + A_3)] \cap (A_1 + A_2) \\
			& = [(C + A_3) \cap V] \cap (A_1 + A_2) = C.
\end{align*}
The obtained inclusion is actually an equality, since
$$
\dim (I_3 + I_2) \cap C_1 = \dim (I_1 \oplus I_3) + \dim C_1 - \dim (I_2 + I_3 + C_1) = 1/3 \dim V = \dim C.
$$
The equations for $A_{12}$ and $A_{23}$ can be checked in a similar manner, and the equation for $A_{31}$ holds by the definition of a frame. 

\item 
Assuming that (\ref{underlying frame}) defines an augmented frame, we need to show that the relations (\ref{sextuple from aug frame}) hold. 

The expressions for $I_1$ and $I_2$ are trivially satisfied, so it remains to show the expressions for $I_3$, $C_1$, $C_2$, and $C_3$. 

We first make some intermediate observations:
\begin{equation}\label{first calc}
A_1 + C = A_1 + C_1 \cap (A_3 + I_3)  = C_1 \cap (A_1 + A_3 + I_3),
\end{equation}
using the modular law for the last equality; from (\ref{first calc}) and the assumption that we have an augmented frame,
\begin{equation}\label{second calc}
\begin{array}{c}
V = A_1 + A_2 + A_3 =  A_1 + C + A_3 = C_1 \cap (A_1 + A_3 + I_3) + A_3 \\
				= (C_1 + A_3) \cap (A_1 + A_3 + I_3), 
\end{array}
\end{equation}
again via the modular law for the last equality. But (\ref{second calc}) implies that 
\begin{equation}\label{third calc}
V = C_1 + A_3  = A_1 + A_3 + I_3,
\end{equation}
and using (\ref{third calc}) we find that 
\begin{equation}\label{fourth calc}
A_2 + A_3 = C_1 \cap C_2 + A_3 = C_2 \cap (C_1 + A_3) = C_2, 
\end{equation}
where the application of modular law is justified since $A_3 = I_2 \subseteq C_2$ by assumption. This is the desire expression for $C_2$. 

The expression for $C_1$ now follows from (\ref{third calc}) and 
\begin{equation}\label{fifth calc}
\begin{array}{c}
A_1 + A_2 = A_1 + C = A_1 + C_1 \cap (I_2 + I_3) = C_1 \cap ( A_1 + I_2 + I_3 ) = C_1,
\end{array}
\end{equation}
using modularity via the fact that $A_1 = I_1 \subseteq C_1$. 

To obtain the expression for $I_3$, we first note that
\begin{equation}\label{sixth calc}
\begin{array}{c}
C + A_3 = (I_2 + I_3) \cap C_1 + A_3 = (I_2 + I_3) \cap (C_1 + A_3) = I_2 + I_3
\end{array}
\end{equation}
since $A_3 = I_2 \subseteq I_2 + I_3$ and, via (\ref{third calc}), $C_1 + A_2 = V$. Now, 
\begin{align*}
(C + A_3) \cap (A_1 + A_{23}) & = (C + A_3) \cap (A_1 + (I_3 + A_1) \cap C_2)  \\
&\overset{(\ref{fourth calc})}{=}  (C + A_3) \cap (A_1 + (I_3 + A_1) \cap (A_2 + A_3)) \\ 
&\overset{\text{mod}}{=} (C + A_3) \cap (I_3 + A_1) \cap (A_1 + A_2 + A_3) \\
&\overset{(\ref{third calc})}{=} (C + A_3) \cap (A_1 + I_3) \overset{(\ref{sixth calc})}{=} I_3 + (C+ A_3) \cap A_1 \\
&= I_3 + (C \cap A_1) = I_3,
\end{align*}
using for the second-to-last equality that $C \subseteq A_1 \oplus A_2$ and $V = A_1 \oplus A_2 \oplus A_3$, and for the last equality that $C \cap A_1 = 0$ by definition of an augmented frame. 

Finally, the desired expression for $C_3$ follows from 
\begin{align*}
A_{12} + I_3 &= (C_1 \cap C_3) + I_3 \overset{\text{mod}}{=} C_3 \cap (C_1 + I_3) \\
	&= C_3 \cap (A_1 + A_2 + (C + A_3) \cap (A_1 + A_{23})) \\
	&\overset{\text{mod}}{=} C_3 \cap (A_2 + (A_1 + C + A_3) \cap (A_1 + A_{23})) \\
	&= C_3 \cap (A_2 + A_1 + A_{23}) = C_3,
\end{align*} 
using in the second line the already-obtained expressions for $C_1$ and $I_3$, and in the last line the fact that $A_1 + C + A_3 = A_1 + A_2 + A_3 = V$. 

\item Suppose first that $f$ is a morphism of augmented frames. Since the subspaces $I_i$,$C_i$ and $I_i', C_i'$ of the respective associated sextuples can be expressed entirely via lattice terms built only of subspaces in the respective augmented frames, it follows that $f(I_i) \subseteq I_i'$ and $f(C_i) \subseteq C_i'$ for each $i$. 

Similarly, by part \ref{augm frames p1} above the underlying augmented frames can be expressed via lattice terms built from subspaces in the respective sextuples, so $f$ is a morphism of augmented frames if it is a morphism of the associated sextuples. 
\item From Lemma \ref{frame facts}, it is clear that any augmented frame $(\bar A, \eta)$ is isomorphic to one of the following form:
\begin{equation}\label{underlying aug frame}
\begin{array}{ll}
V= U \times U \times U  & \\
A_1= U \times 0 \times 0 & \quad A_{12}= \{(x,-x,0)\mid x \in U\} \\
A_2 = 0 \times U \times 0 & \quad A_{23}=\{(0,x,-x)\mid x \in U\} \\
A_3 = 0 \times 0 \times U & \quad  A_{31}=\{(x,0,-x)\mid x \in U\} \\
C = \{ (-\eta x, x, 0) \mid x \in U \} & 
\end{array}
\end{equation}
The sextuple associated to this augmented frame is precisely (\ref{S_eta normal form}), and by part \ref{augm frames p3} of this lemma, it is isomorphic to $S_\eta$. 
\end{enumerate}  \qed \\

\begin{rmk}
Together, parts \ref{augm frames p1} and \ref{augm frames p2}  give a complete characterization of framed sextuples in lattice theoretic terms: a sextuple is a framed sextuple if and only if the expressions (\ref{underlying frame}) define an augmented frame, and augemented frames are themselves defined in lattice-theoretic terms. 
\end{rmk}

%

\begin{prop}\label{endo to sextuple} Let $(U, \eta)$, $(U', \eta')$ be spaces with endomorphisms.  
\begin{enumerate}
\item There is a bijective correspondence between morphisms $(U, \eta) \rightarrow (U', \eta')$ and morphisms $S_\eta \rightarrow S_\eta'$.  When $(U, \eta) = (U', \eta')$, this correspondence is an isomorphism of the respective endomorphism algebras. In particular, $(U, \eta)$ is indecomposable if and only if $S_\eta$ is. 
\item Let $u_1,...,u_n, v_1,...,v_n,w_1,..,w_n$ and $u_1',...,u_n', v_1',...,v_n',w_1',..,w_n'$ be frame bases of the underlying frames of sextuples $S_\eta$ and $S_{\eta'}$. Let $\hat f : S_\eta \rightarrow S_{\eta'}$ be a morphism and  $f$ its restriction $f: A_2 = \langle  v_1,...,v_n \rangle \rightarrow A_2' = \langle  v_1',...,v_n' \rangle$. If the coordinate matrix of $f$ with respect to the bases $\langle  v_1,...,v_n \rangle$ and $\langle  v_1',...,v_n' \rangle$ is $M$, then the coordinate matrix of $\hat f$ with respect to the respective frame bases is 
\begin{equation}\label{block matrix}
\left [
\begin{array}{ccc}
M & 0 & 0 \\
0 & M & 0 \\
0 & 0 & M
\end{array}
\right].
\end{equation}
\end{enumerate}
\end{prop}

\pf
\begin{enumerate}
\item \label{functorial corr} The said correspondence is the one defined in the proof of Lemma \ref{endos frames}. By Lemma \ref{endos frames} and  Lemma \ref{augm frames}, \text{2.}, it maps morphisms $(U, \eta) \rightarrow (U', \eta')$ to morphisms $S_\eta \rightarrow S_\eta'$, and when $(U, \eta) = (U', \eta')$, this correspondence is an isomorphism of algebras. The statement about indecomposability follows then from the fact that $(U,\eta)$ and $S_\eta$ are each indecomposable if and only if their respective endomorphism algebras are local, and the fact that ``local-ness'' is preserved under isomorphism. 
\item 
By definition, $\hat f = (h'_{21} \circ f \circ h_{12}) \oplus f \oplus (h'_{23} \circ f \circ h_{32}) : A_1 \oplus A_2 \oplus A_3 \rightarrow A_1' \oplus A_2' \oplus A_3'$ (c.f. Lemma \ref{augm frames}). The form (\ref{block matrix}) follows now from the fact that, with respect to frame bases, all of the maps $h_{12}$, $h_{23}$, $h_{21}'$, $h_{32}'$ have coordinate matrices which are the identity matrix.
\end{enumerate} \qed 
\begin{rmk}
The correspondence in Proposition \ref{endo to sextuple}, is functorial. 
\end{rmk}

\subsection{Identifying framed sextuples}\label{identifying framed sextuples}

In this section we identify which continuous-type indecomposable sextuples are isomorphic to a framed sextuple. 

From Lemma \ref{continuous dim vectors} and Corollary \ref{char via cont dim vectors}  we see that, for $\eta$ indecomposable, $S_\eta \simeq \Delta_3(k; 0)$ when $\eta$ is nilpotent, and $S_\eta \simeq \Delta_2(k; \infty)$ when $\eta$ has eigenvalue $1$. Indeed, if $\eta$ is nilpotent, then
$$
I_3 \cap C_2 = \{ (-\eta x, x, -x) \mid x \in U \} \cap (0 \times U \times U)  \supseteq \{ (0, x, -x) \mid x \in \text{ker} \eta \} \neq 0.
$$
And if $\eta$ has eigenvalue $1$ with associated eigenspace $U_1$, then 
$$
I_2 \cap C_3 = (0 \times 0 \times U) \cap \{ (y -\eta x, x - y, -x) \mid x, y \in U \}  \supseteq \{ (0, 0, -x) \mid x \in U_1 \} \neq 0.
$$
We will call an indecomposable endomorphism $\textbf{exceptional}$ if it is nilpotent or has $1$ as eigenvalue. Note that any non-exceptional indecomposable endomorphism is an automorphism. 

It remains now to identify the sextuples $S_\eta$, with $\eta$ non-exceptional, in terms of the classification discussed in the previous section.

\begin{prop}\label{classification endo-sextuples} 
On the level of isomorphism classes (and for fixed ambient dimension $3k$), there is a one-to-one correspondence between
the sextuples $S_\eta$, and the sextuples $\Delta(k; \gamma)$, where $\gamma$ and $\eta$ are non-exceptional indecomposable endomorphisms. If one views both $\eta$ and $\gamma$ as endomorphisms of $\fieldk^k$, the correspondence is given by 
\begin{equation}\label{bijection of endos}
\eta = \frac{\gamma}{\gamma - 1}
\end{equation} 
or, in inversely, $\gamma = \tfrac{\eta}{\eta -1}$. 
\end{prop}

\pf 
We start by considering an indecomposable continuous-type sextuple  $\Delta(k; \gamma)$, with $\gamma$ non-exceptional. 
From Section \ref{class cont sextuples}, we know that $\Delta(k; \gamma)$ is isomorphic to a sextuple which has the form
$$
V=X_1\oplus X_2 \oplus X_3, \ \ C_i=X_i+X_{i+1}, \ \ I_i=\im (\alpha_i\oplus \beta_i) = \Gamma (\beta_i \alpha_i^{-1})
$$
where $\alpha_i : Y_i \rightarrow X_i$ and $\beta : Y_i \rightarrow X_{i+1}$ are all invertible, and $\gamma = \beta_1 \circ \alpha_1^{-1}$. 
In particular, we can identify the $Y_i$ with the $I_i$ (since the maps $\alpha_i \oplus \beta_i$ are injective), and we have
\begin{equation}\label{up to perm}
 C_i= X_i\oplus Y_i=X_{i+1}\oplus Y_i \quad i=1,2,3.
\end{equation}

To show that $\Delta(k; \gamma) \simeq S_\eta$ for some endomorphism $\eta$, we ``guess'' an underlying augmented frame (using Lemma \ref{augm frames} to make our ansatz), and we show that this is an augmented frame whose associated sextuple is isomorphic to $\Delta(k; \gamma)$. In this case we know that $\eta$ must be non-exceptional, since $S_\eta \simeq \Delta_3(k; 0)$ or $S_\eta \simeq \Delta_2(k; \infty)$ when $\eta$ is exceptional.

As our ansatz for the underlying augmented frame associated to $(V; C_i, I_i)$, we set
$$
\begin{array}{ll}
A_1=Y_1 & A_{12}= (X_1+X_2)\cap (X_3+X_1) =X_1 \\
A_2=(X_1+X_2)\cap (X_2+X_3)=X_2 \ & A_{23}=(Y_3+Y_1)\cap(X_2+X_3) \\
A_3=Y_2 & A_{31} = (A_1+A_3)\cap (A_{12}+A_{23})  \\
C=(Y_2+Y_3)\cap(X_1+X_2) &  
\end{array}
$$
and check that this defines an augmented frame. In order to verify a $\oplus$-relation, it suffices to
check the $+$- or the $\cap$-relation, if the dimensions
of the subspaces involved are known and add up properly. 
Thus, to see that $V = A_1 \oplus A_2\oplus A_3$, it suffices to note that
\[
A_1 + A_2 + A_3  =  Y_1 + X_2 + Y_2  = X_1+X_2 + Y_2 = X_1 + X_2 + X_3 = V.
\]
Similarly, $A_{12}\oplus A_1=A_{12}\oplus A_2=A_1+A_2$ follows from
$$
\begin{array}{l}
A_{12}+ A_1= X_1+Y_1=X_1+X_2=A_1+A_2 \\
A_{12}+A_2= X_1+X_2=A_1+A_2.
\end{array}
$$
To see $A_{23}\oplus A_2=A_{23}\oplus A_3=A_2+A_3$, note that 
\begin{align*}
A_{23}+A_2= & \ (Y_1+Y_3)\cap (X_2+X_3)+X_2
=(Y_1+Y_3+X_2)\cap (X_2+X_3) \\
&=X_2+X_3 =X_2+Y_2= A_2+A_3 \\
A_{23}+A_3= & \  (Y_1+Y_3)\cap (X_2+X_3)+Y_2
=(Y_1+Y_3+Y_2)\cap (X_2+X_3) \\
&=X_2+X_3=
X_2+Y_2= A_2+A_3
\end{align*}
using modularity (and that $X_2, Y_2 \subseteq X_2+X_3$) to obtain the second equality in each line, respectively. At this point it is not yet clear that the lefthand sums are direct, since the dimension of $A_{23}$ is not yet determined. This can be checked directly, for example:
\[
Y_1+Y_3+X_2=X_1+Y_3+X_2=X_1+X_3+X_2=V 
\mbox{ whence } A_{23}\cap A_2=(Y_1+Y_3)\cap X_2 =0
\]
and
\[ 
A_{23}\cap A_3=(Y_1+Y_3)\cap Y_2 =0.
\]
Finally, $ A_1\oplus C=A_1+A_2$ because
$$
\begin{array}{l}
A_1+C=Y_1+(Y_2+Y_3)\cap(X_1+X_2)=(Y_1+Y_2+Y_3)\cap (X_1+X_2)
=X_1+X_2=A_1+A_2 \\
A_1\cap C=Y_1\cap  (Y_2+Y_3)=0
\end{array}
$$
using modularity in the first line (with the fact that $Y_1 \subseteq X_1+X_2$). 

This establishes that we have an augmented frame. It remains to verify that
the original sextuple $(V,I_i,C_i)$ is the one associated to this augmented
frame, \text{i.e.} we check that the equations (\ref{sextuple from aug frame}) hold:
\begin{itemize}
\item  $I_1=Y_1=A_1 $
\item $I_2=Y_2=A_3 $
\item $I_3=Y_3= (Y_2+Y_3)\cap (Y_1+Y_3)=(C + A_3)\cap
(A_1+A_{23}) $, 

using for the last equality that 
\begin{align*}
Y_2+Y_3 =& \  (Y_2+Y_3)\cap (X_1+X_2+X_3) = (Y_2+Y_3)\cap (X_1+ X_2 + Y_2) \\
	&= (Y_2+Y_3)\cap (X_1+X_2)+Y_2 = C+A_{3}
\end{align*}
and
\begin{align*}
Y_1+Y_3 =& \ (Y_1+Y_3)\cap (X_1+X_2+X_3)= (Y_1+Y_3)\cap (Y_1+X_2+X_3) \\
	&= Y_1+ (Y_1+Y_3)\cap(X_2+X_3)=(A_1+A_{23}) 
\end{align*}
\item $C_1=A_1+A_2=Y_1+X_2=X_1+X_2=C_1$
\item $C_2=A_2+A_3= X_2+Y_2=X_2+X_3=C_2$
\item $C_3=A_{12}+(C +A_3)\cap (A_1+ A_{23})=A_{12}+I_3
= X_1+Y_3=X_1+X_3=C_3 $
\end{itemize}
This establishes that every sextuple $\Delta(k; \gamma)$ is isomorphic to some $S_\eta$, with $\eta$ non-exceptional. It remains now to show that every $S_\eta$ with $\eta$ non-exceptional is isomorphic to some $\Delta(k; \gamma)$. For this it is sufficient to prove the formula (\ref{bijection of endos}) and note that this formula defines a bijection (in fact an involution) of the set of non-exceptional indecomposable endomorphisms of $\fieldk^k$.

Consider again a sextuple of type $\Delta(k; \gamma)$. Set $X = \fieldk^k$. We'll use normal forms which are isomorphic to the ones (\ref{sextuple normal form}): let $V = X \times X \times X$, and
$$
X_1 = X \times 0 \times 0, \ X_2 = 0 \times X \times 0, \ X_3 = 0 \times 0 \times X.
$$ 
Thus normal forms for these sextuple are
\begin{equation}\label{}
\begin{array}{ll}
I_1 = \Gamma (\beta_1 \alpha_1^{-1}) = \{ (x, -\gamma x, 0) \mid x \in X \} & \quad C_1 = X \times X \times 0 \\
I_2 = \Gamma (\beta_2 \alpha_2^{-1}) = \{ (0,  x, x) \mid x \in X \} & \quad C_2 = 0 \times X \times X \\
I_3 = \Gamma (\beta_3 \alpha_3^{-1}) = \{ (x, 0, x) \mid x \in X \}  & \quad C_3 =  X \times 0 \times X.
\end{array}
\end{equation}
Set $g_{i}:= \beta_i \alpha_i^{-1}$ and $g := g_1g_3g_2$. Note that $g = - \gamma$ when these are viewed as maps $X \rightarrow X$.
  
For such a sextuple, we compute the endomorphism $\eta$ underlying the associated augmented frame. By definition, $\eta = h_{12} \circ h$, so we need to compute $h_{12}$ and $h$.  From the first part of this proof we know that for this frame 
\begin{equation*}\label{}
\begin{array}{l}
A_1 = I_1  \\
A_2 = X_2  \\
A_3 = I_2  \\
A_{12} = \Gamma(- h_{12}) = X_1 \\
C = \Gamma(-h) = (I_2 + I_3) \cap (X_1 + X_2).
\end{array}
\end{equation*}
In particular one finds easily that 
\begin{equation*}
h_{12} :  A_1  \rightarrow  A_2, \ (x, g_1x, 0) \mapsto (0, g_1x,0),
\end{equation*}
$
C = \{ (-g_3g_2x, x, 0) \mid x \in X \},
$
and 
\begin{equation*}
h : A_2  \rightarrow A_1, \ (0, x, 0) \mapsto (g_3g_2 (g + 1)^{-1}x, g(g+ 1)^{-1}x,0).
\end{equation*}
It follows that 
\begin{equation*}
\eta :  A_2  \rightarrow  A_2, \ (0, x, 0) \mapsto (0, g(g+ 1)^{-1}x,0).
\end{equation*}
Thus, viewed as endomorphisms of $X$, $\eta = \tfrac{\gamma}{1 - \gamma}$. 
\qed \\

\begin{rmk}
Since $\eta$ and $\gamma$ commute, if $\gamma$ has an eigenvalue, say $\lambda_\gamma$, then $\eta$ will also have an eigenvalue, say $\lambda_\eta$, and any eigenvector for $\lambda_\gamma$ will also be an eigenvector for $\lambda_\eta$. In this case,  $\lambda_\eta = \lambda_\gamma/(\lambda_\gamma - 1)$ and  $\lambda_\gamma = \lambda_\eta/(\lambda_\eta - 1)$.
\end{rmk}

\subsection{Duals of continuous sextuples}\label{duals of continuous sextuples}

In this section, we identify the duals of indecomposable continuous-type sextuples. We recall: 
\begin{itemize}
\item The {\bf dual} of a sextuple $(V;C_i,I_i)$ 
is $(V^*;I_i^\circ,C_i^\circ)$, where, for $U \subseteq V$, the subspace
$U^\circ=\{f \in V^*\mid f(U)=0\}$ is the annihilator of $U$. 
\item A sextuple is {\bf self-dual} if it admits an isomorphism (of poset representations)
to its dual. A pair of sextuples is called \textbf{mutually dual} if each is isomorphic to the dual of the other. 
\item The operation of taking the annihilator obeys the rules 
$$(U_1 + U_2)^\circ  = U_1^\circ \cap U_2^\circ \quad \text{and} \quad (U_1 \cap U_2)^\circ  = U_1^\circ + U_2^\circ,$$
for any subspaces $U_1, U_2$. 
\end{itemize}

From the structure of the dimension vector of indecomposable continuous-type sextuples it follows that the dual of a continuous-type sextuple is again of continuous type. Identifying the duals of the exceptional continuous-type sextuples is easiest. 

\begin{lemma}
The indecomposable sextuples of type $\Delta_1(k; 1)$ and $\Delta_2(k; 1)$ are mutually dual. 
\end{lemma}

\pf
This follows from  Corollary \ref{char via cont dim vectors}: if $(V; C_i, I_i)$ is a sextuple of type $\Delta_1(k; 1)$, then $\dim (I_1 + I_2 + I_3) = \dim V -1$. Thus for the dual $(V'; C_i', I_i')$ will hold $C_1' \cap C_2' \cap C_3' = I_1^\circ \cap I_2^\circ \cap I_3^\circ = (I_1 + I_2 + I_3 )^\circ = 1$. This implies, by Corollary \ref{char via cont dim vectors} that the dual is of type $\Delta_2(k; 1)$ (the dual must be indecomposable and of continuous type).   
\qed \\

\begin{lemma}
The indecomposable sextuples $\Delta_{i}(k; 0)$ and $\Delta_{i-1}(k; \infty)$ are mutually dual. In other words, in each ambient dimension $3k$, we have the following three pairs of mutually dual indecomposable sextuples:
\begin{equation*}
\Delta_{1}(k; 0) \text{ and } \Delta_{3}(k; \infty), \quad \Delta_{2}(k; 0) \text{ and } \Delta_{1}(k; \infty), \quad \Delta_{3}(k; 0) \text{ and } \Delta_{2}(k; \infty).
\end{equation*}
\end{lemma}

\pf
Suppose $(V; C_i, I_i)$ is a sextuple of type $\Delta_{i}(k; 0)$. By Corollary \ref{char via cont dim vectors}, this sextuple will satisfy $ \dim ( I_i \cap C_{i-1}) = 1$. In particular then 
$$
\dim (I_i + C_{i-1}) = \dim I_i + \dim C_{i-1} - \dim ( I_i \cap C_{i-1}) = k + (2k) - 1 = \dim V - 1. 
$$
The dual sextuple $(V'; C_i', I_i')$ will therefore satisfy
$$
 \dim ( I_{i-1}' \cap C_{i}') = \dim ( C_{i-1}^\circ \cap I_{i}^\circ) = \dim (I_{i}+ C_{i-1})^\circ = 1. 
$$ 
This implies, via Corollary \ref{char via cont dim vectors}, that $(V'; C_i', I_i')$ is of type $\Delta_{i-1}(k; \infty)$
\qed \\

\begin{rmk}
From Section \ref{identifying framed sextuples} we know that $\Delta_{3}(k; 0) \text{ and } \Delta_{2}(k; \infty)$ are isomorphic, respectively, to the indecomposable framed sextuples $S_\eta$ where $\eta$ is either nilpotent (in the first case) or has eigenvalue $1$ (in the second case). Thus the above shows that the sole two exceptional framed sextuples are dual to one another. 
\end{rmk}

Now we consider the non-exceptional indecomposable continuous-type sextuples $\Delta(k; \gamma)$. From  Corollary \ref{char via cont dim vectors} it follows that the dual of such a sextuple is again a non-exceptional indecomposable continuous-type sextuple; let $\Delta(k; \gamma')$ be the type of the dual. Our goal now is to determine the relationship between $\gamma$ and $\gamma'$. 

From Proposition \ref{classification endo-sextuples}  we know that $\Delta(k; \gamma)$ and $\Delta(k; \gamma')$ are isomorphic, respectively, to sextuples $S_\eta$ and $S_{\eta'}$, with $\eta$ and $\eta'$ non-exceptional. 


\begin{prop}\label{duals of endo-sextuples}
Consider a non-exceptional indecomposable endomorphism $\eta \in \text{End}(U)$ with  associated sextuple $\mathcal{S}_\eta$. Let 
$ \beta := (u_1,...,u_k, v_1, ..., v_k, w_1,.., w_k)$ be a frame basis for this sextuple.
Then:
\begin{enumerate}
\item The dual of $\mathcal{S}_\eta$ is isomorphic to $\mathcal{S}_{\eta'}$, with $\eta' = (\text{Id} - \eta)^* \in \text{End}(U^*)$. Moreover, 
$$(-w_1^*,...,-w_k^*, v_1^*,..., v_k^*, -u_1^*, ..., -u_k^*) $$ 
is a frame basis for $\mathcal{S}_{\eta'}$, where $\beta^* := (u_1^*, ..., w_k^*)$ is the dual basis of $\beta$. 
\item A bijective correspondence 
$$B : \mathcal{S}_\eta \overset{\sim}{\longrightarrow} \mathcal{S}_{\eta'} \quad \quad \longleftrightarrow \quad \quad b: (U, \eta) \overset{\sim}{\longrightarrow} (U^*, \text{Id} - \eta^*)$$
is given by restricting maps $B$ to $U$. 

With respect to the bases $\beta$ and $\beta^*$, any isomorphism $B : \mathcal{S}_\eta \rightarrow \mathcal{S}_{\eta'}$ has coordinate matrix of the form
$$ H_M :=\left(\begin{array}{ccc} O&O&-M\\O&M&O\\-M&O&O \end{array}\right),$$
where $M \in   \fieldk^{k\times k}$ is the coordinate matrix of the corresponding isomorphism $b: (U, \eta) \rightarrow (U^*, \eta')$ with respect to the respective bases $(v_1,..., v_k)$ and $(v_1^*,..., v_k^*)$ of $U$ and  $U^*$. 

Note, in particular, that $B$ is (skew)-symmetric if and only if the corresponding map $b$ is. 
\end{enumerate}
\end{prop}


\pf
%
%
By definition, $\mathcal{S}_\eta$ has ambient space $V = U \times U \times U$, and the associated augmented frame $(\bar A, \eta)$ is (\ref{underlying aug frame}); in particular $A_1 = U \times 0 \times 0$, $A_2 = 0 \times U \times 0$ and $A_3 = 0 \times 0 \times U$. We view the dual sextuple $\mathcal{S}_\eta^* = S_{\eta'} = (V'; C'_i, I'_i)$ as having ambient space $V' = U^* \times U^* \times U^*$, paired with $V = U \times U \times U$ via $(\xi_1,\xi_1, \xi_3): (u_1, u_2, u_3) \mapsto \xi_1(u_1) +  \xi_2(u_2) + \xi_3(u_3)$.  The sextuple $\mathcal{S}_{\eta'}$, expressed in terms of the underlying augmented frame of $\mathcal{S}_\eta$, is 
\begin{equation}
\begin{array}{ll}
I_1' = C_1^\circ = A_1^\circ \cap A_2^\circ = 0 \times 0 \times U^*  & \ C_1' = I_1^\circ = A_1^\circ = 0 \times U^* \times U^* \\
I_2' = C_2^\circ = A_2^\circ \cap A_3^\circ = U^* \times 0 \times 0  & \ C_2' = I_2^\circ = A_3^\circ = U^* \times U^* \times 0  \\
I_3' = C_3^\circ = A_{12}^\circ \cap [(C^\circ \cap A_3^\circ) + (A_1^\circ \cap A_{23}^\circ)] & \ C_3' = I_3^\circ = (C^\circ \cap A_3^\circ) + (A_1^\circ \cap A_{23}^\circ)
\end{array}
\end{equation}
We compute the underlying augmented frame $(\bar A', \eta')$ of this sextuple, using (\ref{underlying aug frame}) as an aid for the calculations: 
\begin{equation}
\begin{array}{l}
A_1' = I_1' = 0 \times 0 \times U^*  \\
A_2' = C_1' \cap C_2' = 0 \times U^* \times 0  \\
A_3' =  I_2' = U^* \times 0 \times 0 \\
A_{12}'  = C_1' \cap C_3' = (I_1 + I_3)^\circ = \{(y, x, -x) \mid x,y \in U \}^\circ =  \{(0, \xi, \xi) \mid \xi \in U^* \} \\
A_{23}' = (I_3' + I_1') \cap C_2' = [(C_3 \cap C_1) + I_2]^\circ =  \{(\xi, \xi, 0) \mid \xi \in U^* \}  \\
A_{31}' = (A_1'+A_3')\cap (A_{12}'+A_{23}') \\
C' = (I_2' + I_2') \cap C_1' = [(C_2 \cap C_3) + I_1]^\circ = \{(0, \xi, (\text{id} - \eta)^*\xi) \mid \xi \in U^* \} 
\end{array}
\end{equation}
For the computation of $C'$, for example:  
$
C_2 \cap C_3 =\{(y-\eta x, x-y,-x)\mid x,y \in U, y- \eta x=0\},
$
and
$I_1+C_2 \cap C_3 = \{(y, x- \eta x,-x)\mid x,y \in U\}$,
so $C'$ is precisely those $(\xi_1, \xi_2, \xi_3) \in U^* \times U^* \times U^*$ such that $\xi_1 =0$ and $\xi_3 = \xi_2 \circ (\text{id} - \eta)$. 

We can now read off the maps $h'_{i,i+1}$ and $h'$ associated to this augmented frame:
\begin{itemize}
\item $A_{12}' = \Gamma(-h_{12}')$ implies that $h_{12}': A_1' \rightarrow A_2', (0,0, \xi) \mapsto (0, -\xi, 0)$.
\item $A_{23}' = \Gamma(-h_{23}')$ implies that $h_{23}': A_2' \rightarrow A_3', (0, \xi, 0) \mapsto (-\xi,0, 0)$.
\item $h_{31}' = (h_{23}' \circ h_{12}')^{-1} : A_3' \rightarrow A_1', (\xi,0,0) \mapsto (0,0, \xi)$.  
\item $C' = \Gamma(-h')$ implies that $h': A_2' \rightarrow A_1, (0, \xi , 0) \mapsto (0,0, - (\text{id} - \eta)^*\xi )$
\end{itemize}
In particular, $\eta' = h_{12}' \circ h' = (\text{id} - \eta)^*$. We also see that, taking $(-v_1^*, ..., -v_k^*)$ as a basis of $A_2' = 0 \times U^* \times 0$, a frame basis is given by $(w_1^*, ..., w_k^*, -v_1^*, ..., -v_k^*, u_1^*, ..., u_k^*)$, since $w_i^* =  h_{21}'(-v_i^*)$ and $u_i^* = h_{23}'(-v_i^*)$, $i =1,.., k$. The remaining statement of the current proposition follows now from Proposition \ref{endo to sextuple}. In particular,  any isomorphism $\mathcal{S}_\eta \rightarrow \mathcal{S}_{\eta'}$ has, with respect to the bases $\beta$ and $\beta^*$, coordinate matrix of the form $H_M$, where $M$ is the coordinate matrix of an isomorphism $U \rightarrow U^*$ such that $M A M^{-1} = \text{id} - A^t$. 
\qed \\

\begin{cor}\label{cont type duals}
The sextuples $\Delta(k; \gamma)$ and $\Delta(k; (\gamma^{-1})^*)$ are mutually dual.
\end{cor}

\pf
Viewing $\gamma$ and $\eta$ as endomorphisms of the same space, and similarly for $\gamma'$ and $\eta'$, we have by (\ref{bijection of endos}) that $\eta = \tfrac{\gamma}{\gamma - 1}$ and $\gamma' =  \tfrac{\eta'}{\eta' -1}$. Substituting $\eta' = 1- \eta^*$ in the latter equation, using the former equation and simplifying gives $\gamma' =  (\gamma^{-1})^*$. 
\qed \\

\begin{cor}
A sextuple $\Delta(k; \gamma)$ is self-dual if and only if $(\fieldk^k, \gamma)$ and $((\fieldk^k)^*, (\gamma^{-1})^*)$ are isomorphic endomorphisms (in the sense of Section \ref{decomp linear reps}). In this case, if $\gamma$ has an eigenvalue $\lambda_\gamma$, then $\lambda_\gamma = -1$. 
\end{cor}

\pf 
The first part follows from Corollary \ref{cont type duals} and the fact that, by Propositions \ref{classification endo-sextuples}  and \ref{endo to sextuple}, $\Delta(k; \gamma_1)$ and $\Delta(k; \gamma_2)$ are isomorphic sextuples if and only if $\gamma_1$ and $\gamma_2$ are isomorphic endomorphisms. The statement about eigenvalues follows from the fact that $\gamma^*$ has the same eigenvalue as $\gamma$ (when such exists), and that $\gamma$ cannot have eigenvalue $1$ (by the assumption that $\gamma$ is non-exceptional). 
\qed \\


\section{Continuous non-split isotropic triples over $\mathbb{C}$ and $\mathbb{R}$}\label{continuous real and complex}

We turn now to the classification of non-split isotropic triples of continuous type. The special cases where the ground field $\fieldk$ is either $\mathbb{C}$ or $\mathbb{R}$ are of particular interest. We treat these cases first, which allows for a simpler analysis. Furthermore, the case $\fieldk = \mathbb{R}$ provides basic intuition for understanding the more involved classification over perfect fields, which we undertake in the next section. 

Throughout this and the next section, $N$ will always denote a nilpotent square matrix which has ``1'' everywhere on the first  upper off-diagonal and all other entries  ``0'', and whose size will be specified, or clear, depending on the context. We call $N$ the standard indecomposable nilpotent matrix (of a given size). For example, if the size happens to be $3\times3$, then
$$ 
N = 
\left(\begin{array}{ccc} 
0 & 1 & 0    \\ 
 0&  0& 1   \\
 0& 0 & 0 \\
\end{array}\right).
$$

For block matrices, we will use the following convention, which generalizes multiplication of matrices by scalars: if $M$ is a block matrix with square blocks of size $l \times l$, and if $K$ is an $l \times l$ matrix, then $KM$ will denote the block matrix whose blocks consist of the blocks of $M$ each multiplied on the left by $K$. We define $MK$ similarly.

\subsection{Classification in case of eigenvalue in ground field}\label{eig value in ground field}

Given a $3k$-dimensional indecomposable sextuple $\mathcal{S}_\eta$, we call a frame basis 
$$
(u_1, ...,u_k, v_1,..., v_k, w_1,..., w_k)
$$ a {\bf Jordan} basis for $\mathcal{S}_\eta$ if  $(v_1, ..., v_k)$ is a Jordan basis for $\eta$ for some (unique) $\lambda \in \fieldk$, \text{i.e.} $\eta(v_1)=\lambda v_1$, $\eta(v_j)= v_{j-1}+\lambda v_j$  for $j>1$. Note that any Jordan basis for $\eta$ extends to one for $\mathcal{S}_\eta$.

\begin{thm}\label{eig}
Consider an  indecomposable  sextuple $\mathcal{S}_\eta$ of dimension $3k$ with underlying endomorphism $\eta$ having eigenvalue $\lambda \in \fieldk $; fix a Jordan basis for $\mathcal{S}_\eta$. Let $M \in \fieldk^{k\times k}$ be of the form
$$ 
\left(\begin{array}{ccccc} 
 & & \dots& 0 & (-1)^{j+1} \\ 
 & \dots  & \dots & \dots &  \dots \\
 \dots & 0  & 1  & 0 & \dots \\
0 &  -1 &  0& \dots &  \\
1 & 0 &  \dots  &    &   
\end{array}\right)
$$
\text{i.e.} the entries of $M$ satisfy $m_{ij}=(-1)^{j+1}$ for $i= k-j+1$, and $m_{ij}=0$ else.
\begin{enumerate}  
\item\label{eig_part_1} The sextuple is isomorphic to its dual if and only if   $\lambda=\frac{1}{2}$.
\end{enumerate} 
Now assume  $\lambda=\frac{1}{2}$.
\begin{enumerate} \setcounter{enumi}{1}
\item\label{eig_part_2}
The matrix 
$$
H_M =\left(\begin{array}{ccc} O&O&-M\\O&M&O\\-M&O&O \end{array}\right)
$$
defines a compatible form  $B$ which is symplectic if $k$ is even, symmetric if $k$ is odd.

\item\label{eig_part_3} Up to isomorphism and multiplication with a scalar, $B$ is the only such form. For even $k$  there is no symmetric compatible form, for odd $k$ no symplectic one.
\item\label{eig_part_4} For any  $0\neq c \in \fieldk$,  $cB$ and $B$ are isometric via an automorphism of the sextuple if and only if $c$ is a square in $\fieldk$.
\item\label{eig_part_5} A complete list of compatible symplectic resp. symmetric forms is given by the matrices $H_{MQ}$, where $Q  \in \fieldk^{k\times k}$ is of the form $\sum_{i=0}^{k-1} a_i N^i$, with  $a_0\neq 0$ and $a_i=0$ for $i$ even. \end{enumerate}
\end{thm} 

\begin{rmk}
Observe that the particular value $\lambda =\frac{1}{2}$
arises from the way the endomorphism is related to the
sextuple. Other ways would give a different (but also unique) value.

Also note that, setting $\zeta := \eta - \tfrac{1}{2}$, we have that $\eta = \tfrac{1}{2} + \zeta$ underlies an indecomposable self-dual framed sextuple if and only if $\zeta$ is similar to $- \zeta^*$. Furthermore, $\eta$ has an eigenvalue if and only if $\zeta$ does, and in the case of self-duality, the unique eigenvalue of $\zeta$ is $0$. The description ``$\eta = \tfrac{1}{2} + \zeta$''  is helpful to keep in mind in the proof below, and in the subsequent sections. 
\end{rmk}

\noindent \textbf{Proof of Theorem \ref{eig}}
With respect to $v_1, \ldots, v_k$, the matrix of $\eta$ is
$A=\lambda I +N$.  In particular, $A$ is similar to $A^t$ via the permutation matrix $P$ corresponding to
reversing the order of the basis.

\emph{Proof of \ref{eig_part_1} and \ref{eig_part_2}}. Assume that $\mathcal{S}_\eta$ is isomorphic to its dual.
Then by Proposition \ref{duals of endo-sextuples}, $I-A^t$ is similar to
$A$, and hence to $A^t$. Thus, $A$ is similar to $I-A$.
In particular, $I-A$ necessarily also has $\lambda$ as its unique eigenvalue. On the other hand, given an eigenvector $v$ for $A$, we have $(I-A)v = v -\lambda v= (1- \lambda) v$, which means that $1-\lambda$ is an eigenvalue of $I - A$. By indecomposablilty, $I-A$ can only have one eigenvalue, thus $\lambda  = 1 - \lambda$, \text{i.e.} $\lambda=\frac{1}{2}$.

Now, assume $\lambda=\frac{1}{2}$ and define $B$ by $H_M$.
It remains to show that $B$ is compatible,
i.e. that $H_M$ defines, with respect to the given basis and its dual,
an isomorphism from $\mathcal{S}_\eta$ onto its dual.
By Proposition \ref{duals of endo-sextuples} this means to show that $M$ defines an isomorphism $(A_2, \eta) \to (A^*_2, (\id - \eta)^*)$, \text{i.e.} $MAM^{-1}=(I-A)^t$. Here, observe that $M=DP$, where $P=P^{-1}$ has anti-diagonal entries $1$, zero else, and where $D = D^{-1}$ is diagonal with entries $d_{jj}=(-1)^{j+1}$. Now, $PAP=A^t =\frac{1}{2}I+N^t$ and $DN^tD=-N^t$, and hence $MAM^{-1}= \frac{1}{2}I +DN^tD = \frac{1}{2}I - N^t  = I-A^t$.

\emph{Proof of \ref{eig_part_3} -- \ref{eig_part_5}} 
By Remark \ref{loc}, the algebra $E$ of endomorphisms of $A_2$ commuting with $A=\frac{1}{2}I+N$ is local and $E = \fieldk \id \oplus \text{Rad} E$; Proposition \ref{endo to sextuple} establishes an isomorphism of $E$ onto the endomorphism algebra $E'$ of $\mathcal{S}_\eta$ which shows that $E'$ is local, too, with $E' = \fieldk \id \oplus \text{Rad} E'$. Thus, Lemma \ref{uniqueness up to a scalar} applies, proving uniqueness. The proof of 4. can be copied from that of Theorem \ref{unidis}. Moreover, units in the ring $E'$ are given by block matrices as in (\ref{block matrix}). Thus, due to \ref{eig_part_2}, any  compatible form is given by a matrix $H_{MQ}$ where $Q  \in \fieldk^{k\times k}$ is the coordinate matrix of an automorphism of $(A_2, \eta)$. In particular, this means that  $Q$ is of the form $\sum_{l=0}^k a_l N^l$ with  $a_0\neq 0$. Thus $MQ = \sum_{l=0}^k a_l MN^l$.  Note that $MN^l$ is  skew-symmetric if $l = k\mod 2$ and symmetric if $l \neq k\mod 2$:
\begin{align*}
(MN^{l})_{ij} &= \sum_p m_{i, p}(N^l)_{p,j} = m_{i, k-i+1} (N^l)_{k-i+1,j}  = (-1)^{k - i}  \quad  \mbox{ for }    i+j=l+k+1\\
(MN^{l})_{ij} & =0 \quad \mbox{ else},  
\end{align*} 
and so, for   $i+j=l+k+1$, 
$$ (MN^{l})^t_{ij}(MN^{l})_{ij} =  (MN^{l})_{ji}(MN^{l})_{ij} = (-1)^{k - j}(-1)^{k - i} = (-1)^{k-l - 1}. $$
It follows that $MQ = \sum_{i \text{ even}} a_i MN^i + \sum_{i \text{ odd}} a_i MN^i$ is the (unique) decomposition of $MQ$ into its symmetric and skew-symmetric parts: when $k$ is odd, the first summand is the symmetric part and the second summand skew-symmetric; when $k$ is even, the reverse is the case. The second summand is non-invertible, while the first summand is invertible (since $a_0 \neq 0$). Thus all compatible (non-degenerate) forms are of the form $H_{MQ}$, where $MQ = \sum_{i \text{ even}} a_i MN^i$. For $k$ odd they are symmetric; for $k$ even, they are skew-symmetric. 
\qed \\
 
\begin{ex}
Let $k = 4$, and consider the self-dual sextuple $\mathcal{S}_\eta$ with underlying indecomposable endomorphism having eigenvalue $\lambda = \tfrac{1}{2}$.  Then $\mathcal{S}_\eta$ admits only compatible symplectic forms. These are parametrised by the set $\{ (a_0, a_2) \in \fieldk^2 \mid a_0 \neq 0 \}$ and given, with respect to a fixed Jordan basis and its dual basis, by the matrices 
$$
H_A =\left(\begin{array}{ccc} O&O&-A\\O&A&O\\-A&O&O \end{array}\right),
$$
where
$$ 
A = 
\left(\begin{array}{ccccc} 
0  & 0  & 0 & -a_0  \\
 0 & 0  & a_0  & 0  \\
0 &  -a_0 &  0& -a_2   \\
a_0 & 0 &  a_2  &  0     
\end{array}\right).
$$
\end{ex}

\subsection{Classification over $\mathbb{R}$}

\begin{thm}\label{real}
Let $\mathcal{S}_\eta$  be an indecomposable sextuple over $\mathbb{R}$ of dimension $3k$ with underlying $\eta$ having no real eigenvalue. 

\begin{enumerate}
\item The sextuple admits an isomorphism onto its dual if and only if $\eta$ has, in $\mathbb{C}$, eigenvalues $\lambda =\frac{1}{2} \pm \sqrt{-1} r $, with real  $r > 0$ (both of geometric multiplicity $1$).

\item Assume the case of self-duality. Then $k=2l$ and both symmetric and symplectic compatible forms exist. Furthermore:
\begin{enumerate}
\item  
There is a frame basis $(u_i,v_i,w_i)$ such that the matrix $A$ of $\eta$ is in real Jordan normal form with respect to the $v_i$, and there are ``canonical'' compatible forms $B$ and $B'$, where ``$B' = \sqrt{-1}B$''. If $l$ is odd, $B$ is symmetric and $B'$ is symplectic; if $l$ is even, $B$ is symplectic and $B'$ is symmetric. 

With respect to the frame basis, $B$ is given by the matrix $H=H_M$
where $M \in \fieldk^{k \times k}$ has $2 \times 2$-blocks $M_{ij}$ with
$$
M_{ij} =(-1)^{j+1} \left(\begin{array}{cc}1&0\\0&1 \end{array}\right) \;\mbox{ for } i+j=l+1,\; M_{ij}=\left(\begin{array}{cc}0&0\\0&0 \end{array}\right) \mbox{ else } 
$$
and the matrix $H'$ of $B'$ is obtained from $H$ by setting $H' = \Im H$, where $\Im$ is the matrix 
$$
\Im =\left(\begin{array}{lr} 0&-1\\1&0 \end{array}\right).
$$
Multiplication with $\Im$ is understood in sense of ``scalar multiplication by `$\sqrt{-1}$' for $2 \times 2$ block matrices''. 
\item\label{part_b}All symmetric resp. symplectic compatible forms are given,
up to isomorphy and multiplication with $\pm 1$, respectively by
\[ 
H \mbox{  resp. } \Im H \mbox{ if } l \mbox{ is odd, } \quad 
 \Im H \mbox{  resp. }  H \mbox{ if } l \mbox{ is even. }
\]
In particular, multiplying with $\Im$ changes the symmetry of a compatible form. 
There is no automorphism which is an isometry from $H$ to $-H$ or from $\Im H$ to $- \Im H$.

\item\label{part_c} A complete list of compatible symmetric, respectively symplectic, forms is given by the matrices
$H_{MQ }$ resp. $\Im H_{MQ}$ where $M$ is as above and where $Q  \in \mathbb{R}^{k\times k}$ is of the form $\sum_{i=0}^{l-1} a_i N^{2i}$, $a_i \in \mathbb{R}$,  $a_0\neq 0$, and $a_i = 0$ for $i$ even.  
These are symmetric, respectively symplectic, if $l$ is odd; symplectic, respectively symmetric, if $l$ is even.
\end{enumerate}
\end{enumerate}
\end{thm}

\pf
With respect to a suitable basis of $A_2$, $\eta$ has coordinate matrix in real Jordan form 
\[
A=ZI+N^2, \quad  \mbox{ where } Z=
\left(\begin{array}{lr} a&-b\\b&a \end{array}\right)
\] 
and $a\pm \sqrt{-1} b $ are the complex eigenvalues of $\eta$ $(b \neq 0$ by hypothesis).
Here, $ZI$ is the block diagonal matrix with diagonal
blocks $Z$ and $N$ the standard indecomposable
nilpotent matrix. 
In view of Proposition \ref{duals of endo-sextuples},
self-duality means that $(I-A)^t$ is similar to $A$. 
Now, $I-A$ and $(I-A)^t$ have the same complex eigenvalues,
and these are also those of $A$, by similarity. Thus, $\lambda$
is an eigenvalue of $A$ if and only if so is $1-\lambda$ and so
both have the same real part which then must be $\frac{1}{2}$.
This proves that
$$
 Z=
\left(\begin{array}{lr} \frac{1}{2}&-r\\r&\frac{1}{2} \end{array}\right), r \neq 0$$ 
On the other hand, any such $Z$ satisfies $(I-Z)^t=Z$, and thus corresponds to a self-dual sextuple. 

Now, to prove that $H_S$ is a compatible form, we can mimic the proof  of \ref{eig_part_2} in Theorem \ref{eig}, replacing scalar matrix entries there by $2\times 2$-blocks here: $0$, $1$, and $-1$ are replaced by  zero, unit, and negative unit matrix, $\frac{1}{2}$ by $B$, and $m_{ij}$ by $M_{ij}$. Thus $H_S$ defines an isomorphism onto the dual. The further  statements   are then obvious.

\emph{Proof of} (\ref{part_b}).  Observe that matrices 
\[
Z=
\left(\begin{array}{lr} x&-y\\y&x \end{array}\right)
\]
belong to a subring $F$ of $\mathbb{R}^{2 \times 2}$ which is isomorphic to $\mathbb{C}$.  $F$ is the field in $\mathbb{R}^{2 \times 2}$ obtained by adjoining, for example, 
$\left(\begin{smallmatrix} 
0&-1 \\
1&0 \end{smallmatrix} \right)$ 
to the field 
$\{ \left(\begin{smallmatrix}
x&0 \\
0&x \end{smallmatrix}
\right) \mid x \in \mathbb{R} \} \subseteq \mathbb{R}^{2 \times 2}$.
For any $Z \in F$, let $ZI$ denote the block diagonal matrix in $\mathbb{R}^{k \times k}$ with diagonal blocks $Z$. Now, $C$ is in the (coordinatized) endomorphism algebra $E$ of $(A_2,\eta)$ if and only if it commutes with the matrix $A$ of $\eta$. By Proposition \ref{endo alg}  below, $E$ consists of the matrices $C$ of the form
\[ 
C=\sum_{i=0}^{l-1} Z_iN^{2i},\quad Z_i \in F. 
 \] 
By Proposition \ref{endo to sextuple}, $\mathcal{S}_\eta$ has endomorphism algebra consisting of the block-diagonal matrices $CI$, with $C \in E$ (\text{i.e.} there are three diagonal blocks, each one a copy of a given $C \in E$). In particular, these matrices commute with $H_M$. Now the proof of Lemma \ref{uniqueness up to a scalar}  generalizes\footnote{For a more thorough discussion, see Subsection \ref{uniq:gen}.} to yield uniqueness
of compatible forms up to isomorphism and multiplication with $Z \in F$; that is,  up to isomorphism, compatible forms are of the form $ZH_M = H_{ZM}$ with 
\[
Z= 
\left(\begin{array}{cc} 
x&-y \\
y&x \end{array}
\right) 
\in F \backslash \{0\}.
 \]
If $y=0$, then $ZM$ has blocks 
\[
\left(\begin{array}{cc} 
x&0 \\
0&x \end{array}\right), \;
\left(\begin{array}{cc} 
-x&0\\
0&-x 
\end{array}\right),\;
\left(\begin{array}{cc} 
x&0\\
0&x 
\end{array}\right) 
\ldots
\]
along its anti-diagonal, and so $ZH_M = xH_M$. In this case, the scaling map $1/\sqrt{\vert x \vert} \cdot \text{Id}$ is an isometric isomorphism to from $xH_M$ to $\text{sign}(x)H_M$. 

If $x=0$, we have blocks 
\[\left(\begin{array}{cc} 
0&-y \\ y&0 \end{array}\right),\;
\left(\begin{array}{cc} 0&y\\-y&0 \end{array}\right),\;
\left(\begin{array}{cc} 0&-y\\y&0 \end{array}\right) \ldots
\]
on the anti-diagonal of $M$, and so $ZH_M = y \Im H_M$. The scaling map $1/\sqrt{\vert y \vert} \cdot \text{Id}$ is an isometric isomorphism to from $y \Im H_M$ to $\text{sign}(y) \Im H_M$. 

If $x\neq 0$ and $y \neq 0$ then we have
\[\left(\begin{array}{cc} x&-y\\y&x \end{array}\right),\;
\left(\begin{array}{cc} -x&y\\-y&-x \end{array}\right),\;
\left(\begin{array}{cc} x&-y\\y&x \end{array}\right) \ldots
\]
on the diagonal of $ZM$, and neither symmetry nor skew-symmetry. 

It remains to show that there is no isometric isomorphism from $H_M$ to $-H_M$ (such a map would also give an isometric isomorphism between $\Im H$ and $H$).
Assume that $f$ were such an automorphism and $C$ the matrix of $f |_{A_2}$. From the above description of $C\in E$ one reads off that $C$ is upper block triangular with 
\[ 
c_{11}= c_{2k-1,2k-1},\quad  c_{21}= c_{2k,2k-1}.  
\] 
On the other hand, by inspection of $H_M$ we have
\[ 
B(v_1,v_{2k-1})=B(v_2,v_{2k})=1,\quad
B(v_1,v_i)=0 \mbox{ for } i\neq 2k-1,\quad B(v_2,v_i)=0 \mbox{ for } i\neq 2k.   
\] 
Thus, one would get the contradiction
\begin{align*}
-1&= -B(v_1,v_{2k-1})=B(f v_1, f v_{2k-1})=
B(c_{11}v_1+ c_{21}v_2,\,\sum_{i=1}^{2k} c_{i,2k-1}v_i) \\
&=B(c_{11}v_1,c_{1,2k-1} v_{2k-1})+ 
B(c_{21}v_2,c_{2k,2k-1} v_{2k}) =c_{11}^2+c_{21}^2.
\end{align*}
 
\emph{Proof of} (\ref{part_c}).
As in the proof of \text{4.} in Theorem \ref{eig} one obtains matrices as stated above, but initially with $a_i\in F$. To have symmetry or skew-symmetry in $MQ$, though, all $a_i$ have to be in $\mathbb{R}I$ or  $\mathbb{R} \Im I$. Conversely, this grants symmetry resp. skew symmetry, depending on the parity of $l$.
\qed \\

\begin{ex}
Let $k=6$, and let $\eta$ be an indecomposable endomorphism over $\mathbb{R}$, with complex eigenvalues $\tfrac{1}{2} \pm \sqrt{-1}r$, $r > 0$. The corresponding sextuple $\mathcal{S}_\eta$ is self-dual (and it lives in ambient dimension $3k = 18$); compatible symmetric resp. symplectic forms are given by matrices $H_M$, where $M$ is of the form
\[
\left(\begin{array}{cccccc} 
0 & 0& 0&0& a_0&0\\
0 & 0 & 0&0&0&a_0\\
0&0&-a_0&0 & 0 & 0 \\
0&0&0&-a_0 & 0 & 0\\
a_0&0&0&0 & a_1 & 0\\
0&a_0&0&0 & 0 & a_1\\
\end{array}\right)
\quad \text{resp.} \quad
\left(\begin{array}{cccccc} 
0&0&0&0&0 &-a_0\\
0&0&0&0&a_0& 0 \\
0&0&0 &a_0 & 0&0\\
0&0&-a_0& 0 & 0&0 \\
0 &-a_0&0&0 & 0&-a_1\\
a_0&0&0&0 & a_1&0\\
\end{array}\right)
\]
for $(a_0, a_1) \in \mathbb{R}^2$, $a_0 \neq 0$. 
\end{ex}

\subsection{Hamiltonian vector fields over $\mathbb{R}$ from non-split framed sextuples}

In this section we return to the connection, discussed in Section \ref{Hamiltonian vector fields}, between linear hamiltonian vector fields and isotropic triples. We show here that each non-split isotropic triple (over $\mathbb{R}$) has an associated linear hamiltonian vector field. 

Let $S_\eta$ be an indecomposable self-dual continuous-type sextuple over $\mathbb{R}$.  We know that the spectrum in $\mathbb{C}$ of $\eta$ must be $\{\tfrac{1}{2} \pm \sqrt{-1}r\}$ for some value of $r \geq 0$. Suppose that we make $S_\eta$ into an isotropic triple by choosing a frame basis and choosing a compatible symplectic form given by a matrix $H_{MQ}$ as in Theorem \ref{eig} or Theorem \ref{real} (depending on whether $\eta$ has the eigenvalue $\tfrac{1}{2}$ or not). Set $T := MQ$ and let $A$ denote the coordinate matrix of $\eta$ with respect to the frame basis. Recall that when $H_{T}$ is skew-symmetric, as we have assumed, then so is $T$. From the normal form given in Lemma \ref{augm frames}, part \ref{augm frames part 3}, we can assume that $S_\eta$ has the form
\begin{equation}\label{normal again}
\begin{array}{ll}
V = U \times U \times U &   \\
I_1 = U \times 0 \times 0 & \quad C_1 = U \times U \times 0 \\
I_2 = 0 \times 0 \times U & \quad C_2 = 0 \times U \times U \\
I_3 = \{(-\eta x,x,-x)\mid x \in U\} & \quad C_3 = \{(x,-x,0)\mid x\in U\} + I_3,
\end{array}
\end{equation}
and from Proposition \ref{duals of endo-sextuples} it follows that $T$ is the coordinate matrix of an isomorphism $(U,\eta) \rightarrow (U^*, 1- \eta^*)$. In other words, $TA = (1-A^t)T$. 

We claim now that an isotropic triple whose underlying sextuple is (\ref{normal again}) satisfies the hypotheses of Proposition \ref{non-exceptional hamiltonian iso triples}, which gives sufficient conditions for constructing an associated hamiltonian vector field. To verify the hypotheses, note first that clearly $V = I_1 \oplus I_2 \oplus I_3$. Second, we must check that $I_I + I_j$ is a symplectic subspace for all $i \neq j$. For this, we check that 
$$
(I_i + I_j) \cap (I_i + I_j)^\perp = (I_i + I_j) \cap C_i \cap C_j = 0.
$$
For $i=1, j=2$, 
\begin{align*}
(I_1 + I_2) \cap C_1 \cap C_2 &= (U \times 0 \times U) \cap (0 \times U \times 0) = 0.
\end{align*}
For $i=3, j=1$, 
\begin{align*}
(I_3 + I_1) \cap C_3 \cap C_1 &= \{ (-\eta x + y, x, -x) \mid x, y \in U \} \cap \{ (z, -z, 0) \mid y \in U \} = 0, 
\end{align*}
since an element $(-\eta x + y, x, -x)$ of this intersection would need to satisfy $-x = 0$, and hence $-\eta x =0$. But then $(y, 0, 0) = (z,-z,0)$ holds for some $z \in U$ only if $y =0$. 

Finally, for $i=2, j=3$, 
\begin{align*}
(I_2 + I_3) \cap C_2 \cap C_3 &= \{ (-\eta x, x, y) \mid x, y \in U \} \cap \{ (0, -\eta z + z, -z) \mid z \in U \} \\
	&= \{ (0, x, y) \mid x, y \in U, \eta x = 0, x = \eta y - y \} = 0.
\end{align*}
Indeed, $\eta x = 0$ implies that $x = 0$, since $\eta$ is invertible, and so $\eta y = y$ must hold. By assumption $\eta$ cannot have eigenvalue $1$, so also $y = 0$. 

Now we will follow the proof of Proposition \ref{non-exceptional hamiltonian iso triples} in order to find the hamiltonian vector field associated to (\ref{normal again}). We have the symplectic decomposition $V = (I_1 \oplus I_2) \oplus (I_1 \oplus I_2)^\perp = (U \times 0 \times U) \oplus (0 \times U \times 0)$.  $I_3$ is the graph of the map $g : (I_1 \oplus I_2)^\perp \rightarrow I_1 \oplus I_2$ given by 
$$
0 \times U \times 0 \rightarrow U \times 0 \times U, \quad (0,x,0) \longmapsto (-\eta x, 0, -x),
$$
so the image of $g$ is the graph of the map $I_2 \rightarrow I_1, (0,0,x) \mapsto (\eta x, 0,0)$. With respect to the chosen frame basis, the coordinate matrix of this map is $A$, and the matrix of the identification of $I_1$ with $I_2^*$ is $T$. Thus the image of $g$ corresponds to a map $f: I_2 \rightarrow I_2^*$ whose coordinate matrix is $TA$. Using $TA = (1-A^t)T$ and that $T$ is skew-symmetric, we find that the antisymmetric part $f_a$ is given by 
$$
\tfrac{1}{2}(TA - A^tT^t) = \tfrac{1}{2}(TA - T^t(1-A)) = \tfrac{1}{2}T(A  + (1-A)) = \tfrac{1}{2}T,
$$
and for the symmetric part $f_s$ 
$$
\tfrac{1}{2}(TA + A^tT^t) = \tfrac{1}{2}(TA + T^t(1-A)) = \tfrac{1}{2}T(A  - (1-A)) = \tfrac{1}{2}T(2A-1).
$$
Thus we obtain the hamiltonian vector field $X = f_a^{-1}f_s$ given in coordinates by the matrix $2A-1$. 
We have proved the following:

\begin{prop}
Non-split isotropic triples satisfy the hypotheses of Proposition \ref{non-exceptional hamiltonian iso triples}: we can construct an associated hamiltonian vector field. 

If $\varphi$ is such a triple, we may assume that its underlying sextuple is a framed sextuple $S_\eta$, built from an endomorphism $(U, \eta)$. In this case, the symplectic form $\omega$ of $\varphi$ induces an isomorphism $U \rightarrow U^*$ which defines a symplectic form $\omega_U$ on $U$. The hamiltonian vector field associated to $\varphi$ is $(U, \omega_U, 2\eta - \text{id})$. 
\end{prop}

\begin{rmk}
Normal forms for indecomposable linear hamiltonian vector fields over $\mathbb{R}$ are given, for example, in \cite{laub-meyer}. These come in both ``split'' and ``non-split'' types, and are labeled by their (complex) eigenvalues. For a non-split isotropic triple as above with underlying indecomposable endomorphism $\eta$, the associated hamiltonian vector field $2\eta - \text{id}$ is also indecomposable and non-split. The possible such $\eta$ are parametrized by $\tfrac{1}{2} \pm \sqrt{-1}\cdot r$ with $r \in [ 0,\infty )$; the corresponding hamiltonian vector fields $2\eta - \text{id}$ correspond in \cite{laub-meyer} to the non-split ones labeled by complex eigenvalues $\pm \sqrt{-1} \cdot \nu$ (setting $\nu= 2r$ in order to use their notation), where $\nu \in [ 0,\infty )$. 
\end{rmk}

\begin{ex}
Let $k=4$, and let $\eta$ be an indecomposable endomorphism over $\mathbb{R}$ with complex eigenvalue $\tfrac{1}{2}$, \text{i.e.} with respect to a Jordan basis $\eta$ has the coordinate matrix
$$
A=
\left(\begin{array}{cccc} 
1/2 &1 & 0 & 0 \\
0&1/2 & 1 & 0\\
0 & 0&1/2 & 1 \\
0 & 0& 0 & 1/2 \\
\end{array}\right).
$$
We can make the corresponding self-dual sextuple $\mathcal{S}_\eta$ into an isotropic triple by choosing the compatible symplectic form given, with respect to a Jordan frame basis, by the matrix $H_T$, with 
$$
T = 
\left(\begin{array}{cccc}  
0 & 0 & 0 &1\\
0 & 0 & -1&0 \\
0 &1 & 0 & 0 \\
-1&0 & 0 & 0 \\
\end{array}\right).
$$
The associated hamiltonian vector field on the symplectic space $(\mathbb{R}^2, T)$ is given, with respect to a Jordan basis for $\eta$, by the matrix 
$$
\left(\begin{array}{cccc} 
0 &1 & 0 & 0 \\
0& 0 & 1 & 0\\
0 & 0& 0 & 1 \\
0 & 0& 0 & 0 \\
\end{array}\right).
$$
\end{ex}
\begin{ex}
Let $k=4$, and let $\eta$ be an indecomposable endomorphism over $\mathbb{R}$ with complex eigenvalue $\tfrac{1}{2} \pm \sqrt{-1}\tfrac{1}{2}\nu$, where $\nu >0$, \text{i.e.} with respect to a (real) Jordan basis $\eta$ has the coordinate matrix
$$
A=\frac{1}{2}
\left(\begin{array}{cccc} 
1&-\nu & 1 & 0\\
\nu&1 & 0 & 1 \\
0 & 0& 1&-\nu\\
0 & 0& \nu&1 \\
\end{array}\right).
$$
Again choosing the compatible symplectic form for $\mathcal{S}_\eta$ given by $H_T$ with 
$$
T = 
\left(\begin{array}{cccc}  
0 & 0 & 0 &1\\
0 & 0 & -1&0 \\
0 &1 & 0 & 0 \\
-1&0 & 0 & 0 \\
\end{array}\right),
$$
we obtain the associated hamiltonian vector field on $(\mathbb{R}^2, T)$ given by
$$
\left(\begin{array}{cccc} 
0&-\nu & 1 & 0\\
\nu&0 & 0 & 1 \\
0 & 0& 0&-\nu\\
0 & 0& \nu&0 \\
\end{array}\right).
$$
\end{ex}

\section{Continuous non-split isotropic triples over perfect fields}\label{continuous perfect}

We continue working with a ground field $\fieldk$ with $\text{char}(\fieldk) \neq 2$, and additionally we assume that $\fieldk$ is perfect\footnote{A field $\fieldk$ is perfect if every algebraic extension of $\fieldk$ is separable. Examples of perfect fields include all finite fields, and all fields of characteristic zero.}, in order that we may use the normal forms discussed in Proposition \ref{fnf} below. No other assumptions are made on $\fieldk$.  

According to Proposition \ref{duals of endo-sextuples}, isomorphisms of a sextuple $\mathcal{S}_\eta$ onto its dual arise from matrices $M$ such that $M^{-1}AM=(I-A)^t$ where $A$ is the matrix of $\eta$ with respect to some basis; and the induced compatible form is (skew-)symmetric if and only if $M$ is. As indicated by the results over the complex or real number field, the crucial case is when $\eta$ has irreducible minimal polynomial $q(x)$. We deal with this case first.

\subsection{Irreducible characteristic polynomial}\label{cay}

Let $q(x) =\sum_{i=0}^l a_ix^i \in  \fieldk[x]$ be an irreducible polynomial, with $a_l=1$. Adjunction to $\fieldk$ of a zero $\lambda$ of $q$ yields an extension field  $\fieldk(\lambda) = \fieldk[\lambda] \cong \fieldk[x]/q(x)$ which as a $\fieldk$-vector space has  basis $1,\lambda,\lambda^2, \ldots, \lambda^{l-1}$. With respect to this basis, the $\fieldk$-linear map $m_\lambda$ defined by $m_\lambda(r) :=\lambda r$ 
has as its coordinate matrix the {\bf Frobenius matrix} $N^t+C$
$$
\left(\begin{array}{ccccc} 
0& 0 & \dots & 0  & -a_0\\
1 & 0 & \ddots  & 0 & -a_1 \\
0 & 1 & \ddots & &  \vdots \\
\vdots & & \ddots& 0 & -a_{l-2} \\
0 & \dots & 0 & 1 & -a_{l-1} \\
\end{array}\right),
$$
where $C = C_q$ denotes the matrix whose last column is $-a_0,-a_1,-a_2, \ldots, -a_{l-1}$ and has all other entries zero.

In particular, by the Cayley-Hamilton theorem the above applies to any endomorphism $\eta$ of an vector space $U$ having characteristic polynomial $q(x)$, where $\eta$ plays the role of $\lambda$ above\footnote{Recall that for indecomposable endomorphisms, the characteristic and minimal polynomials coincide. Also note that, when the characteristic polynomial is irreducible, the corresponding endomorphism is necessarily indecomposable.}. Here we view $\fieldk(\eta)$ as a subring of $\End(U)$ with subfield $\{a {\sf id}\mid a \in \fieldk\} \cong \fieldk$. In this situation, any $v_1 \neq 0$ extends (uniquely) to a basis $\{v_1, \eta v_1, ...,\eta^{l-1}v_1\}$ such that the corresponding coordinate matrix of $\eta$ is the Frobenius matrix of $q(x)$. In other words, if $A \in \fieldk^{l\times l}$ has $q(x)= \det(xI-A)$ irreducible, then $ \fieldk(A) = \fieldk [A]$ is a subring of $\fieldk^{l \times l}$ which is an extension field of $\{aI\mid a \in \fieldk\} \cong \fieldk$ with primitive element $A$, a zero of $q(x)$.

\begin{lemma}\label{extensions and involution} 
Let $q(x) \in \fieldk [x]$ be an irreducible monic polynomial with $\deg q = l > 1$. Let $\lambda$ be a zero of $q$ in some extension field of $\fieldk$, and set $\bf{E} = \fieldk(\lambda)$. \

Suppose $q(1- \lambda) = 0$. Then  there exists $\mu \in \bf{E}$ and an automorphism $g$ of $\bf{E}$ over $\fieldk$ such that $\fieldk(\mu) = \bf{E}$ and $g(\mu) = -\mu$. 
This implies that:
\begin{enumerate}
\item Only even powers of $x$ occur in the minimal polynomial $r(x)$ of $\mu$ over $\fieldk$; in particular $\deg r = l$ must be even.
\item With respect to the basis $\{ 1, \mu, ..., \mu^{l-1}\}$ of $\bf{E}$ over $\fieldk$, $g$ has diagonal coordinate matrix $D$, with entries $d_{ii} = (-1)^{i}$, for $i = 0,..., l-1$.
\item $g$ is a $\fieldk$-isomorphism $({\bf{E}}, m_{- \mu}) \rightarrow ({\bf{E}}, m_\mu)$, i.e. $m_\mu g = g m_{- \mu}$.
\end{enumerate}  
\end{lemma}

\pf 
Since $q(\lambda)=q(1-\lambda)=0$, there exists an automorphism $g$ of $\bf{E} = \fieldk(\lambda)$ over $\fieldk$ such that $g(\lambda)=1-\lambda$. Setting $\mu:=\lambda -\frac{1}{2}$, we have $\lambda= \frac{1}{2} +\mu$ and $1-\lambda= \frac{1}{2} -\mu$; in particular $\fieldk(\lambda)=\fieldk(\mu)$ and $g(\mu)=-\mu$. Note that $-\mu \neq \mu$, since $-\mu = \mu$ would imply $\lambda = 1/2 \in \fieldk$, and hence that $[\fieldk (\lambda): \fieldk] = \text{deg}(q) = l = 1$, contrary to the assumption $l>1$. 
For the minimal polynomial 
$$
r(x) = a_lx^l + a_{l-1}x^{l-1} + ... + a_1x + a_0 \quad \quad (\text{with } a_l=1)
$$
of $\mu$ over $\fieldk$ it follows that $r(-\mu)= r(g (\mu)) =  g (r(\mu)) = 0$, using that $g$ fixes $\fieldk$.  This implies that only even powers of $x$ occur in $r(x)$, by comparing the coefficients of $r(-\mu)$ and $r(\mu)$. Indeed,  
$$
  \sum_{k \text{ odd }} a_k \mu^k = \frac{1}{2}(r(\mu)- r(-\mu))=0 
$$ 
since $r(\mu) = r(-\mu)=0$, so $\mu$ is a zero of $\delta(x):= \sum_{k \text{ odd }} a_k x^k$, and hence $r(x) \mid \delta(x)$. If $r(x)$ were to have odd degree, then $\delta(x)$ and $r(x)$ would have the same degree and hence would be equal (since they are both monic). But this would contradict the irreducibility of $r(x)$, because for a polynomial with only odd powers of $x$, one can always factor out the polynomial $p(x)=x$.  Therefore it must be that $\delta(x) \equiv 0$, \text{i.e.} $a_k = 0$ for all odd $k$. 

The statement about the coordinate matrix of $g$ is clear. For point \text{3.}, we check on the basis elements $\mu^i$ of $\bf{E}$:
$$  
m_\mu (g (\mu^i)) = m_\mu ((-1)^{i} \mu^{i}) = (-1)^i \mu^{i+1} = -(-1)^{i+1}\mu^{i+1} = g( m_{-\mu} (\mu^i)). 
$$  
\qed \\


\begin{prop}\label{sim}
Let $0\neq A\in \fieldk^{l \times l}$ such that $q(x)=\det(xI-A)$
is irreducible. Then the following are equivalent:
\begin{enumerate}
\item $A$ and $(I-A)^t$ are similar to each other
\item $A$ and $I-A$ are similar to each other
\item  $q(\lambda)=q(1-\lambda)=0$ for some $\lambda$ in some extension field of $\fieldk$.
\item  For any $\lambda$ in any extension field of $\fieldk$: if $q(\lambda)=0$, then $q(1-\lambda)=0$.
\end{enumerate}
Suppose now that any of these equivalent conditions holds. Then $\deg q(x)$ is even. Moreover,  setting $\tilde A =A-\frac{1}{2}I$, the minimal polynomial $r(x)$ of $\tilde A$ satsifies $r(-\tilde A)=0$ and, for any invertible $T$,  
$$TAT^{-1}=(I-A)^t \ \Leftrightarrow \ T\tilde AT^{-1}=-\tilde A^t.$$
\end{prop}

\pf
That \text{1.} and  \text{2.} are equivalent follows from the fact that a matrix and its transpose are always similar, and that similarity is a transitive relation.  \text{2.} implies that $q(I-A)=0$, and thus that \text{3.} holds in the field $\fieldk(A) = \fieldk [A] \in \fieldk^{l\times l}$. 
To see that \text{3.} implies \text{4.}, it is enough to consider extensions $\fieldk(\lambda)$ of $\fieldk$ such that $q(\lambda)=0$; these are all isomorphic via isomorphisms fixing $\fieldk$, and hence in each of them also the equation $q(1- \lambda)= 0$ is satisfied. Finally we show that \text{4.} implies \text{2.}. 
Since $q(A)=0$, \text{4.} implies $q(I-A)=0$ and thus that $q(x)$ is also the unique elementary divisor of $I-A$. The multiset of elementary divisors is a complete invariant for similarity, so \text{2.} follows. That $\deg q$ is even and that $r(x)$ satisfies $r(-\tilde A)$ follows from Lemma \ref{extensions and involution}. 
\qed \\


\begin{rmk}\label{non-degeracy}
In the subsequent, we'll use the following fact. Let $\bf{E} \supseteq \fieldk$ be a finite field extension, so $\bf{E}$ is a finite-dimensional vector space over $\fieldk$, and let $\tau : \bf{E} \rightarrow \fieldk$ be a non-zero $\fieldk$-linear map. Then the bilinear form on $\bf{E}$ defined by 
$$
(u,v) \longmapsto \tau(uv)
$$
is non-degenerate. Indeed, if for fixed $v \in \mathbf{E} \backslash \{ 0 \}$ we have $\tau(uv) = 0$ for all $u \in \bf{E}$, then it would follow that $\tau(u') = 0$ for all $u' \in \bf{E}$, since any $u'$ can be written as $u'v^{-1}v$. This would contradict the assumption that $\tau \neq 0$. 
\end{rmk}

\begin{prop}\label{mc}
Let $l > 0$, and let $0\neq A \in \fieldk^{l\times l}$ be a  Frobenius matrix with irreducible characteristic polynomial $r(x)$ such that $r(-A)=0$. Then there are skew-symmetric as well as symmetric invertible $T\in \fieldk^{l\times l} $ such that $TAT^{-1}=-A^t$.
\end{prop}

\pf
We may restate the task as follows. We are given $V=\fieldk(\mu)$, $\dim V=l$ even, with irreducible monic $r(x)=\sum_{i=0}^{l} a_ix^i$ such that $r(\mu)=r(-\mu)=0$, in particular $a_i=0$ for odd $i$. Let $m_\mu(v)=\mu v$ for $v \in V$. The task is to find an isomorphism $\beta:V \to V^*$ such that 
\begin{equation}\label{compatible form beta}
 \beta m_{\mu} = -m_{\mu}^* \beta, 
\end{equation}
and such that the matrix $T$ of $\beta$ with respect to the basis $1,\mu, \ldots, \mu^{l-1}$ and its dual basis $1^*,\mu^*, \ldots$ is skew-symmetric resp. symmetric. Because of these bases, within the scope of this proof we use indices which always live between $0$ and $l-1$.

We consider the skew-symmetric task first. From  \cite{qu-sc-sc:quadratic}, page 276, we know there exists symmetric $S$
such that $SAS^{-1}=A^t$. Namely, one has the linear form $\tau \in V^*$ defined on the basis $1, \mu, \dots, \mu^{l-1}$ by
\[ 
\tau(\mu^i)= 1 \ \text{ if } i=l-1, \quad \tau(\mu^i)=0 \text{ else }
\]
and we take $S$ to be the matrix of the non-degenerate symmetric bilinear form $(u,v) \mapsto \tau(uv)$, \text{i.e.} the entries of $S$ are 
\[ 
s_{ij}= \tau(\mu^{i+j}) \quad \quad 0 \leq i, j \leq l-1 .
\]
Let $V_0$ be the $\fieldk$-subspace of $V$ spanned by $\{ \mu^i  \mid 0\leq i <l, \  i \text{ even }\}$;  let $V_1$ be the $\fieldk$-subspace spanned by $\{ \mu^i \mid 0 \leq i <l, \  i \text{ odd }\}$.

\begin{claim} For any $i \in \mathbb{N}$,
$\mu^i \in V_0$ if $i$ is even, $\mu^i\in V_1$ if $i$ is odd. In particular $\tau(\mu^i) =0$ for all even $i \in \mathbb{N}$. 
\end{claim}

\noindent
\emph{Proof of the claim}. We use induction. Clearly, the claim holds for $i<l$. Also, note that $\mu^l=-\sum_{j<l}a_j\mu^j$ where $a_j=0$ if $j$ is odd, and so $\mu^l\in V_0$. In particular, we see that $\mu V_0 \subseteq V_1$ and $\mu V_1 \subseteq V_0$. Now let $i>l$. If $i$ is odd, then by the induction hypothesis $\mu^{i-1} \in V_0$, and so $\mu^i = \mu \mu^{i-1} \in \mu V_0 \subseteq V_1$. The analogous argument applies for $i$ even.   
\qed \\

It follows that $S$ has an  following properties: using indices $i, j \in \{0,\ldots, l-1\}$, the entries $s_{ij}$ of $S$ only depend on $i+j$, with $s_{ij} =0$ for $i + j \leq l-2$ or $i+ j$ even, and $s_{ij}=1$ for $i+j = l-1$ . 

Now consider $T:=SD$, where $D$ is the diagonal matrix defined in Lemma \ref{extensions and involution}, \text{i.e.} the diagonal entries are $d_{ii}= (-1)^i$ for $i = 0,..., l-1$. Thus $T$ has entries $t_{ij}=s_{ij}$ if $j$ is even, $t_{ij}=-s_{ij}$
if $j$ is odd. In particular then, $t_{ij}=0$
if $i+j$ is even (since $s_{ij} = 0$ if $i+j$ is even) and $t_{ij}=-t_{ji}$ if $i+j$ is odd (since $i$ and $j$ must have different parity when $i+j$ is odd). So $T$ is skew symmetric.

With respect to the basis $1,\mu, \mu^2, \ldots $ and its dual basis,
$S$ provides (according to \cite{qu-sc-sc:quadratic} \text{Prop.} 4.4) an isomorphism
$(V, m_\mu) \to (V^*, m_{\mu}^*)$,  while 
by Lemma \ref{extensions and involution} $D$ provides an isomorphism
$(V, -m_{\mu}) \to (V, m_{\mu})$. Thus $T=SD$ gives an
isomorphism $(V, -m_{\mu}) \to (V^*, m_{\mu}^*)$, as desired.

\medskip

With the skew-symmetric case done, we now turn to the symmetric case. Recall that we want to find $\beta$ such that (\ref{compatible form beta}) holds; for this we work with the coordinate matrix $T$ of $\beta$ with respect to the basis $1,\mu, \ldots, \mu^{l-1}$ and its dual basis. 

For convenience, we write $T=(t_{ij})_{i,j= 0,\ldots ,l-1}$.   
Summations will be for $0,\ldots, l-1$ unless stated otherwise. 

\begin{claim}\label{claim for commuting}
Equation (\ref{compatible form beta}) holds for $T$ if and only if all of the following are satisfied:
\begin{enumerate}
\item $ t_{k,j+1} =-t_{k+1,j} \mbox{ for } j,k <l-1$
\item $\sum_h-a_ht_{kh} =-t_{k+1,l-1}   \mbox{ for } j=l-1,\; k<l-1
$
\item $ t_{l-1,j+1}=\sum_h a_ht_{hj} \mbox{ for } j<l-1,\; k=l-1
$
\item $ 
-\sum_h a_ht_{l-1,h} = \sum _h a_ht_{h,l-1} \mbox{ for } j=k=l-1$
\end{enumerate}
\end{claim}

\noindent
\emph{Proof of Claim.}
Let $j,k\leq l-1$. The matrix entries corrsponding to $\beta \hat{\mu}$ and $-\hat \mu^*\beta$, respectively,  are
\begin{align*} 
x_{jk}:=&(\beta \hat{\mu}(\mu^j))(\mu^k) =  (\beta(\mu^{j+1}))(\mu^k) \\
y_{jk}:=&(-\hat{\mu} \beta(\mu^j))(\mu^k)= -\sum_i t_{ij} \hat{\mu}^*(\mu^{i*}) (\mu^{k})=  -\sum_i t_{ij}\mu^{i*}(\hat{\mu}(\mu^{k})) =  -\sum_i t_{ij}\mu^{i*}(\mu^{k+1}),
\end{align*}
and  (\ref{compatible form beta}) holds if and only if $x_{jk}=y_{jk}$ for all $j,k\leq l-1$. Now,
\begin{align*}
x_{jk}&= \sum_i t_{i,j+1}\mu^{i*}(\mu^k)= t_{k,j+1} \quad \mbox{ if } j<l-1 \\
y_{jk} &=  -\sum_i t_{ij}\mu^{i*}(\mu^{k+1}) = -t_{k+1,j} \quad \mbox{ if } k<l-1.   \\
x_{l-1,k} &=\beta(\mu^l) (\mu^k) =\beta(\sum_h -a_h\mu^h )(\mu^k) =\sum_ {i,h} -a_ht_{ih}\mu^{i*}(\mu^k)=\sum_h-a_ht_{kh} \\
y_{j,l-1} &=  -\sum_i t_{ij}\mu^{i*}(\mu^{l}) = -\sum_i t_{ij}\mu^{i*}(\sum_h -a_h\mu^h) = \sum_{i,h} t_{ij}a_h\mu^{i*}(\mu^h) =\sum_h t_{hj}a_h
\end{align*}
Thus, that $x_{jk}=y_{jk}$ for all $j$, $k$ is equivalent to the equations stated in the claim.
\qed \\

Define $\tilde \tau \in V^*$ by
\[ 
\tilde \tau (\mu^i)=1 \mbox{ if } i=l-2, \quad \tilde \tau (\mu^i) = 0 \text{ else }
\] 
(note that we use that $l > 0$).
By Remark \ref{non-degeracy}, the matrix $\tilde S$ whose entries are $\tilde s_{ij} = \tilde \tau (\mu^{i+j})$ is invertible. 
Note that  $\tilde \tau (\mu^k)=0$ for odd $k \in \mathbb{N}$, since $\tilde \tau \vert_{V_1} = 0$. 
In particular, $(-1)^{i+1}\tilde \tau(\mu^{i+j})= (-1)^{j+1} \tilde \tau (\mu^{i+j})$ for all $0 \leq i, j \leq l-1$, since both sides are zero when $i$ and $j$ have different parity. 
Thus, if we define 
\[ 
t_{ij} := (-1)^{i+1}\tilde \tau(\mu^{i+j})= (-1)^{j+1} \tilde \tau (\mu^{i+j}) 
\]
we obtain a symmetric matrix $T$. Note that $T = - \tilde S D$, so $T$ is invertible. The entries $t_{ij}$ only depend on $i+j$, up to sign. For for $i+j < l-2$ they are zero; for $i+j= l-2$, $t_{ij} = (-1)^{i}= (-1)^{j}$.

It remains now only to check that this $T$ fulfills the conditions  
of Claim \ref{claim for commuting}: 

\begin{enumerate}
\item $ t_{k,j+1} = (-1)^{k+1}\tilde \tau(\mu^{k+j+1}) =-t_{k+1,j} \mbox{ for } j,k <l-1$
\item $t_{k+1,l-1}= (-1)^k\tilde \tau(\mu^{k+l})= (-1)^k \sum_h -a_h \tilde \tau(\mu^{k+h}) =  \sum_h a_h(-1)^{k+1}\tilde \tau(\mu^{k+h}) = \sum_ha_ht_{kh}  \\  \mbox{ for } j=l-1,\; k<l-1 $
\item $ t_{l-1,j+1}= (-1)^{j} \tilde \tau(\mu^{j+l}) = (-1)^j \sum_h -a_h \tilde \tau(\mu^{j+h}) = \sum_h a_h (-1)^{j+1}\tilde \tau(\mu^{j+h}) =\sum_h a_ht_{hj}\\ \mbox{ for } j<l-1,\; k=l-1$
\item  $ -\sum_h a_ht_{l-1,h} =-\sum_{h\; \mbox{\tiny  even}} a_h t_{l-1,h} =0$ since $l-1+h$ is odd for even $h$. \\ Similarly, $\sum_h a_h t_{h,l-1} =\sum_{h\; \mbox{\tiny even}} a_h t_{h,l-1} =0$ since $l-1+h$ is odd for even $h$.
\end{enumerate}
\qed \\

\subsection{Generalized Jordan blocks}\label{section genJNF}

In this section we recall some well-known facts about normal forms for endomorphisms. Let $m, l \in \mathbb{N}$ be natural numbers, and consider $ml\times ml$-matrices $A$ with  $l\times l$-blocks $A_{ij}$. Let $N$ be the standard nilpotent matrix of size $ml$. Then $N^l$ has blocks $I_l$ as first upper off-diagonal of blocks, and all other entries zero, \text{i.e.}
\[
 N^l=\left(\begin{array}{cccc} O&I_l&O\\O&O&I_l\\&\ddots&\ddots&\ddots
 \end{array}\right).
\]
Note that for any invertible matrix $K \in \fieldk^{l \times l}$, one has $(KI) N^l (K^{-1}I)=N^l$, where
\[
KI=\left(\begin{array}{cccc} 
K&O&O\\O&K&O\\&\ddots&\ddots&\ddots
 \end{array}\right).
\]

\begin{prop}\label{fnf}
Let $\eta$ be an indecomposable endomorphism of the $n$-dimensional vector space $V$ over a perfect field $\fieldk$. Then 
\begin{enumerate}
\item There is unique $m \in \mathbb{N}$ and irreducible monic $q(x) \in \fieldk[x]$  such that $q(x)^m$ is the minimal (= characteristic) polynomial of $\eta$; that is, $q(x)^m$ is the unique elementary divisor of $\eta$. In particular, $ml=n$ where  $l=\deg q(x)$. Moreover, as an $\fieldk[x]$-module, $V$ is isomorphic to $\fieldk[x]/q(x)^m$; in particular it is cyclic.
\item\label{gen Jordan} Let $q(x)^m$ be the unique monic elementary divisor of $\eta$. There is a basis $\bar v$ of $V$
with respect to which  $\eta$ has matrix 
\[
A=ZI+N^l=
\left(\begin{array}{cccc} 
Z&I_l&O&\cdots\\
O&Z&I_l&\ddots\\
O&O&Z&\ddots\\
\vdots&&\ddots&\ddots 
\end{array}\right).
\]
where $q(x)=\det(xI_i-Z)$.
Conversely, for any basis with respect to which the coordinate matrix of $\eta$ is of the form $A=ZI+N^l$, one has $\det(xI_l-Z)=q(x)$ and $\det(xI-A)=q(x)^m$.
\end{enumerate}
\end{prop}

\pf 
1. This follows from the theory of a linear operator on a finite-dimensional vector space. In general, the $k[x]$-module $V$ is isomorphic to the direct sum of the $\fieldk[x]/d_i(x)$ where the $d_i(x)=q_i(x)^{m_i}$ are the elementary divisors of $\eta$, the $q_i(x)$ being coprime irreducibles. The characteristic polynomial is $\prod_i d_i(x) =\det(xI-A)$. That $\eta$ is indecomposable means that there is only a single elementary divisor $q(x)^m$, and so the characteristic polynomial and minimal polynomial of $A$ coincide with $q(x)^m$.  

2. For the existence of such a basis and corresponding normal form we use our assumption that $\fieldk$ is perfect; see \cite{Malcev:foundations} for a discussion of this normal form, and \cite{Robinson:genJordan} and \cite{Dalalyan:genJordan} regarding necessary and sufficient conditions for its existence. Conversely, given such a normal form for $\eta$, it follows from the behaviour of block matrices and determinants that $q(x)^m = \det (xI - A) = (\det(xI_l - Z))^m$. By the uniqueness of factorization of monic polynomials into irreducible polynomials, and by comparing degrees, it follows that $\det(xI_l - Z) = q(x)$. 
\qed \\

\begin{rmk}
We refer to the normal form given in part \ref{gen Jordan} above as ``generalized Jordan normal form''. The matrix $Z$ is, in coordinate form, a ``generalized eigenvalue'' of $\eta$. 
\end{rmk}

\subsection{Admissible forms for self-duals}

\begin{thm}\label{hom}
Given an indecomposable  sextuple $\mathcal{S}_\eta$ with no eigenvalue in $\fieldk$, the following statements are equivalent: 
\begin{enumerate}
\item it is self-dual
\item it admits compatible symplectic forms
\item it admits compatible symmetric forms.
\end{enumerate}
In the case of self-duality, compatible $\varepsilon$-symmetric forms can be given, with respect to suitable bases, by matrices of the form  
\[ 
H_M=\left(
\begin{array}{ccc} 
O&O&-M\\
O&M&O\\
-M&O&O 
\end{array}\right),
\quad \ \text{ with } \ M=\left(
\begin{array}{cccc} 
\vdots &\vdots&&\\
O&O&T&\\
O&-T&O&\cdots\\
T&O&O&\cdots 
\end{array}\right)
\in \fieldk^{ml \times ml},
\] 
where $T\in \fieldk^{l\times l}$ and $q(x)^m$ is the characteristic polynomial of $\eta$, with $q(x)$ irreducible and   $\deg q(x)=l$. Clearly,  $\varepsilon(H_M)=\varepsilon(M)=(-1)^{m+1}\varepsilon(T)$. Suitable bases are frame bases $\{ \bar u, \bar v, \bar w \} $ where the basis $\bar v = v_1,\ldots ,v_{ml}$ renders $\eta$ in normal form $(\frac{1}{2}I_l+\tilde Z)I+N^l$  with  Frobenius matrix $\tilde Z \in \fieldk^{l\times l}$. Here, $T$ can be chosen according to Proposition \ref{mc}. 
\end{thm}

\pf
Given compatible forms, self-duality is obvious. Conversely, consider self-dual indecomposable  $\mathcal{S}_\eta$ and consider a basis $\bar v$ such that the corresponding matrix of $\eta$ is $A=ZI+N^l$ as in 2. of Proposition \ref{fnf}. Let $q(x)^m$ be the unique elementary divisor of $\eta$, where $q$ is irreducible over $\fieldk$ with $\deg q =l$. 

By Proposition \ref{duals of endo-sextuples}, $A$ is similar to $(I-A)^t$. Thus, $q(x)^m$ is also the unique elementary divisor of $(I-A)^t$, implying in particular that $\det(xI_l- (I_l-Z)^t) = q(x)$ (\text{c.f.} the proof of part 2. in Proposition \ref{fnf}). So $Z$ and $(I_l-Z)^t$ share $q(x)$ as their unique elementary divisor, and are hence similar. Now by Proposition \ref{sim}, $\tilde Z=Z-\frac{1}{2}I_l$  has irreducible characteristic polynomial $r(x)$ with $r(-\tilde Z)=0$. 

Since $r$ is irreducible, there exists invertible $K$ such that $K\tilde ZK^{-1}$ is Frobenius. Note that $KZK^{-1}= \frac{1}{2}I_l+K\tilde Z K^{-1}$ and $KAK^{-1}=  \frac{1}{2}I +K\tilde ZK^{-1}I +N^l$. Thus we may assume directly that $\eta$ has matrix $A=ZI+N^l$ with $Z=\frac{1}{2}I_l+\tilde Z$ where $\tilde Z$ is a Frobenius matrix with irreducible characteristic polynomial $r(x)$ satisfying $r(-\tilde Z)=0$. Moreover,  $TZT^{-1}=(I_l-Z)^t$ if and only if $T\tilde ZT^{-1}=-\tilde Z^t$ for any invertible $T$ (c.f. Proposition \ref{sim}). By Proposition \ref{mc}, there exist $T$, skew-symmetric as well as symmetric, which satisfy $T\tilde ZT^{-1}=-\tilde Z^t$. Choose such a $T$ and let $M$ be defined as above in the statement of this proposition. By Proposition \ref{duals of endo-sextuples}, it suffices to show that $MAM^{-1}=(I-A)^t$.  

For this, observe first that $M=TQJ=QTJ$ where  $J$ is block-diagonal with diagonal blocks  $J_{jj}=(-1)^{j+1}I_l$ and $Q$ is block-anti-diagonal with blocks $Q_{ij}=I_l$ for $i+j=m+1$, and all other entries being $0$. Also, observe that  $Q^2=I$ and that  $Q$ acts as permutation of blocks: acting on the right it exchanges block-columns $j$ and $m+1-j$, acting on the left it exchanges block-rows $i$ and $m+1-i$. 
Furthermore, the following properties are easily seen:
\begin{itemize}
\item $QDQ=D$ if $D$ is block-diagonal of the form $XI$ for some $X \in \fieldk^{l \times l}$
\item $QN^lQ=(N^l)^t$
\item $J^2=I$, and $JDJ=D$ if $D$ is block-diagonal 
\item $JN^lJ=-N^l$
\item $T(N^l)^tT^{-1}=(N^l)^t$ since $(N^l)^t$ has only unit and zero blocks. 
\end{itemize}

%

It follows that
\[ 
MZIM^{-1} =QTJ ZIJT^{-1}Q=QT ZIT^{-1}Q= 
Q(I_l-Z)^tIQ = (I_l-Z)^tI
\]
\[ 
MN^lM^{-1} =TQJ N^lJQT^{-1}= TQ(-N^l)QT^{-1}= -T(N^l)^tT^{-1}= -(N^l)^t
\]
and hence
\[ 
MAM^{-1}=MZIM^{-1}+ MN^lM^{-1} = 
(I_l-Z^t)I -(N^l)^t = I-A^t, 
\]
as desired. 
\qed \\

\begin{cor}\label{paramet underlying}
Up to isomorphism, the indecomposable sextuples underlying non-split continuous-type isotropic triples are paremetrized by indecomposable linear endomorphisms $\eta$ of the form 
$$
\eta = \tfrac{1}{2} + \zeta
$$
where $\zeta$ is an indecomposable endomorphism which lives in an even-dimensional space and  is such that $\zeta$ is similar to $-\zeta^*$. 
Such endomorphisms $\zeta$ are themselves parametrized by the set
$$
\{ r^m \mid m \in \mathbb{Z}_{>0}, \ r \text{ monic irreducible and } r = p(x^2) \text{ with } p \in \fieldk[x] \} \cup \{ (x^2)^m \mid m \in \mathbb{Z}_{>0} \}.
$$
\end{cor}

\pf  
From Proposition \ref{duals of endo-sextuples} we know that (isomorphism classes of) indecomposable self-dual continuous-type sextuples are in one-to-one correspondence with (isomorphism classes of) indecomposable endomorphisms $\eta$ such that $\eta$ is similar to $1 - \eta^*$.  Setting $\zeta := \tfrac{1}{2} - \eta$ it easily seen that $\zeta$ is indecomposable if and only if $\eta$ is, and that 
$$
\eta \text{ similar to } 1- \eta^* \quad \Leftrightarrow \quad \zeta \text{ similar to } - \zeta^*.
$$
Moreover, Theorem \ref{hom} and Theorem \ref{eig} imply that the underlying sextuples of non-split isotropic triples, up to isomorophism, are parametrized by all endomorphisms $\eta$ as above, with the exception of those which have an eigenvalue in the ground field and live in an odd-dimensional space. 

In the case when $\eta$ has no eigenvalue in the ground field, it follows from the above proof of Theorem \ref{hom} and from Lemma \ref{extensions and involution} that the corresponding endomorphism $\zeta$ has a minimal polynomial of the form $r(x)^m$, with $r$ irreducible and with only even powers of the variable appearing, \text{i.e.} $r(x) = p(x^2)$ for some $p \in \fieldk[x]$. (The irreducible polynomial $r$ here is the minimal polynomial of the ``eigenvalue'' $\tilde Z$ of $\zeta$ which appears in the proof of Theorem \ref{hom}; for Lemma \ref{extensions and involution} we use $\tilde Z$ in the role of $\mu$). 

In the case when $\eta$ does have an eigenvalue in the ground field, the corresponding endomorphism $\zeta$ is nilpotent, and hence, if it lives in an even-dimensional space, it has a minimal polynomial of the form $r(x) = (x^2)^{m}$ for some $m \in \mathbb{Z}_{>0}$ 
\qed

\subsection{Endomorphism algebras}\label{endo alg details}

In order to address the question of uniqueness of compatible forms, we first recall the structure of endomorphism
algebras of indecomposable sextuples $\mathcal{S}_\eta$. We've seen that these are isomorphic to endomorphism algebras of vector spaces with indecomposable endomorphism, \text{i.e.} endomorphism algebras of $\fieldk [x]$-modules $V_{\fieldk[\eta]}$, where $\eta$ is an indecomposable endomorphism of a vector space $V$ over $\fieldk$. These are well-known; we recall some basic facts.

\begin{prop}\label{endo alg}
Given $\eta,Z,A,\bar v$ as in Proposition \ref{fnf},  the following hold:
\begin{enumerate}
\item The endomorphism algebra $E$ 
of $V_{\fieldk[\eta]}$ is isomorphic to the $\fieldk$-algebra $\fieldk[x]/q(x)^m$,
which is of dimension $n=ml = \dim_\fieldk V$, where $l$ is the degree of the irreducible minimal polynomial $q(x)$ of $\eta$. In particular, $E$ is local with radical
$E\cdot q(x)$.
\item With respect to $\bar v$, $E$ is given by the matrices $C$ such that $CA=AC$. We have $C\in E$ if and only if  $C=\sum_{i=0}^{k-1}Z_iN^{li}$ with $Z_i\in \fieldk(Z)$, and this representation is unique.  In particular, $C$ is upper block triangular with diagonal blocks $C_{11}= \ldots = C_{kk} = Z_0$ and  $C$ is invertible in $E$ if and only if $Z_0\neq 0$.
\end{enumerate} 
\end{prop} 

\pf 
1. Rephrasing the definitions, $E$ is the collection of $\fieldk$-linear endomorphisms of the vector space $V$ which commute with $\eta$. For $\eta$ indecomposable, $E$ is simply $\fieldk[\eta] \subseteq \text{End}(V)$, and hence is isomorphic to $\fieldk[x]/q(x)^m$. One way to see this is to note that since $\eta$ is indecomposable, there exists a basis of $V$ of the form $\{v, \eta v, ..., \eta^{n-1}v \}$ for some $v \in V$. For any $f \in E$, observe that $fv = \sum_{i=0}^{n-1} c_i\eta^iv$ for some $c_i \in \fieldk$;  it follows then that in fact $f = \sum_{i=0}^{n-1} c_i\eta^i$ since the latter expression coincides with $f$ on the above given basis of $V$.

2. The first statement is clear. To prove the first ``if'' it suffices to observe that, for $W \in \fieldk(Z)$, one has $WZ=ZW$ since $\fieldk(Z)$ is commutative and $WIN^l=N^lWI$ since $N^l$ has zero and unit
blocks only. Thus, matrices $C=\sum_{i=0}^{k-1}Z_iN^{li}$ form a subalgebra $E'$ of $E$.  

To prove the ``only if'', we show that $E'$ and $E$ have the same dimension over $\fieldk$. For this, we compute the dimension of $E'$ as an $\fieldk(Z)$-algebra. 
Since $N^{lm}=0$, induction on $m-h>0$ shows that $I, N^{l}, N^{l2}, \ldots ,N^{l(m-1)}$ are independent over $\fieldk(Z)$: if $\sum_{j=h}^{m-1} Z_jN^{lj}=0$ then, applying $N^l$, we get $\sum_{j=h}^{m-2} Z_jN^{l(j+1)}=0$ and hence, using the induction hypothesis, that $Z_h = ... = Z_{m-1} = 0$.  Thus, $E'$ has dimension $m$ over $\fieldk(Z)$ and so dimension $ml=n=\dim E$ over $\fieldk$, since $\fieldk(Z)$ has dimension $l$ over $\fieldk$. Thus, $E'=E$.
\qed \\

\subsection{Uniqueness: General result}\label{uniq:gen}

In this section we generalize Lemma \ref{uniqueness up to a scalar} in a way which applies to all poset representations and also to endomorphisms $(V, \eta)$.
Recall (see Section \ref{decomp linear reps}) that for the latter there is a natural notion of direct sum, morphism, etc., and that for an indecomposable endomorphism $(V, \eta)$, its endomorphism algebra
$$
\text{End}(V,\eta) = \{ f : V \rightarrow V \mid f\eta = \eta f \}
$$ 
is always local. 
For our application, we will use the following notion of the \textbf{dual} of an endomorphism: we define
$$
(V, \eta)^* := (V^*, \text{Id}- \eta^*). 
$$
Analogous to the definition for poset representation, a \text{compatible form} for $(U, \eta)$ is an isomorphism $B: (V, \eta) \rightarrow (V^*, \text{Id}-\eta^*)$ which defines either to a symmetric or skew-symmetric bilinear form on $V$. 


Within the current subsection, the word ``representation'' and the notation $``\psi "$ will be used to denote a poset representation or an endomorphism (unless further specification is given) and similarly for compatible forms, etc.. Given a bilinear form $B$ we write $\varepsilon(B)=1$ if $B$
is symmetric, and $\varepsilon(B)=-1$ if $B$ is skew-symmetric. Similar notation also applies for matrices,  we set $A^{-t}=(A^{-1})^t$, and we use $A^*=A^t$ interchangeably.

Consider now the following situation. Assume that we are given an indecomposable  representation $\psi$ in $V$ of dimension $ml$, a basis of $V$, and a subfield $F$ of $\fieldk^{l\times l}$  such that the endomorphism algebra $E$
of $\psi$  is given by block matrices (with $l\times l$-blocks) of the form
\[ 
\sum_{i=0}^{m-1} Z_iN^{li} 
\]  
with unique  $Z_i \in F$ and such that $ZN^l=N^lZ$ for all $Z\in F$.
In particular, $E$ is commutative and has radical ${\sf rad}E=EN^l =\sum_{i=1}^{m-1} X_iN^{li}$.

Further, assume that we are given a compatible form for $\psi$, which we call $B_1$. Denote the corresponding coordinate matrix by $H_1$, which is itself necessarily also a block matrix as above.  

\begin{dfn}
Let $\psi$, $F$, $H_1$ be as above. For $A \in \fieldk^{l\times l}$, set $A^\dagger := H_1^{-1}A H_1$. The anti-involution ``$(-)^\dagger$'' restricts to one on $F$ via $Z^\dagger := (ZI)^\dagger$. We define
$$
F_{H_1}^+ := \{Z \in F \mid Z^\dagger = Z \} \quad \quad \quad \quad F_{H_1}^-  := \{ Z \in F \mid Z^\dagger = - Z \}.
$$
Note that 
$$
F_{H_1}^+ = \{Z \in F \mid Z^tH_1=H_1Z\} \quad  \quad \quad \quad F_{H_1}^-  = \{ Z \in F \mid Z^tH_1=-H_1Z\}.
$$
\end{dfn}
\begin{rmk}
$F = F_{H_1}^+ \oplus F_{H_1}^-$ by the usual argument: any $Z \in F$ has a unique decomposition 
$$
Z= \tfrac{1}{2}(Z + Z^\dagger) + \tfrac{1}{2}(Z - Z^\dagger)
$$
into selfadjoint and anti-selfadjoint parts. 
\end{rmk}
\begin{rmk}
If $0\neq Z \in F_{H_1}^+ $, then $H := H_1Z$ determines a compatible form $B$ 
such that $\varepsilon(B)=\varepsilon(B_1)$, since in this case
$$
(H_1Z)^t= Z^tH_1^t = \varepsilon_1  Z^t  H_1= \varepsilon_1 H_1Z.
$$
Similarly, given $0\neq Z' \in F_{H_1}^- $, then $H' := H_1Z$ defines an compatible form $B'$ 
such that $\varepsilon(B')=- \varepsilon(B_1)$.

\end{rmk}

\begin{lemma}\label{uniqu} 
Consider, in addition to $B_1$, another compatible form $B_2$ for $\psi$, with associated coordinate matrix $H_2$.
Set $\varepsilon_1 := \varepsilon(B_1)$, $\varepsilon_2 := \varepsilon(B_2)$ and $\epsilon = \varepsilon_1 \varepsilon_2$. 

If $\epsilon=1$, there exists $0\neq Z\in F^+_{H_1} $ and an automorphism $f$  of $\psi$ which is an isometry from $B_2$ to the compatible form given by $H_1Z$. 
\end{lemma}

\noindent \textbf{Notation}: In this case we say that $B_2$ and $B_1$ are \emph{equivalent up to  automorphisms of $\psi$ and multiplication with ``scalars'' in $F$}. 

\

\pf 
Let $^\dagger$ denote the antiautomorphism  given by the operation of adjoint with respect to $H_1$, i.e. 
$A^\dagger = H _1^{-1} A^t H_1$. Note that when $A$ is in the endomorphism algebra $E$ of $\psi$, then so is $A^\dagger$. Note also that $(H_1^{-1})^\dagger = H_1^{-t}$. 

Observe that $H_1^{-1} H_2$ determines an automorphism of $\psi$, so
\[  H_1^{-1}H_2=C_0I-R_0  \mbox{ for some invertible } C_0 \in F \mbox{ and }
R_0\in {\sf rad}E. 
\]
It follows that
\begin{align*}
(C_0I)^\dagger -R_0^\dagger =& \ (H_1^{-1} H_2)^\dagger = H_2^\dagger (H_1^{-1})^\dagger=H_1^{-1} H_2^t H_1 H_1^{-t} \\
& = H_1^{-1}\varepsilon_2 H_2H_1 \varepsilon_1 H_1^{-1} = \epsilon H_1^{-1} H_2 = \epsilon C_0I- \epsilon R_0.
\end{align*}
Since $E = FI  \oplus \rad E$ and this decomposition is preserved under taking adjoints,  
$ (C_0I)^\dagger=\epsilon C_0I$ and $R_0^\dagger = \epsilon R_0$; in particular $C_0\in F^\epsilon_{H_1} $. 

Let $\epsilon = 1$. Since $A \mapsto A^\dagger$ is an anti-automorphism of $E$, $(C_0^{-1}I)^\dagger= \epsilon C_0^{-1}I = C_0^{-1}I$. Set
\[
R :=C_0^{-1}R_0=R_0C_0^{-1},\quad  H_3 := H_2C_0^{-1},\quad C := H_1^{-1}H_3 = I-R 
\] 
Then $R^\dagger = (C_0^{-1})^\dagger R_0^\dagger = C_0^{-1}R_0 = R$ and since $R$ is nilpotent we can proceed as in Lemma \ref{uniqueness up to a scalar} and construct a unit $h \in E$ such that  $h^* H_1 h = H_3$ (where $H_3$ here plays the role of $H_2$ in that Lemma). Setting $f :=h^{-1}$ and using that $E$ is commutative, we obtain
\[
f^* H_2f = f^* H_3C_0f = f^* H_3f C_0  =H_1C_0. 
\]
\qed \\ 

\begin{lemma}\label{uniqu2}
Let $B_1$ be as above, and assume that its matrix $H_1$ is zero above the $l\times l$-block anti-diagonal. Let $B_0$ be another compatible form, with matrix $H_1Z$ for some non-zero $Z \in F$. 

There is an automorphism $f$ of $\psi$ which is an isometry from $B_0$ to $B_1$  if and only if there exist $X\in F^+_{H_1}$ and $Y\in F^-_{H_1}$ such that  $Z=X^2-Y^2$. 
\end{lemma}

\pf
If $Z=X^2-Y^2$ with $X \in F^+_{H_1}$ and $Y \in F^-_{H_1}$ then
define $f$ by $(X+Y)I$. In this case, 
\begin{align*}
f^*H_1Zf =& \ (X+Y)^tH_1Z(X+Y)= (X^t+Y^t)H_1(X+Y)Z \\
	&=(X^t+Y^t)(X^t-Y^t)H_1Z = Z^tH_1Z = H_1,
\end{align*}
using in the last step that $Z = X^2 - Y^2$ implies that $Z \in F^+_{H_1}$. 

Conversely, let $f \in \text{End} (\psi)$ be an isometry from $B_0$ to $B_1$. The matrix $A$ of $f$ is of the form $A= \sum_i Z_i N^i$ with $Z_i \in F$  and $A^tH_1A=H_1Z$. Observe that $N$ commutes with the $Z_i$ and that applying $N$ on the right of a matrix moves all columns of blocks one step to the right (with overspill) and sets the first column to zero. Similarly, application of $N^t$ on the left shifts block rows downward. Thus, since $H_1$ is zero above the  block anti-diagonal,
 $(N^t)^j H_1 N^i$ has zero block anti-diagonal if $i+j>0$. Now
\[
H_1Z=A^tH_1A= \sum_{i,j}Z^t_j(N^t)^jH_1N^iZ_i 
\]
and $H_1Z$ is zero above the block anti-diagonal, and on the block anti-diagonal the blocks are
\[
T_{ij}Z=Z_0^tT_{ij}Z_0,\ \text{ with  }i+j=m+1
\]
where the $T_{ij}$ denote the blocks of $H_1$.
By hypothesis, $Z_0=X+Y$ for some $X\in F^+_{H_1}$ and $ Y\in F^-_{H_1}$. It follows that
\[ 
Z= T_{ij}^{-1} (X^t+Y^t)T_{ij}(X+Y) =
(X-Y)(X+Y)=X^2-Y^2. 
\] 
\qed \\

\subsection{Uniqueness for sextuples $\mathcal{S}_\eta$}\label{unihom}


 For indecomposable self-dual framed sextuples $\mathcal{S}_\eta$ we obtain the following uniqueness result for compatible forms by applying Lemma \ref{uniqu} and Lemma \ref{uniqu2} to the underlying endomorphism $(U, \eta)$. 
minimal polynomial 
\begin{thm}\label{uniqueness general fields}
Let $\mathcal{S}_\eta$ be an indecomposable self-dual framed sextuple such that $\eta$ has no eigenvalue in $\fieldk$,  and let $q(x)^m$ be the minimal polynomial of $\eta$, with $q(x)$ irreducible, $\deg q(x)=l$. 

Let $B_1$ and $B_2$ be compatible (skew)symmetric forms for $\mathcal{S}_\eta$. Set $\varepsilon_1 := \varepsilon(B_1)$, $\varepsilon_2 := \varepsilon(B_2)$, and $\epsilon := \varepsilon_1 \varepsilon_2$. Furthermore: 
\begin{itemize}
\item choose a basis of $U$ which is as in Proposition \ref{fnf} and extend this to a frame basis for $\mathcal{S}_\eta$;
\item define a subfield $F \subseteq \fieldk^{l \times l}$ by $F=\fieldk(Z)$, with $Z$ from Proposition \ref{fnf}.
\end{itemize}
In terms of the respective coordinate matrices $H_1$ and $H_2$ of the compatible forms $B_1$ and $B_2$, we have:
\begin{enumerate}
\item\label{unique pt1} If $\epsilon = 1$, then $H_1$ and $H_2$ are equivalent up to automorphisms of $\mathcal{S}_\eta$ and multiplication with ``scalars'' in $F^+_{H_1}$. 
\item\label{unique pt2}  Given $0 \neq K \in F$, there is an automorphism $f$ of $\mathcal{S}_\eta$ which is an isometry from $H_1$ to $H_1K$  if and only if there are $X\in F^+_{H_1}$ and $Y\in F^-_{H_1}$ such that  $K=X^2-Y^2$. 
\end{enumerate}
\end{thm}

\pf 
In view of Proposition \ref{fnf} and Propositions \ref{endo alg}, we can apply Lemma \ref{uniqu} to the underlying endomorphism $(U, \eta)$. It follows then from Proposition \ref{endo to sextuple} and Proposition \ref{duals of endo-sextuples} that we can transfer the uniqueness statement from Lemma \ref{uniqu} to the corresponding statement in part (\ref{unique pt1}) above for compatible forms for $\mathcal{S}_\eta$. 

We turn to proving part \ref{unique pt2}. If $K=X^2-Y^2$ for some $X\in F^+_{H_1}$ and $Y\in F^-_{H_1}$, then an isometry is given by $(X+Y)\text{Id}$ as in the proof of Lemma \ref{uniqu2}. So assume that $H_1$ is equivalent (isometric) to $H_1K$. 

Without loss of generality we can assume that $H_1$ is given by one of the canonical compatible forms defined in Theorem \ref{hom}. Indeed, suppose that the statement to be proved were true for those canonical forms (and let $H_1$ represent, for a moment, a form which is not necessarily such a canonical one).  Suppose moreover that there exists an isometry $f$ between $H_1$ and $H_1K$. Choose the appropriate canonical compatible form $H_0$ such that $\varepsilon(H_0)= \varepsilon(H_1) = \varepsilon(H_1K)$. Then, by part \text{2.} above, there exists some $C \in F^+_{H_1}$ and an isometry $g : H_0 \rightarrow H_1C$. This will also be an isometry $g: H_0K \rightarrow H_1CK = H_1KC$. From all this we obtain the isometry
$$
g^{-1}f g : H_0 \rightarrow H_0K. 
$$
By assumption, this implies that $K = X^2 - Y^2$ for some $X\in F^+_{H_1}$ and $Y\in F^-_{H_1}$. 

So we can assume $H_1$ is a canonical compatible form. Now we wish to proceed in an analogous manner as we did in proving parts 1 and 2, but this time using Lemma \ref{uniqu2}. It remains only to check that the hypotheses of Lemma \ref{uniqu2} are satisfied by $(U, \eta)$. 

Note that from Theorem \ref{hom} and Proposition \ref{endo to sextuple} it follows that the coordinate matrix $H_1 \vert_U$ is of the form

\[ 
\left(
\begin{array}{cccc} 
\vdots &\vdots&&\\
O&O&T&\\
O&-T&O&\cdots\\
T&O&O&\cdots 
\end{array}\right)
\in \fieldk^{ml \times ml},
\] 
where $T\in \fieldk^{l\times l}$ is chosen as in the proof of Theorem \ref{hom}. In particular, $H_1$ is zero above the $l \times l$ block anti-diagonal, as required by Lemma \ref{uniqu2}. 
%
\qed \\

\end{document}